\documentstyle[11pt]{article}
  
\begin{titlepage}

\title{The Spectrum of Kleinian Manifolds}

\author{Ulrich Bunke\thanks{Mathematisches Institut, Universit\"at G\"ottingen, Bunsenstr. 3-5, 37073 G\"ottingen, GERMANY, E-mail:bunke@cfgauss.uni-math.gwdg.de} and
Martin
Olbrich\thanks{Mathematisches Institut, Universit\"at G\"ottingen, Bunsenstr. 3-5, 37073 G\"ottingen, GERMANY, E-mail:bunke@cfgauss.uni-math.gwdg.de  }
}

\end{titlepage}

\newcommand{\proof}{{\it Proof.$\:\:\:\:$}}

\newcommand{\cM}{{\cal M}}

\newcommand{\dist}{{\rm dist}}
\newcommand{\kaaa}{{\bf k}}
\newcommand{\paaa}{{\bf p}}
\newcommand{\taaa}{{\bf t}}
\newcommand{\haaa}{{\bf h}}

\newcommand{\R}{{\bf R}}
\newcommand{\Q}{{\bf Q}}
\newcommand{\Z}{{\bf Z}}
\newcommand{\C}{{\bf C}}

\newcommand{\gaaa}{{\bf g}}
\newcommand{\maaa}{{\bf m}}
\newcommand{\aaaa}{{\bf a}}
\newcommand{\naaa}{{\bf n}}

\newcommand{\res}{{\rm res}}

\newcommand{\cZ}{{\cal Z}}
\newcommand{\cH}{{\cal H}}
\newcommand{\interi}{{\rm int}}

\newcommand{\cE}{{\cal E}}

\newcommand{\cD}{{\cal D}}

\newcommand{\cA}{{\cal A}}

\newcommand{\cU}{{\cal U}}
\newcommand{\Hom}{{\mbox{\rm Hom}}}

\newcommand{\End}{{\mbox{\rm End}}}

\newcommand{\im}{{\mbox{\rm im}}}

\newcommand{\cF}{{\cal F}}
\newcommand{\Ree}{{\rm Re }}

\newcommand{\Imm}{{\rm Im}}

\newcommand{\clo}{{\rm clo}}

\newcommand{\Mat}{{\rm Mat}}
\newcommand{\ee}{{\rm e}}

\newcommand{\tr}{{\mbox{\rm tr}}}

\newcommand{\Ad}{{\mbox{\rm Ad}}}
\newcommand{\ad}{{\mbox{\rm ad}}}

\newcommand{\id}{{\mbox{\rm id}}}

\newcommand{\nat}{{\bf  N}}
\newcommand{\supp}{{\mbox{\rm supp}}}

\newcommand{\aca}{{\aaaa_\C^\ast}}
\newcommand{\acag}{{\aca(\sigma)}}

\newcommand{\cR}{{\cal R}}
\def\hB{\hspace*{\fill}$\Box$ \newline\noindent}

\def\hB{\hspace*{\fill}$\Box$ \\[0.5cm]\noindent}

 \newcommand{\cG}{{\cal G}}

\newtheorem{prop}{Proposition}[section]
\newtheorem{lem}[prop]{Lemma}
\newtheorem{ddd}[prop]{Definition}
\newtheorem{theorem}[prop]{Theorem}
\newtheorem{kor}[prop]{Corollary}
\newtheorem{ass}[prop]{Assumption}

\newtheorem{fact}[prop]{Fact}

\begin{document}

\maketitle

\tableofcontents
\section{Introduction}

This paper is devoted to the spectral decomposition of the space
of sections $L^2(Y,V_Y(\gamma))$ of a locally homogenous bundle $V_Y(\gamma)$ over a locally symmetric
space $Y$ of rank one of infinite volume with respect to locally invariant
differential operators.

Let $G$ be a real semisimple linear connected Lie group of real rank one
with finite center. Let $K\subset G$ be a maximal compact subgroup.
Then $X:=G/K$ is a Riemannian symmetric space of negative curvature.

Let $\Gamma\subset G$ be a convex-cocompact, torsion-free, discrete
subgroup and $Y:=\Gamma\backslash X$ be the corresponding locally symmetric space. 
Let $(\gamma,V_\gamma)$ be a finite-dimensional unitary representation
of $K$. Then we form the (locally)
homogeneous vector bundle $V(\gamma):=G\times_KV_\gamma$ ($V_Y(\gamma):=\Gamma\backslash G\times_KV_\gamma$).

Let $\gaaa$ denote the Lie algebra of $G$,
$\cU(\gaaa)$ the universal enveloping algebra of $\gaaa$ and $\cZ$ its center.
Through the left regular action of $\cU(\gaaa)$ on $C^\infty(X,V(\gamma))$ 
any  $A\in \cZ$ gives rise to an $G$-invariant differential
operator $A_\gamma$. This operator descends to $C^\infty(Y,V_Y(\gamma))$.

Let $\cZ_\gamma$ denote the algebra $\{A_\gamma\:|\: A\in \cZ\}$.
By $\nabla$ we denote the canonical invariant connection
of $V(\gamma)$. 
The algebra $\cZ_\gamma$ is a finite extension of the algebra $\C(\Delta)$, where  $\Delta=\nabla^*\nabla$ is the
Bochner Laplace operator on $V(\gamma)$.
We employ a suitable invariant scalar product on $\gaaa$ in order to normalize the Riemannian metric of $X$ and the Casimir operators
$\Omega_G$ and $\Omega_K$ of $G$ and $K$. The Casimir operators are related with the Laplacian by $\Delta=-\Omega_G+\gamma(\Omega_K)$.

We view $\cZ_\gamma$ as an algebra of unbounded operators 
on the Hilbert spaces $L^2(X,V(\gamma))$ (resp. $L^2(Y,V_Y(\gamma))$) with common domain $C_c^\infty(X,V(\gamma))$ (resp. $C_c^\infty(Y,V_Y(\gamma))$). 
Since $X$ (resp. $Y$) are complete any formally selfadjoint locally invariant elliptic operator is essentially selfadjoint on the domain $C^\infty_c(X,V(\gamma))$ (resp. $C_c^\infty(Y,V_Y(\gamma))$).
Note that $\Omega_G$ is elliptic and formally selfadjoint.
Using this it is easy to see that the algebra $\cZ_\gamma$ can be generated
by essentially selfadjoint elements with commuting resolvents.
Thus there exist spectral decompositions
of $L^2(X,V(\gamma))$ (resp. $L^2(Y,V_Y(\gamma))$) with respect to $\cZ_\gamma$
(which can equivalently be considered as a spectral decompositions with respect to $\cZ$ as we will do).
We will also consider spectral decompositions of these Hilbert spaces
with respect to finite abelian extensions of $\cZ_\gamma$ which are obtained
by adjoining further essentially selfadjoint differential operators.

Since $\cZ_\gamma$ is a quotient of $\cZ$ we can parametrize
characters of $\cZ_\gamma$ using the Harish-Chandra isomorphism.
Let $\haaa\subset \gaaa$ denote a Cartan algebra of $\gaaa$,
$W=W(\gaaa,\haaa)$ the Weyl group, and $\haaa_\C^*$ the complexified 
dual of $\haaa$.
The Harish-Chandra isomorphism identifies characters of
$\cZ$ with points in $\haaa^*_\C/W$. Let $\lambda\in\haaa^*_\C$
represent some $W$-orbit. Then we denote the corresponding character
of $\cZ$ by $\chi_\lambda$.

The abstract spectral decomposition gives a measurable
field of Hilbert spaces $\{H_\lambda\}_{\lambda\in \haaa^*_\C/W}$,
a measure $\kappa$ on $\haaa^*_\C/W$, and an isometry
$$\alpha:L^2(Y,V_Y(\gamma))\cong \int_{\haaa^*_\C/W} H_\lambda \kappa(d\lambda)\ .$$  
If we let $\cZ$ act on $H_\lambda$ by the character $\chi_\lambda$,
then $\alpha$ is compatible with the action of $\cZ$.
It is of course apriori known that $\kappa$ is supported
on the set of $\lambda$ with the property that $\chi_\lambda$
factors through $\cZ_\gamma$. 

We want to describe this set in greater detail.
The sphere bundle of $X$ can be identified with the homogeneous space
 $G/M$, where $M\subset K$. 
Let $\maaa$ denote the Lie algebra of $M$. 
We choose a Cartan algebra $\taaa$ of $\maaa$ and let $\aaaa$ be a one-dimensional subspace of the orthogonal complement of $\kaaa$ in $\gaaa$.
Then $\aaaa\oplus\taaa=:\haaa$ is a Cartan algebra of $\gaaa$.
We further choose a positive root system of $\taaa$.
Let $\rho_m$ denote half of the sum of the positive roots of $(\maaa,\taaa)$.
For $\sigma\in \hat{M}$ let $\mu_\sigma\in\taaa^*$ be its highest weight.
Let $\aca$ denote the complexification of the dual of $\aaaa$.
A pair $(\sigma\in \hat{M},\lambda\in\aca)$ determines the character
$\chi_{\mu_\sigma+\rho_m-\lambda}$ of $\cZ$.
Then representation theory of $G$ implies that
$$\supp(\kappa)\subset \{\chi_{\mu_\sigma+\rho_m-\lambda}|[\gamma_{|M}:\sigma]\not=0,\lambda\in\aca\}\ .$$ 
In fact there are more restrictions since  
characters contributing to the spectral decomposition must be selfadjoint.
This restriction implies that if
$\chi_{\mu_\sigma+\rho_m-\lambda}\in\supp(\kappa)$, then $\lambda$ has to be either real or
imaginary. Thus we apriori know that the support of $\kappa$ is contained
in the projection to $\haaa^*_\C/W$ of a finite union of lines in $\haaa^*_\C$.
Speaking about the absolute-continuous part of the spectrum, we have in mind that the measure $\kappa$ restricted to the corresponding one-dimensional
set is absolute continuous with respect to the one-dimensional Lebesgue measure.

We do not employ this apriory knowledge about the support of $\kappa$ in the proofs, it rather follows from our arguments. 

The goal of this paper is to describe 
in detail spectral decomposition of $L^2(Y,V_Y(\gamma))$ with respect to 
the algebra $\cZ$.
The Eisenstein series is used to 
identify the absolute continuous
spectrum. We then show that $\kappa$ has no singular continuous
component. Finally obtain the finiteness of the point spectrum
and a description of all infinite-dimensional eigenspaces.
In analogy to the spectral decomposition in the finite volume case
we obtain a decomposition 
$$L^2(Y,V_Y(\gamma))=L^2(Y,V_Y(\gamma))_c\oplus L^2(Y,V_Y(\gamma))_{res}\oplus L^2(Y,V_Y(\gamma))_{cusp}\oplus L^2(Y,V_Y(\gamma))_{scat} .$$
Here $L^2(Y,V_Y(\gamma))_c$ is the continuous part given by wave packets
of Eisenstein series. The space  $L^2(Y,V_Y(\gamma))_{res}$ is the finite-dimensional  residual part which is essentially 
generated by the residues of Eisenstein series. The cuspidal
part $L^2(Y,V_Y(\gamma))_{cusp}$ consists of a finite number of
infinite dimensional eigenspaces and is related to the discrete
series representations of $G$ occuring in $L^2(X,V(\gamma))$.
The scattering part $L^2(Y,V_Y(\gamma))_{scat}$ consists
of finite-dimensional eigenspaces at the boundary of the continuous spectrum.
In contrast to the finite-volume case this part can be non-trivial
as we demonstrate by an example.
  
The motivation for studying the spectral decomposition with respect to $\cZ$
(and larger commutative algebras) 
instead of $\Delta$ is that these algebras encode additional symmetries.
If only the spectral decomposition of $L^2(Y,V_Y(\gamma))$ with respect
to the Laplacian $\Delta$ is considered, then one encounters
embedded eigenvalues. Their "stability"
is explained by the additional symmetries since
they are isolated with respect to the larger algebras. 

For locally symmetric manifolds of the sort considered in present paper
the spectral decomposition of $L^2(Y,V_Y(\gamma))$ with respect to the Laplacian
(respectively partial results) were obtained by 
\begin{itemize} 
\item Patterson \cite{patterson75} for trivial $\gamma$ and surfaces
\item Lax-Phillips \cite{laxphillips82}, 
\cite{laxphillips841}, \cite{laxphillips842}, \cite{laxphillips85}, Perry \cite{perry87} for higher dimensional hyperbolic manifolds and trivial $\gamma$
\item Mazzeo-Phillips
\cite{mazzeophillips90} for differential forms on hyperbolic manifolds  
\item Epstein-Melrose-Mendoza \cite{epsteinmelrosemendoza91}, Epstein-Melrose \cite{epsteinmelrose90} for differential forms on complex-hyperbolic manifolds.
\end{itemize}
There is related work on Eisenstein series and the scattering matrix in the real  hyperbolic case for trivial $\gamma$ (e.g.
\cite{patterson76}, \cite{patterson761}, \cite{patterson89}, \cite{mandouvalos86}, \cite{mandouvalos89}, \cite{perry89}).

At the end of this introduction let us make some remarks concerning the methods.
Once the Eisenstein series are constructed the
realization of the absolute continuous part of the spectrum is almost standard.
One of the most important steps in proving the spectral decomposition is
to show the absence of the singular continuous spectrum.
Usually, the limiting absorption principle (e.g. \cite{perry87}) or commutator methods (see e.g. \cite{froesehislopperry91}) are employed at this point. Here we use  a completely different method (proposed in \cite{bernstein88}) which is based on an apriori knowledge of all relevant generalized eigenfunctions of $\cZ$. Our discussion of the point spectrum is based on the asymtotic expansion
of eigenfunctions and boundary value theory. 

Before starting with the main topic of the paper in Section \ref{esess}
we analyse the boundary values of generalized eigenfunctions along
the geodesic boundary $\partial X$ of $X$. In particular we are interested in the
space of $\Gamma$-invariant distributional sections of homogeneous
bundles over $\partial X$ with support in the limit set of $\Gamma$.
Sections \ref{fiert} to \ref{invvv} are devoted to the analysis
on $\partial X$.

Our results also have a more representation theoretic interpretation.
Let $L^2(\Gamma\backslash G)_K$ denote the space of all $K$-finite
vectors on $L^2(\Gamma\backslash G)$. Then combining our results for
all $\gamma\in K$ one can obtain a decomposition of 
$L^2(\Gamma\backslash G)_K$ into unitarizable $(\gaaa,K)$-modules
and a decomposition of $L^2(\Gamma\backslash G)$ as a direct
integral of unitary representations of $G$. The classification of the unitary dual of $G$ then leads to further restrictions of the location of the residual part of the spectrum.

\noindent
{\it Acknowledgement: We thank R. Mazzeo and P. Perry for discussing
of parts of this work.}

\section{Geometric preparations}\label{fiert}

Let $G$ be a connected, linear, real semisimple Lie group of rank one, $G=KAN$ be an Iwasawa decomposition
of $G$, $\gaaa=\kaaa\oplus\aaaa\oplus\naaa$ be the corresponding Iwasawa decomposition of the Lie algebra $\gaaa$,
$M:=Z_K(A)$ be the centralizer of $A$ in $K$ and $P:=MAN$
be a minimal parabolic subgroup. 
The group $G$ acts isometrically on the rank-one symmetric space $X:=G/K$.
Let $\partial X:=G/P=K/M$ be its geodesic boundary. We consider $X\cup\partial X$ as a compact manifold with boundary.

By the classification of symmetric spaces with strictly negative
sectional curvature $X$ is one of the following spaces:
\begin{itemize}
\item a real hyperbolic space,
\item a complex hyperbolic space,
\item a quaternionic hyperbolic space,
\item or the Cayley hyperbolic plane, 
\end{itemize}
and $G$ is a linear group finitely covering of the orientation-preserving isometry 
group of $X$. 

Let $\Gamma \subset G$ be a torsion-free, discrete subgroup.
\begin{ass}\label{asss}
We assume that there is a $\Gamma$-invariant partition $\partial X =\Omega\cup \Lambda$, where  $\Omega\not=\emptyset$ is open and 
$\Gamma$ acts freely and cocompactly on 
$X\cup\Omega$.
\end{ass}
The locally symmetric space $Y:=\Gamma\backslash X$ is a complete Riemannian manifold of infinite volume without cusps. It can be compactified by adjoining
the geodesic boundary $B:=\Gamma\backslash \Omega$.
We call $\Lambda$ the limit set of $\Gamma$.

A group $\Gamma$ satisfying \ref{asss} is also called convex-cocompact
since it acts cocompactly on the convex hull of the limit set $\Lambda$.
The quotient $Y$ can be called a Kleinian manifold in generalizing
the corresponding notion for three-dimensional hyperbolic manifolds.

We now consider some geometric consequences of \ref{asss}
which eventually allow us to define the exponent $\delta_\Gamma$ of $\Gamma$.

Let $g=\kappa(g)a(g)n(g)$, $\kappa(g)\in K$, $a(g)\in A$, $n(g)\in N$ be
defined with respect to the given Iwasawa decomposition.
By $\aca$ we denote the comlexified dual of $\aaaa$.
If $\lambda\in \aca$, then we set $a^\lambda:=\ee^{\langle \lambda,\log(a)\rangle}\in\C$.
The roots of $\aaaa$ corresponding to $\naaa$ distinguish a positive cone $\aaaa^*_+$.
Define $\rho\in \aaaa_+^*$ as usual by $\rho(H):=\frac{1}{2}\tr(\ad(H)_{|\naaa})$, $\forall H\in\aaaa$.

We adopt the following conventions about the notation for points of $X$ and $\partial X$.
A point $x\in \partial X$ can equivalently be denoted by a subset $kM\subset K$
or $gP\subset G$ representing this point in $\partial X=K/M$ or $\partial X=G/P$.
If $F\subset \partial X$, then $FM:=\bigcup_{kM\in F}kM\subset K$.
Analogously, we can denote $b\in X$ by $gK\subset G$, where $gK$ represents
$b$ in $X=G/K$.
 
\begin{lem}\label{no}
For any compact $F\subset \Omega$ we have $\sharp(\Gamma\cap FMA_+K)<\infty$.
\end{lem}
\proof
The compact set $FMA_+K\cup F\subset X\cup\Omega$ contains at most a finite
number of points of the orbit $\Gamma K$ of the origin of $X$. \hB

Let $A_+:=\exp(\aaaa_+)$.
Any element $g\in G$ has a decomposition $g=ka_gh$, $k,h\in K$, $a_g\in A_+\cup \{1\}$, where $a_g$ is uniquely determined by $g$.

\begin{lem}\label{no1}
Let $k_0M \in \partial X$. For any compact $W\subset (\partial X\setminus k_0M)M$ there exists a neighbourhood $U\subset K$ of $k_0M$ and constants $c>0$, $C<\infty$, such that for all $g=ha_gh^\prime \in WA_+K$ and $k\in U$  
\begin{equation}\label{unm}
 c a_g \le  a(g^{-1}k) \le C a_g \ .
 \end{equation}
\end{lem}
\proof 
The set $W^{-1}k_0M$ is compact and disjoint from $M$.
Let $w\in N_K(M)$ represent the non-trivial element of the Weyl group of $(\gaaa,\aaaa)$.
Set   $\bar{\naaa}=\theta(\naaa)$,
where $\theta$ is the Cartan involution of $G$ fixing $K$ and define $\bar{N}:=exp(\bar{\naaa})$.
There is a precompact open $V\subset \bar{N}$
such that $W^{-1}k_0M\subset w \kappa(V)M$.
By enlarging $V$  we can  assume that $V$ is $A_+$-invariant, where $A$ acts on $\bar{N}$
by  $(a,\bar{n})\mapsto a\bar{n}a^{-1}$.
Moreover,  there exists an open neighbourhood $U\subset K$ of $k_0M$
such that $w^{-1}W^{-1}U M\subset \kappa(V)M$.

Let $k\in U$ and $g=ha_gh^\prime \in WA_+K$.
Then we have $h^{-1}k=w\kappa(\bar{n})m $ for $\bar{n}\in V$, $m\in M$.
Furthermore,
\begin{eqnarray*}
a(g^{-1}k)&=&a(h^{\prime -1}a_g^{-1}h^{-1}k)\\
&=& a(a_g^{-1}w\kappa(\bar{n})m) \\
&=&a(a_g \kappa(\bar{n}))\\
&=&a(a_g\bar{n}n(\bar{n})^{-1}a(\bar{n})^{-1})\\
&=&a(a_g\bar{n}a_g^{-1}) a(\bar{n})^{-1} a_g\ .
\end{eqnarray*}
Now $a_g\bar{n}a_g^{-1}\in V$. Set 
\begin{eqnarray*}
c&:=&\inf_{\bar{n}\in V} a(\bar{n}) \inf_{\bar{n}\in V} a(\bar{n})^{-1}\\
C&:=&\sup_{\bar{n}\in V} a(\bar{n}) \sup_{\bar{n}\in V} a(\bar{n})^{-1}\ . 
\end{eqnarray*}
Since $V$ is precompact we have $0 < c \le C<\infty$.
It follows that $c a_g \le a(g^{-1}k )\le Ca_g\ .$
\hB
   
\begin{lem}\label{poi} 
The series 
$$\sum_{g\in\Gamma} a_g^{-2\rho}$$ 
converges.
\end{lem}
\proof
We represent $\partial X=G/P$. There is a $P$-invariant splitting $\gaaa=\bar{\naaa}\oplus\paaa$
and an identification $\bar{\naaa}^*=\naaa$ by the invariant bilinear form on $\gaaa$.  In particular $\Lambda^{max}\bar{\naaa}^*=\Lambda^{max}\naaa$ as a 
$P$-module. 
It follows that \begin{equation}\label{sad}\Lambda^{max}T^*\partial X=G\times_P \Lambda^{max}\naaa \ . \end{equation}
This bundle is trivial as a $K$-homogeneous bundle.  
Fix an orientation of $\partial X$ and let
$\omega\in C^\infty(\partial X, \Lambda^{max}T^*\partial X)$
be a positive $K$-invariant volume form,
which is unique up to a positive scalar factor.
 
Let $g\in G$. Then $g$ acts as a diffeomorphism on $\partial X$.
Using (\ref{sad}) and the $K$-invariance of $\omega$ we obtain $(g^*\omega)(kM)=a(g^{-1}k)^{-2\rho}\omega(kM)$.
Let $F\subset \Omega$ be a compact set with non-trivial interior such that
$gF\cap F=\emptyset$ for all $1\not=g\in\Gamma$. Such $F$ exists by Assumption \ref{asss}.
Then 
\begin{eqnarray*}
\infty&>& \int_{\partial X} \omega 
>  \int_{\cup_{g\in\Gamma}\: gF} \omega\\
&=& \sum_{g\in\Gamma} \int_{gF} \omega 
 = \sum_{g\in\Gamma} \int_F g^*\omega\\
&=& \sum_{g\in\Gamma} \int_F a(g^{-1}k)^{-2\rho} \omega(kM).
\end{eqnarray*}
Let $F_1\subset \Omega$ be a compact neighbourhood of $F$.
By Lemma \ref{no} the set $\Gamma\cap F_1MA_+K$ is finite.
We apply the Lemma \ref{no1} taking for $W$ the closure of 
$\Gamma\setminus (\Gamma\cap F_1MA_+K)$. Then we can cover $FM$
with finitly many sets $U$ the existence of which was asserted
in that lemma. 
Thus here is a constant $C\in A$ such that for all $g\in \Gamma\setminus (\Gamma\cap F_1MA_+K)$ and $k\in FM$
$$a(g^{-1}k)\le Ca_g\ .$$
It follows that 
\begin{eqnarray*}
\sum_{g\in\Gamma}  a_g^{-2\rho}  \int_F \omega&=& \sum_{g\in\Gamma} \int_F a_g^{-2\rho} \omega\\&
\le&  C^{2\rho} \sum_{g\in\Gamma\setminus(\Gamma\cap F_1 MA_+K)} \int_F a(g^{-1}kM)^{-2\rho} \omega(kM)\ \  + \sum_{g\in\Gamma\cap F_1MA_+K}   a_g^{-2\rho}\\[0.5cm]
&<& \infty\ . 
\end{eqnarray*}
This implies the lemma since
$$  \int_F  \omega\not=0\ .$$
\hB

\begin{ddd}
Let $\delta_\Gamma\in \aaaa^*$ be the smallest element such that $\sum_{g\in\Gamma} a_g^{- \rho-\lambda}$ converges for all $\lambda\in\aaaa^*$ with $\lambda>\delta_\Gamma$.
\end{ddd}
$\delta_\Gamma$ is called the exponent of
$\Gamma$. By Proposition \ref{poi} we have $\delta_\Gamma \le \rho$.
If $\Gamma$ is non-trivial, then $\delta_\Gamma\ge -\rho$.
For the trivial group we have $\delta_{\{1\}}=-\infty$.

If $X$ is the $n$-dimensional hyperbolic space and we identify $\aaaa^*$ with $\R$ such that $\rho=\frac{n-1}{2}$, then it was shown by Patterson \cite{patterson762} and Sullivan \cite{sullivan79}, that $\delta_\Gamma+\frac{n-1}{2}=\dim_H(\Lambda)$,
where $\dim_H$ denotes the Hausdorff dimension
(the Hausdorff dimension of the empty set is
by definition $-\infty$). It was also shown in \cite{sullivan79} that in the real hyperbolic case $\sum_{g\in\Gamma} a_g^{- \rho-\delta_\Gamma}$ diverges.
Hence $\delta_\Gamma<\rho$.

In the proof of the meromorphic continuation of the Eisenstein
series we employ at a certain place 
that $X$ belongs to a series of symmetric
spaces with increasing dimensions.
It is there where we need the following assumption
\begin{ass}\label{caly}
If $X$ is the Cayley hyperbolic plane, then we assume that $\delta_\Gamma<0$.
\end{ass}
We believe that this assumption is only of technical nature
and can be dropped using other methods for the meromorphic continuation
of the Eisenstein series.

\section{Analytic preparations}\label{anaprep}

The goal of the following two sections is to construct 
$\Gamma$-invariant vectors in principal series representations of $G$.
For this reason we introduce the extension map $ext$ and the
scattering matrix, and construct their meromorphic continuations.
The principal series representations of $G$ are realized
on spaces of distribution sections of bundles over $\partial X$.
Roughly speaking the extension map extends a $\Gamma$-invariant
distribution sections on $\Omega$ across the limit set $\Lambda$.
The space of $\Gamma$-invariant distributions on $\Omega$
is easily described as the space of distribution sections
on bundles over $B$. The majority of the $\Gamma$-invariant vectors
of the principal series is then obtained by extension.

Consider a finite-dimensional unitary representation $\sigma$ of $M$  on $V_\sigma$. If $w\in W(\gaaa,\aaaa)$ is the non-trivial element
of the Weyl group, then we can form 
the representation 
$\sigma^w$ of $M$ on $V_\sigma$ by conjugating
the argument of $\sigma$ with a representative of $w$ in $N_K(M)$.
In this section we will assume that $\sigma$ is irreducible.
Note that we will change this convention later in the case that $\sigma$ and $\sigma^w$ are non-equivalent by considering the sum $\sigma\oplus \sigma^w$ instead.

For $\lambda\in \aca$  we  form the representation $\sigma_\lambda$
of $P$ on $V_{\sigma_\lambda}:=V_\sigma$, which is given by
$\sigma_\lambda(man):=\sigma(m)a^{\rho-\lambda}$.
Let $V(\sigma_\lambda):=G\times_P V_{\sigma_\lambda}$ be the associated
homogeneous bundle.
Set $V_B(\sigma_\lambda):=\Gamma\backslash V(\sigma_\lambda)$.

Let $\tilde{\sigma}$ be the dual representation to $\sigma$.
Then there are  natural pairings 
\begin{eqnarray*}
V(\tilde{\sigma}_{-\lambda})\otimes V(\sigma_\lambda)&\rightarrow& \Lambda^{max}T^*\partial X\\
V_B(\tilde{\sigma}_{-\lambda})\otimes V_B(\sigma_\lambda)&\rightarrow& \Lambda^{max}T^*\partial B\ .
\end{eqnarray*}
The orientation of $\partial X$ induces one of $B$.
Employing these pairings and integration with respect to the fixed
orientation we obtain identifications
\begin{eqnarray*}
C^{-\infty}(\partial X,V(\sigma_\lambda))&=&C^{\infty}(\partial X,V(\tilde{\sigma}_{-\lambda}))^\prime\\
C^{-\infty}(B,V_B(\sigma_\lambda))&=&C^{\infty}(B,V_B(\tilde{\sigma}_{-\lambda}))^\prime\ .
\end{eqnarray*}

As a $K$-homogeneous bundle we have a canonical identification $V(\sigma_\lambda)\cong K\times_M V_\sigma$. Thus $\bigcup_{\lambda\in\aca} V(\sigma_\lambda)\rightarrow \aca\times \partial X$
has the structure of a trivial holomorphic family of bundles.

Let $\pi^{\sigma,\lambda}$ denote the representation of $G$ on the space of sections of $V(\sigma_\lambda)$ given by the left-regular representation. 
Then $\pi^{\sigma,\lambda}$ is called a principal series representation
of $G$. Note that there are different globalizations of this
representation which are distinguished by the regularity
of the sections (smooth, distribution  e.t.c.).

For any small open subset $U\subset B$
and diffeomorphic lift $\tilde{U}\subset \Omega$  the restriction $V_B(\sigma_\lambda)_{|U}$
is canonically isomorphic to $V(\sigma_\lambda)_{|\tilde{U}}$.
Let $\{U_\alpha\}$ be a cover of $B$ by open sets as above.
Then  
$$\bigcup_{\lambda\in\aca} V_B(\sigma_\lambda)\rightarrow \aca\times B$$
can be given the structure of a holomorphic family of bundles by glueing the trivial families
$$\bigcup_{\lambda\in\aca}V_B(\sigma_\lambda)_{|U}\cong \bigcup_{\lambda\in\aca}V(\sigma_\lambda)_{|\tilde{U}}$$ together
using the holomorphic families of glueing maps induced by $\pi^{\sigma,\lambda}(g)$, $g\in\Gamma$.
Thus it makes sense to speak of holomorphic or smooth or continuous families
of sections
$\aca\ni\mu\mapsto f_\mu\in C^{\pm\infty}(B,V_B(\sigma_\mu))$.

When dealing with holomorphic families of vectors in topological vector
spaces we will employ the following functional analytic facts.
Let $\cF,\cG,\cH \dots$ be complete locally convex
topological vector spaces.
A locally convex vector space is called a Montel space if its   
closed bounded subsets are compact. 
A Montel space is reflexive, i.e., the canonical map into its bidual is an isomorphism.
Moreover, the dual space of a Montel space is again a Montel space.
\begin{fact} The space of smooth sections of a vector bundle
and its topological dual are Montel spaces.
\end{fact}
We equip $\Hom(\cF,\cG)$ with the topology of uniform convergence
on bounded sets.
  Let $V\subset \C$ be open.
A map $f:V\rightarrow \Hom(\cF,\cG)$ is called holomorphic
if for any $z_0\in V$ there is a sequence $f_i\in \Hom(\cF,\cG)$ such that
 $f(z)=\sum_{n=0}^\infty f_i (z-z_0)^i$ converges for all $z$ close to $z_0$.
Let $f:V\setminus \{z_0\} \rightarrow \Hom(\cF,\cG)$ be holomorphic
and  $f(z)=\sum_{n=-N}^\infty f_i (z-z_0)^i$ for all  $z\not=z_0$ close to $z_0$.
Then we say that  $f$ is meromorphic and has a pole of order $N$ at $z_0$.
If $f_i$, $i=-N,\dots,-1$, are finite dimensional, then $f$
has, by definition, a finite-dimensional singularity.
We call a subset  $A\subset \cF\times \cG^\prime$ sufficient if for    
$B\in \Hom(\cF,\cG)$ the condition 
 $<\phi,B \psi >=0$, $ \forall (\psi,\phi)\in A$,  implies  $B=0$.
\begin{fact}\label{holla}
The following assertions are equivalent :
\begin{enumerate}
\item (i)  $\:\:f:V\rightarrow \Hom(\cF,\cG)$ is holomorphic.
\item (ii) $\:\:\: f$ is continuous and 
there is a sufficient set $A\subset \cF\times \cG^\prime$
such that for all $(\psi,\phi)\in A$ the function $V\ni z\mapsto \langle \phi,f(z)\psi\rangle $
is holomorphic.
\end{enumerate}
\end{fact}
 
\begin{fact}\label{seq}
Let $f_i:V\rightarrow Hom(\cF,\cG)$ be a sequence
of holomorphic maps. Moreover let $f :V\rightarrow Hom(\cF,\cG)$ be continuous such that for a sufficient set
$A \subset \cF\times \cG^\prime$ the functions $\langle \phi,f_i \psi\rangle $,
$(\psi,\phi)\in A$,    
converge locally uniformly
in $V$ to $\langle\phi,f \psi\rangle$. Then $f$ is holomorphic, too.
\end{fact}
 
\begin{fact}\label{adjk}
Let $f:V\rightarrow \Hom(\cF,\cG)$
 be continuous.
Then the adjoint $f^\prime:V\rightarrow \Hom(\cG^\prime,\cF^\prime)$
is continuous.
If $f$ is holomorphic, then so is $f^\prime$. 
\end{fact}

\begin{fact}\label{comp}
Assume that $\cF$ is a Montel space.
Let $f:V\rightarrow \Hom(\cF,\cG)$
and $f_1:V\rightarrow \Hom(\cG,\cH)$
be continuous.
Then $f_1\circ f : V\rightarrow \Hom(\cF,\cH)$
is continuous.
If $f,f_1$ are holomorphic, so is $f_1\circ f$. 
\end{fact}
 The following lemma will be employed in Section \ref{extttttt}.
Since it is of purely functional analytic nature we consider it at this place.

Let $\cH$ be a Hilbert space and $\cF\subset \cH$
be a Fr\'echet space such that the embedding is continuous and compact.
In the application we have in mind $\cH$ will be some $L^2$-
space of sections of a vector bundle over a compact closed manifold
and $\cF$ be the Fr\'echet space of smooth sections of this bundle.
The continuous maps $\Hom(\cH,\cF)$ will  be called smoothing operators.

Let $V\subset \C$ be open and connected, and $V\ni z\rightarrow R(z)\in \Hom(\cH,\cF)$
be a meromorphic family of smoothing operators  with at most finite-dimensional singularities.
Note that $R(z)$ is a meromorphic family of compact operators on $\cH$ in a natural way.

\begin{lem} \label{merofred}
If $1-R(z)$ is invertible at some point $z\in V$ where $R(z)$ is regular,
then 
$$(1-R(z))^{-1}=1-S(z)\ ,$$ 
where $V\ni z\rightarrow S(z)\in \Hom(\cH,\cF)$
is a meromorphic family of smoothing operators  with at most finite dimensional singularities.
\end{lem}
\proof
We apply Reed-Simon IV \cite{reedsimon78}, Theorem XIII.13
in order to conclude that $(1-R(z))^{-1}$ is a meromorphic
family of operators on $\cH$ having at most finite-dimensional singularities.
Making the ansatz $(1-R(z))^{-1}= 1-S(z)$, where apriori $S(z)$ is a meromorphic
familiy of bounded operators on $\cH$ with finite dimensional singularities,
we obtain $S=-R-R\circ S$. This shows that $S$ is a meromorpic family in $\Hom(\cH,\cF)$.
\hB
 
This finishes the functional analytic preparations and we now turn to the construction of the extension map. In fact, we first introduce its adjoint
which is the push-down
$$\pi_\ast:C^\infty(\partial X,V(\sigma_\lambda))\rightarrow  C^\infty(B,V_B(\sigma_\lambda))\ .$$ 
 
Using the identification  $C^\infty(B,V_B(\sigma_\lambda))= {}^\Gamma C^\infty(\Omega,V(\sigma_\lambda))$
we define $\pi_\ast$ by
\begin{equation}\label{summ}\pi_*(f)(kM)= \sum_{g\in\Gamma} (\pi(g)f)(kM),\quad kM\in\Omega \ ,\end{equation}
if the sum converges. Here $\pi(g)$ is the action
induced by $\pi^{\sigma,\lambda}(g)$.

Note that the universal enveloping algebra 
$\cU(\gaaa)$ is a filtered algebra. Let $\cU(\gaaa)_m$, $m\in\nat_0$,
be the space of elements of degree less or equal than $m$.
For any $m$ and bounded subset $A\subset \cU(\gaaa)_m$
we define the seminorm $\rho_{m,A}$ on $C^\infty(\partial X,V(\sigma_\lambda))$
by $$\rho_{m,A}(f):=\sup_{X\in A, k\in K}|f(\kappa(kX))|\ .$$
These seminorms define the Fr\'echet topology of $C^\infty(\partial X,V(\sigma_\lambda))$
(in fact a countable set of such seminorms is sufficient).

In order to define the Fr\'echet topology on $C^\infty(B,V_B(\sigma_\lambda))$
we fix an open cover $\{U_\alpha\}$ of $B$ such that each $U_\alpha$ has a diffeomorphic
lift $\tilde{U}_\alpha \subset\Omega$.
Then we have canonical isomorphisms 
$$C^\infty(\tilde{U}_\alpha ,V(\sigma_\lambda))\cong C^\infty(\ U_\alpha,V_B(\sigma_\lambda))\ .$$
For any $U\in \{U_\alpha\}$ we define the topology of $C^\infty(\tilde{U},V(\sigma_\lambda))$
using the seminorms 
$$\rho_{U,m,A}(f):=\sup_{X\in A, k\in \tilde{U}M }|f(\kappa(kX))|\ ,$$
where  $m\in\nat_0$ and $A\subset \cU(\gaaa)_m$ is bounded.
Since $C^\infty(B,V_B(\sigma_\lambda))$ maps to $C^\infty( U_\alpha,V_B(\sigma_\lambda))$
by restriction for each $\alpha$ we obtain a system of seminorms defining the Fr\'echet topology of $C^\infty(B,V_B(\sigma_\lambda))$.

\begin{lem}\label{anal1}
If $\Ree(\lambda)<-\delta_\Gamma$, then the sum (\ref{summ}) converges for $f\in C^\infty(\partial X,V(\sigma_\lambda))$ and defines
a  continuous map $$\pi_\ast: C^\infty(\partial X,V(\sigma_\lambda))\rightarrow  C^\infty(B,V_B(\sigma_\lambda))\ .$$
Moreover, $\pi_*$ depends holomorphically on $\lambda$.
\end{lem}
\proof 
Consider $U\in \{U_\alpha\}$. We want to estimate  
$$C^\infty(\partial X,V(\sigma_\lambda))\ni f\mapsto res_{\tilde{U}}\circ \pi(g) f \in 
C^\infty(\tilde{U},V(\sigma_\lambda))\ .$$
Let $\Delta:\cU(\gaaa)\rightarrow \cU(\gaaa)\otimes \cU(\gaaa)$
be the coproduct and write $\Delta(X)=\sum_i X_i\otimes Y_i$.
Fix $l\in \nat_0$ and a bounded set $A\in \cU(\gaaa)_l$.
Then there is another bounded set $A_1\subset \cU(\gaaa)_l$
depending on $A$ such that
\begin{eqnarray*}
\rho_{U,l,A}(res_{\tilde{U}M}\circ \pi(g) f)&=&\sup_{X\in A,k\in \tilde{U}}|(\pi(g)f)(\kappa(kX))|\\
&=& \sup_{X\in A,k\in \tilde{U}M}|\sum_ i a(g^{-1}\kappa(kX_i))^{\lambda-\rho}f(\kappa(g^{-1}\kappa(kY_i)))|\\
&\le & 
\sup_{X\in A_1,k\in \tilde{U}M}| a(g^{-1}kX)^{\lambda-\rho}| \sup_{X\in A_1,k\in \tilde{U}M}|   f(\kappa(g^{-1}kX))| \ .
\end{eqnarray*}

The Poincar\'e-Birkhoff-Witt theorem gives a  decomposition $\cU(\gaaa)=\cU(\bar{\naaa})\cU(\maaa)\cU(\aaaa)\oplus \cU(\gaaa)\naaa$.
Let $q:\cU(\gaaa)\rightarrow \cU(\bar{\naaa})\cU(\maaa)\cU(\aaaa)$
be the associated projection. Then for $g\in G$ and $X\in\cU(\gaaa)$ we have 
$\kappa(gX)=\kappa(gq(X))$, $a(gX)=a(gq(X))$.

Let $U_1\subset \Omega$ be an open neighbourhood of $\tilde{U}$.
Then by Lemma \ref{no} the intersection $\Gamma\cap U_1MA_+K$ is finite.
Let $W:=(\partial X\setminus U_1)M$. Then by Lemma \ref{no1}
we can find a compact $A_+$-invariant set $V\subset \bar{\naaa}$
such that $W^{-1}\tilde{U}M\subset w\kappa(V)M$.
For $g=ha_gh^\prime\in WA_+K$ and $k\in \tilde{U}M$ we obtain $h^{-1}k=w\kappa(\bar{n})m$ for some $\bar{n}\in V$, $m\in M$.

Let $X\in \cU(\gaaa)$. Then
\begin{eqnarray*}
\kappa(g^{-1}kX)&=&\kappa(h^{\prime -1}a_g^{-1}h^{-1}kX)\\
&=&h^{\prime -1}\kappa(a_g^{-1}w\kappa(\bar{n})mX)\\
&=&h^{\prime -1}w\kappa(a_g\bar{n}n(\bar{n})^{-1}a(\bar{n})^{-1}mX) \\
&=&h^{\prime -1}w\kappa(a_g\bar{n}a_g^{-1} a_g[n(\bar{n})^{-1}a(\bar{n})^{-1}mXm^{-1}a(\bar{n})n(\bar{n})]a_g^{-1})m\\
&=&h^{\prime -1}w\kappa(a_g\bar{n}a_g^{-1} a_gq(n(\bar{n})^{-1}a(\bar{n})^{-1}mXm^{-1}a(\bar{n})n(\bar{n}))a_g^{-1})m\ .
\end{eqnarray*}
Since $V$ is compact the sets  $n(V)^{-1}a(V)^{-1}MA_1Ma(V)n(V)=:A_2\subset \cU(\gaaa)_l$
and  $q(A_2)$ are bounded. 
Conjugating $q(A_2)$ with $A_+$ gives
clearly another bounded set $A_3\subset \cU(\gaaa)_l$.
We can find a bounded set $A_4\subset \cU(\gaaa)_l$ such that
$\kappa(a_g\bar{n}a_g^{-1}A_3)\subset \kappa(\kappa(a_g\bar{n}a_g^{-1})A_4)$
for all $a_g\in A_+$.
This implies for $g\in WA_+K$ that
\begin{equation}\label{kkll1}\sup_{X\in A_1,k\in \tilde{U}}|   f(\kappa(g^{-1}kX))|\le \rho_{l,A_4}(f)\ .\end{equation}

We also have 
\begin{eqnarray*}
a(g^{-1}kX)&=&a(h^{\prime -1}a_g^{-1}h^{-1}kX)\\
&=& a(a_g^{-1}w\kappa(\bar{n})mX)\\
&=&a(a_g \kappa(\bar{n})mXm^{-1})\\
&=&a(a_g\bar{n}n(\bar{n})^{-1}a(\bar{n})^{-1}mXm^{-1})\\
&=&a(a_g\bar{n}a_g^{-1}a_g n(\bar{n})^{-1}a(\bar{n})^{-1}mXm^{-1}a(\bar{n})n(\bar{n})a_g^{-1}) a(\bar{n})^{-1} a_g\ .\\
&=&a(a_g\bar{n}a_g^{-1}a_g q(n(\bar{n})^{-1}a(\bar{n})^{-1}mXm^{-1}a(\bar{n})n(\bar{n}))a_g^{-1}) a(\bar{n})^{-1} a_g\ . 
\end{eqnarray*}
Again there is a constant $C<\infty$ such that 
$$ |a(a_g\bar{n}a_g^{-1}a_g q(n(\bar{n})^{-1}a(\bar{n})^{-1}mXm^{-1}a(\bar{n})n(\bar{n}))a_g^{-1})^{\lambda-\rho}| a(\bar{n})^{ \rho- \lambda } <C$$
for all $a_g\in A_+$, $\bar{n}\in V$, $m\in M$, and $X\in A_1$.
It follows that
\begin{equation}\label{hj1}\sup_{X\in A_1,k\in \tilde{U}M}| a(g^{-1}kX)^{\lambda-\rho}|\le C a_g^{\lambda-\rho}\end{equation}
for almost all $g\in \Gamma$.
The estimates (\ref{kkll1}) and (\ref{hj1}) together imply that the sum
$$C^l(\partial X,V(\sigma_\lambda))\ni f\mapsto \sum_{g\in\Gamma} res_{\tilde{U}}\circ \pi(g) f \in  C^l(\tilde{U},V(\sigma_\lambda))$$
converges for $\Ree(\lambda)<-\delta_\Gamma$ 
and defines a continuous map of Banach spaces.
This map depends holomorphically on $\lambda$
by Fact \ref{seq}.

Combining  these considerations for all  
$U\in \{U_\alpha\}$ and $l\in\nat_0$ we obtain that 
$$\pi_*:C^\infty(\partial X,V(\sigma_\lambda))\rightarrow C^\infty(B,V_B(\sigma_\lambda))$$
is defined and continuous for  $\Ree(\lambda)<-\delta_\Gamma$. 
Moreover it is easy to see that $\pi_*$ depends holomorphically
on $\lambda$. \hB

Still postponing the introduction of the extension map we
now consider its left-inverse, the restriction 
$$res:{}^\Gamma C^{-\infty}(\partial X,V(\sigma_\lambda))\rightarrow C^{-\infty}(B,V_B(\sigma_\lambda))\ .$$
In fact the space ${}^\Gamma C^{-\infty}(\Omega,V(\sigma_\lambda))$ 
of $\Gamma$-invariant vectors in $C^{-\infty}(\Omega,V(\sigma_\lambda))$
can be canonically identified with the corresponding space 
$C^{-\infty}(B,V_B(\sigma_\lambda))$.
Composing this identification with the restriction
$res_\Omega:C^{-\infty}(\partial X,V(\sigma_\lambda))\rightarrow C^{-\infty}(\Omega,V(\sigma_\lambda))$
we obtain the required restriction map $res$. 
 
\begin{lem}\label{lopi}
There exists a continous map 
$$\widetilde{res}: C^{-\infty}(\partial X ,V(\sigma_\lambda))\rightarrow C^{-\infty}(B,V_B(\sigma_\lambda))\ ,$$  
which depends holomorphically on $\lambda$ and coincides with $res$ on
${}^\Gamma C^{-\infty}(\Omega,V(\sigma_\lambda))$.
\end{lem}
\proof
We exhibit $\widetilde{res}$ as the adjoint of a continuous map
$$\pi^*:C^\infty(B,V_B(\tilde{\sigma}_{-\lambda}))\rightarrow C^\infty (\partial X ,V(\tilde{\sigma}_{-\lambda}))$$
which depends holomorphically on $\lambda$. Then the lemma follows from Fact
\ref{adjk}.

Let $\{U_\alpha\}$ be a finite open cover of $B$  such that each $U_\alpha$ has a diffeomorphic
lift $\tilde{U}_\alpha \subset\Omega$. Choose a subordinated partition of unity
$\{\chi_\alpha\}$. Pulling $\chi_\alpha$ back to $\tilde{U}_\alpha$ and extending the resulting function by $0$ we obtain a function $\tilde \chi_\alpha\in  C^\infty (\partial X ,V(\tilde{\sigma}_{-\lambda}))$. We define 
$$\pi^*(f):=\sum_\alpha \tilde\chi_\alpha f, \quad f\in C^\infty(B,V_B(\tilde{\sigma}_{-\lambda}))\ ,$$
where we consider $f$ as an element of 
${}^\Gamma C^{-\infty}(\partial X,V(\tilde{\sigma}_{-\lambda}))$.
Then we set $\widetilde{res}:=(\pi^*)^\prime$. \hB

The extension map $ext$ will be defined as an right inverse to $res$.
\begin{ddd}\label{defofext}
For $\Ree(\lambda)>\delta_\Gamma$ we define 
the extension map 
$$ext: C^{-\infty}(B,V_B(\sigma_\lambda))\rightarrow {}^\Gamma C^{-\infty}(\partial X,V(\sigma_\lambda))$$
to be the adjoint of 
$$\pi_*:C^{ *}(\partial X,V(\tilde{\sigma}_{-\lambda})) \rightarrow C^{*}(B,V_B(\tilde{\sigma}_{-\lambda}))\ .$$
\end{ddd}
This definition needs a justification. In fact, 
by Lemma \ref{anal1} the extension exists, is continuous, and  by Fact \ref{adjk} it  depends holomorphically on $\lambda$. Moreover, it is easy to see that the range of the adjoint of $\pi_*$ consists of $\Gamma$-invariant vectors.
\begin{lem}\label{iden}
We have 
$res\circ ext = \id$. \end{lem}
\proof
Recall the definition of $\pi^*$ from the proof of Lemma \ref{lopi}. Then
$\pi_*\pi^*$ is the identity on $C^\infty(B,V_B(\tilde{\sigma}_{-\lambda}))$.
We obtain
$$res \circ ext=\widetilde{res}\circ ext=  (\pi^*)^\prime \circ (\pi_*)^\prime
=(\pi_*\pi^*)^\prime=\id\ .$$ 
\hB
Let $C^{-\infty}(\Lambda,V(\sigma_\lambda))$ denote the space of distribution sections of $V(\sigma_\lambda)$ with support in the limit set $\Lambda$.
Since $\Lambda$ is $\Gamma$-invariant
$C^{-\infty}(\Lambda,V(\sigma_\lambda))$ is a subrepresentation of 
the representation $\pi^{\sigma,\lambda}$ of $\Gamma$ on $C^{-\infty}(\partial X,V(\sigma_\lambda))$.
Note that ${}^\Gamma C^{-\infty}(\Lambda,V(\sigma_\lambda))$ is exactly the kernel of $res$.

\begin{lem}\label{mainkor}
If  ${}^\Gamma C^{-\infty}(\Lambda ,V(\sigma_\lambda))=0$ and if $ext$ is defined, then  we have $ext\circ res =\id$.
\end{lem}
\proof 
The assumption implies that $res$ is injective.
By Lemma \ref{iden} we have
$res(ext\circ res - \id)=0$. \hB

In order to apply this lemma we have to check its assumption.
In the course of the paper we will prove several vanishing results for
${}^\Gamma C^{-\infty}( \Lambda ,V(\sigma_\lambda))$.
One of these is available at this point.
\begin{lem}\label{pokm}
If $\Ree(\lambda)>0$ and $\Imm(\lambda)\not=0$, then   
${}^\Gamma C^{-\omega}( \Lambda ,V(\sigma_\lambda))=0$.
\end{lem}
\proof
We employ the Poisson transform.
Let $\gamma$ be a finite-dimensional representation of $K$ on $V_\gamma$
such that there exists an injective $T\in \Hom_M(V_\sigma,V_\gamma)$.
We will view sections of $V(\gamma)$ as functions from $G$ to
$V_\gamma$ satisfying the usual $K$-invariance condition. 
\begin{ddd}\label{defofpoi} The Poisson transform 
$$P:=P^T_\lambda:C^{-\infty}(\partial X,V(\sigma_\lambda))\rightarrow C^\infty(X,V(\gamma))$$
is defined by  
$$(P \phi)(g):=\int_K a(g^{-1}k)^{-(\rho+\lambda)}\gamma(\kappa(g^{-1}k)) T \phi(k) dk\ .$$
Here $\phi\in C^{-\infty}(\partial X, V(\sigma_\lambda))$ and the integral is a formal notation
meaning that the distribution $\phi$ has to be applied to the smooth integral kernel.
\end{ddd}
In \cite{olbrichdiss} it is shown that the Poisson transform $P$ intertwines
the left-regular representations of $G$ on $C^{-\infty}(\partial X,V(\sigma_\lambda))$ and $C^\infty(X,V(\gamma))$.
For $\Imm(\lambda)\not=0$ the principal series representation $\pi^{\sigma,\lambda}$ of $G$ on $C^{-\infty}(\partial X, V(\sigma_\lambda))$
is topologically irreducible. Since the Poisson transform $P$ does not vanish 
and is continuous it is injective.
There is a real constant $c_\sigma$ 
(see \cite{bunkeolbrich955})
such that $(-\Omega_G+c_\sigma+\lambda^2)P \phi =0$ 
(where $\lambda^2 :=\langle \lambda,\lambda\rangle $ with respect to the $\C$-linear scalar product induced on $\aca$ by the the invariant form on $\gaaa$).

 Let $V\subset \partial X$ and $U\subset X$ such that $\clo(U)\cap V=\emptyset$,
where we take the closure of $U$ in $X\cup\Omega$.
Then for $\Ree(\lambda)>0$ the integral kernel of the Poisson transform
$(g,k)\rightarrow p_\lambda(g,k):= a(g^{-1}k)^{-(\rho+\lambda)} \gamma(\kappa(g^{-1}k))T$
is a smooth function from $VM$ to $L^2(U,V(\gamma))\otimes V_{\tilde{\sigma}}$
(a much more detailed analysis of the Poisson kernel is given below in the proof of Lemma \ref{th43}).

If  $\phi\in {}^\Gamma C^{-\infty}( \Lambda ,V(\sigma_\lambda)) $,
then $P \phi $ is $\Gamma$-invariant, and since $\Ree(\lambda)>0$
it descends to a section in $L^2(Y,V_Y(\gamma))$.
Moreover it is annihilated by $(-\Omega_G+c_\sigma+\lambda^2)$. 
Since $Y$ is complete $\Omega_G$ is essentially selfadjoint on the domain $C_c^\infty(Y,V_Y(\gamma))$.
Its selfadjoint closure
has the domain of definition
$\{f\in L^2(Y,V_Y(\gamma))| \Omega_Gf\in L^2(Y,V_Y(\gamma))\}$.  
In particular, $\Omega_G$ can not have non-trivial eigenvectors in $L^2(Y,V_Y(\gamma))$
to eigenvalues with non-trivial imaginary part.
Since $\Imm(\lambda^2)\not=0$ we conclude that
$P\phi =0$ and hence $\phi=0$ by the injectivity of the Poisson transform.
\hB

\section{Meromorphic continuation of $ext$}\label{extttttt}

The extension $ext$ forms a holomorphic
family of maps depending on $\lambda\in\aca$
(we have omitted this dependence in order to simplify
the notation) which is defined now for $\Ree(\lambda)>\delta_\Gamma$.
In the present section we consider the meromorphic continuation of $ext$
to (almost) all of $\aca$. Our main tool is the scattering matrix
which we introduce below. The scattering matrix for the trivial
group $\Gamma=\{1\}$ is the Knapp-Stein intertwining operator
of the corresponding principal series representation.
We first recall basic properties of the Knapp-Stein intertwining operators.
Then we define the scattering matrix
using the extension and the Knapp-Stein intertwining operators. 
We simultaneously obtain the meromorphic continuations of the scattering matrix and the extension map.

If $\sigma$ is a representation of $M$, then 
we define its Weyl-conjugate $\sigma^w$ by $\sigma^w(m):=\sigma( w^{-1} mw)$,
where $w\in N_K(M)$ is a representative of the non-trivial element of the Weyl group $\cong \Z_2$ of $(\gaaa,\aaaa)$.
The Knapp-Stein intertwining operators form  meromorphic families of $G$-equivariant 
operators (see  \cite{knappstein71})
$$\hat{J}_{\sigma,\lambda}:C^{*}(\partial X,V(\sigma_\lambda)) \rightarrow C^{*}(\partial X ,V(\sigma^w_{-\lambda})),\quad *= -\infty,\infty\ .$$
Here  $\:\hat{}\:$ indicates that $\hat{J}_{\sigma,\lambda}$ is unnormalized.

In order to fix our conventions we give a definition of $\hat{J}_{\sigma,\lambda}$
as an integral operator acting on smooth functions for $\Ree(\lambda)<0$.
Its continuous extension to distributions is obtained by duality. 
For $\Ree(\lambda)\ge 0$ it is defined by meromorphic continuation.

Consider $f\in C^{\infty}(\partial X,V(\sigma_\lambda))$
as a function on $G$ with values
in $V_{\sigma_\lambda}$ satisfying the usual invariance condition with respect to $P$.
For $\Ree(\lambda)<0$ the intertwining operator is defined by the
convergent integral
\begin{equation}\label{furunkel}
(\hat{J}_{\sigma,\lambda}f)(g):=\int_{\bar{N}} f(g w\bar{n}) d\bar{n}\ .
\end{equation}

For all irreducible $\sigma\in\hat{M}$ we fix a minimal $K$-type (see \cite{knapp86}, Ch. XV for all that) of the principal series representation $C^{\infty}(\partial X,V(\sigma_\lambda))$. 
Let $c_{\sigma}(\lambda)$ be the value of $\hat{J}_{\sigma^w,-\lambda }$
on this minimal $K$-type.
Then $c_{ \sigma}(\lambda)$ is a meromorphic function on $\aca$ and we define
the normalized intertwining operators by
$$J_{ \sigma,\lambda}:= c_{  \sigma^w}^{-1}(-\lambda )\hat{J}_{\sigma,\lambda}\ .$$
Let $P_\sigma(\lambda):=c_\sigma(\lambda)^{-1}c_\sigma(-\lambda)^{-1}$
be  the Plancherel density.
Then the  intertwining operators satisfy the following functional equation. 
\begin{equation}\label{spex} \hat{J}_{\sigma,\lambda}\circ \hat{J}_{\sigma^w,-\lambda}=\frac{1}{P_\sigma(\lambda)},\quad \quad J_{\sigma^{w },-\lambda }\circ J_{ \sigma,\lambda}=\id\ .\end{equation}

Our next goal is to show that the intertwining operators form
a meromorphic family of operators in the sense defined in Section \ref{anaprep}.
This is an easy application of the approach to the intertwining operators
developed by Vogan-Wallach (see \cite{wallach92}, Ch. 10). The additional point
we have to verify is that in the domain of convergence of (\ref{furunkel}) the
operators $\hat{J}_{\sigma,\lambda}$ indeed form a continuous family.
\begin{lem}\label{iny1}
For $\Ree(\lambda)<0$ the intertwining operators 
$$\hat{J}_{\sigma,\lambda}:C^{\infty}(\partial X,V(\sigma_\lambda)) \rightarrow C^{\infty}(\partial X ,V(\sigma^w_{-\lambda}))$$
form a holomorphic family of continuous operators.
\end{lem}
\proof
Let $X_i$, $i=1,\dots,\dim(\kaaa)$, be an orthonormal base of $\kaaa$.
For any multiindex $r=(i_1,\dots,i_{\dim(\kaaa)})$ we set $X_r=\prod_{l=1}^{\dim(\kaaa)} X_l^{i_l}$, $|r|=\sum_{l=1}^{\dim(\kaaa)} i_l$,   and for $f\in C^{\infty}(K,V_{\sigma_\lambda})$
we define the seminorm
$$\|f\|_r=\sup_{k\in K} |f(X_rk)|\ .$$
It is well known that the system $\{\|.\|_r\}$, $r$ running over all multiindices,
defines the Fr\'echet topology of $C^{\infty}(K,V_{\sigma_\lambda})$ and by restriction the topology
of $C^{\infty}(\partial X,V(\sigma_\lambda))$.
 
We extend $f\in C^{\infty}(K,V_{\sigma_\lambda})$ to a function $f_\lambda$ on $G$ by setting
$f_\lambda(kan):=f(k)a^{\lambda-\rho}$.
Then we can define 
$$\hat{J}_{\sigma,\lambda} (f)(k)=\int_{\bar{N}} f_\lambda (k w\bar{n}) d\bar{n}\ .$$
For any $\lambda_0\in \aca$ with $\Ree(\lambda)<0$ and $\delta>0$ we can find
an $\epsilon>0$ such that for $|\lambda-\lambda_0|<\epsilon $
$$ \int_{\bar{N}} |a(\bar{n})^{\lambda_0-\rho}-a(\bar{n})^{\lambda-\rho}| d\bar{n} <\delta\ .$$
We then have
\begin{eqnarray*}
\|\hat{J}_{\sigma, \lambda_0}  f -\hat{J}_{\sigma,\lambda}  f \|_r&=&\sup_{k\in K} \int_{\bar{N}} (f_{\lambda_0} (X_rk w\bar{n}) - (f_{\lambda } (X_rk w\bar{n})  ) d\bar{n}\\
&=& \sup_{k\in K} \int_{\bar{N}} f(X_rkw\kappa(\bar{n})) (a(\bar{n})^{\lambda_0-\rho}-a(\bar{n})^{\lambda-\rho})d\bar{n}\\
&\le& \|f\|_r \int_{\bar{N}}   |a(\bar{n})^{\lambda_0-\rho}-a(\bar{n})^{\lambda-\rho}|d\bar{n}\\
&\le& \delta \|f\|_r
\end{eqnarray*}
This immediately implies that $\lambda\mapsto \hat{J}_{\sigma,\lambda}$
is a continuous family of operators on the space of smooth functions.
The fact that the family $\hat{J}_{\sigma,\lambda}$, $\Ree(\lambda)<0$, depends holomorphically on $\lambda$
is now easy to check (apply \cite{wallach92}, Lemma 10.1.3 and Fact \ref{holla}).
\hB 

\begin{lem}\label{jjj}
The family of intertwining operators 
$$\hat{J}_{\sigma,\lambda}:C^{\infty}(\partial X,V(\sigma_\lambda)) \rightarrow C^{\infty}(\partial X ,V(\sigma^w_{-\lambda}))$$
extends   meromorphically to all of $\aca$. 
\end{lem}
\proof
We employ \cite{wallach92}, Thm. 10.1.5, which states that there are polynomial
maps $b:\aca\rightarrow \C$ and $D:\aca\rightarrow \cU(\gaaa)^K$,
such that
\begin{equation}\label{shifty}
b(\lambda)\hat{J}_{\sigma,\lambda}=\hat{J}_{\sigma, \lambda-4\rho}\circ \pi^{\sigma,\lambda-4\rho}(D(\lambda))\ .
\end{equation}
This formula requires some explanation.
We  identify $$C^{\infty}(\partial X, V(\sigma_\lambda))\cong C^\infty(K,V_\sigma)^M$$
canonically. Then all operators act on the same space $C^\infty(K,V_\sigma)^M$.
 
If we know that $\hat{J}_{\sigma,\lambda}$ is meromorphic up to $\Ree(\lambda)<\mu$,
then we conclude that
$$\hat{J}_{\sigma,\lambda}= \frac{1}{b(\lambda) } \hat{J}_{\sigma,\lambda-4\rho}\circ \pi^{\sigma,\lambda-4\rho}(D(\lambda))$$
is meromorphic up to $\Ree(\lambda)<\mu+4\rho$.
Thus the lemma follows from Lemma \ref{iny1}.
\hB

We call $\lambda\in \aca$ integral if
$$ 2 \frac{\langle \lambda,\alpha\rangle }{\langle\alpha,\alpha\rangle} \in\Z\quad  \mbox{for some root $\alpha$ of $(\gaaa,\aaaa)$}\ ,$$
 and non-integral otherwise.
The integral points form a discrete lattice in $\aaaa^*$.
It is known that the poles of $\hat{J}_{\sigma,\lambda}$ have at most order one  and appear only at integral $\lambda$ (\cite{knappstein71}, Thm. 3 and Prop. 43).
\begin{lem}\label{off}
Let $\chi,\phi\in C^\infty(\partial X)$ such that $\supp(\phi)\cap \supp(\chi)=\emptyset$.
Then $\chi\hat{J}_{\sigma,\lambda}\phi$ is a holomorphic family of smoothing operators.
In particular, the residues of $\hat{J}_{\sigma,\lambda}$ are differential
operators.
\end{lem}
\proof 
Since $\supp(\phi)\cap\supp(\chi)=\emptyset$,
there exists a compact set $V\subset \bar{N}$ such that
$$\kappa(\supp(\chi)M w (\bar{N}\setminus V)M \subset (\partial X\setminus \supp(\phi))M\ .$$
For $\Ree(\lambda)<0$ and $f\in C^\infty(\partial X, V(\sigma_\lambda))$
we have (viewing $f$ as a function on $K$ with values in $V_{\sigma_\lambda})$
\begin{eqnarray*}
(\chi \hat{J}_{\sigma,\lambda} \phi f)(k) &=& \int_{\bar{N}} \chi(k)f(\kappa(kw\bar{n}))\phi(\kappa(kw\bar{n})) a(\bar{n})^{ \lambda-\rho}d\bar{n}\\
&=& \int_{ V }\chi(k) f(\kappa(kw\bar{n})) \phi(\kappa(kw\bar{n})) a(\bar{n})^{ \lambda-\rho}d\bar{n}\ .
\end{eqnarray*}
The right-hand side of this equation extends to all of $\aca$
and defines a holomorphic family of operators.
This proves the first part of the lemma.
It in particular implies that the residues of $\hat{J}_{\sigma,\lambda}$ are local operators.
Hence the second assertion follows.
\hB

We need the following consequence of Lemma \ref{off}.
Let $W\subset \partial X$ be a closed
subset and let
$$\cG_\lambda:=\{f\in C^{-\infty}(\partial X, V(\sigma_\lambda))|
f_{|\partial X\setminus W}\in  C^{\infty}(\partial X\setminus W, V(\sigma_\lambda))\}\ .
$$
We equip $\cG_\lambda$ with the weakest topology such that the maps
$\cG_\lambda\hookrightarrow C^{-\infty}(\partial X, V(\sigma_\lambda))$
and $\cG_\lambda\rightarrow  C^{\infty}(\partial X\setminus W, V(\sigma_\lambda))$
are continuous.
Let $U\subset \bar{U}\subset \partial X\setminus W$ be open.
\begin{lem}\label{dissmo}
The composition 
$$res_U \circ \hat{J}_{\sigma,\lambda}:\cG_\lambda\rightarrow C^\infty(U,V(\sigma^w_{-\lambda}))$$
is well-defined and depends meromorphically on $\lambda$.
\end{lem}

We introduce a notational convention concerning $\sigma$.
Below $\sigma$ shall always denote a Weyl-invariant representation
of $M$ which is either irreducible or of the
form $\sigma^\prime\oplus \sigma^{\prime w}$ with $\sigma^\prime$ irreducible 
and not Weyl-invariant.
In the latter case $c_\sigma(\lambda):=c_{\sigma^\prime}(\lambda)=c_{\sigma^{\prime w}}(\lambda)$, $P_\sigma(\lambda):=P_{\sigma^\prime}(\lambda)=P_{\sigma^{\prime w} }(\lambda)$.
We omit the subscript $\sigma$ in the notation of the intertwining operators.

We  now turn to the definition of the (normalized) scattering matrix
as a family of operators
$$\hat{S}_\lambda \:( S_{\lambda})\: :C^{ *}(B,V_B(\sigma_\lambda))\rightarrow C^{*}(B,V_B(\sigma_{-\lambda })),\quad *=\infty,-\infty \ .$$

\begin{ddd}\label{scatdef}
For $\Ree(\lambda)>\delta_\Gamma$ we define
\begin{equation}\label{scatde}
\hat{S}_{ \lambda} :=res\circ \hat{J}_{\lambda}\circ ext\quad, \quad 
S_{ \lambda} :=res\circ J_{\lambda}\circ ext \ .\end{equation}
\end{ddd}

\begin{lem}\label{kkk}
For $\Ree(\lambda)>\delta_\Gamma$ 
the scattering matrix forms a meromorphic family of operators
$$C^{\pm\infty}(B,V_B(\sigma_\lambda))\rightarrow C^{\pm\infty}(B,V_B(\sigma_{-\lambda }))\ .$$ If $\hat{S}_\lambda$ is singular and  $\Ree(\lambda)>\delta_\Gamma$, then $\lambda$ is integral and the residue of $\hat{S}_\lambda$ is a differential operator.
\end{lem}
\proof
The assertion for the scattering matrix acting on distributions follows from 
Lemma \ref{jjj}, Lemma \ref{off} and   Fact \ref{comp}.
The fact that the scattering matrix restricts to smooth sections follows
from Lemma \ref{dissmo}.\hB

\begin{lem}\label{lok}
If  $\Ree(\lambda)>\delta_\Gamma$, then the adjoint
$${}^tS_{ \lambda}: C^{\infty}(B,V_B(\tilde{\sigma}_\lambda))\rightarrow C^{\infty}(B,V_B(\tilde{\sigma}_{-\lambda }))$$
of
$$S_{ \lambda}: C^{-\infty}(B,V_B(\sigma_\lambda))\rightarrow C^{-\infty}(B,V_B(\sigma_{-\lambda }))$$
coincides with the restriction of $S_\lambda$ to $C^{\infty}(B,V_B(\tilde{\sigma}_\lambda))$.
\end{lem}
\proof
We employ the fact that the corresponding relation holds for the intertwining operators (step (\ref{zx2}) below).
Let $F\subset \Omega$ be a fundamental domain of $\Gamma$ and $\chi\in C_c^\infty(\Omega)$
be a cut-off function with $F\subset \supp(\chi)$ and $\sum_{g\in\Gamma} g^*\chi\equiv 1$ on $\Omega$.
Let $\phi\in C^{\infty}(B,V_B(\sigma_\lambda))$,
$f\in C^{\infty}(B,V_B(\tilde{\sigma}_\lambda))$,
and consider $\phi$ as a distribution section.
Then
\begin{eqnarray}
\langle \phi ,{}^tS_\lambda f\rangle &=&\langle  S_\lambda\phi ,  f\rangle 
 = \langle  res\circ J_\lambda\circ ext \phi ,  f\rangle\nonumber\\
&=&\langle J_\lambda\circ ext \phi , \chi f   \rangle 
 = \langle  ext \phi , {}^t J_\lambda \chi f  \rangle \label{zx1}\\
&=&\langle   ext \phi , J_\lambda \chi f  \rangle\label{zx2}\\
&=&\langle   \chi  \phi ,  \sum_{g\in\Gamma } res_\Omega \circ \pi(g) J_\lambda \chi   f\rangle 
 =  \sum_{g\in\Gamma } \langle   \chi  \phi ,  res_\Omega \circ \pi(g) J_\lambda \chi   f\rangle\label{zx3}\\
&=& \sum_{g\in\Gamma } \langle   \chi  \phi ,   J_\lambda \pi(g)   \chi f  \rangle 
 =  \sum_{g\in\Gamma } \langle  \pi(g)    J_\lambda \chi  \phi ,    \chi  f\rangle\nonumber\\
&=&   \langle \phi , S_\lambda f\rangle\ . \label{zx4}
\end{eqnarray}
Here in (\ref{zx1}) and (\ref{zx3}) we view $\chi f$ and $\chi \phi$ as
sections over $\Omega$, respectively, and $\pi(g)$ is induced by the corresponding principal 
series representations.
In order to obtain (\ref{zx4}) from the preceding line we do the transformations
backwards with the roles of $\phi$ and $f$ interchanged.
\hB

\begin{lem}\label{pfun}
If $|\Ree(\lambda)|<- \delta_\Gamma$,
then the scattering matrix satisfies the functional equation
(viewed as an identity of meromorphic families of operators)
$$S_{-\lambda}\circ S_\lambda = \id\ .$$
\end{lem}
\proof
We employ Lemmas \ref{mainkor}, \ref{pokm}, and (\ref{spex}) in order to compute for $\Imm(\lambda)\not= 0$, $\Ree(\lambda)\not=0$, 
\begin{eqnarray*}
S_{ -\lambda }\circ S_{ \lambda}&=&res\circ J_{ -\lambda }\circ ext\circ res \circ J_{ \lambda}\circ ext\\
&=&res\circ J_{ -\lambda }\circ J_{ \lambda}\circ ext \\
&=&res\circ ext\\
&=& \id\ .
\end{eqnarray*}
This identity now extends meromorphically to all of $\{|\Ree(\lambda)|< -\delta_\Gamma\}$.
\hB

Now we start with the main topic of the present section, the meromorphic continuation of $ext$.
We first invoke the meromorphic Fredholm theory  Lemma \ref{merofred}
in order to provide a meromorphic continuation of the scattering matrix 
to almost all of $\aca$ under the condition $\delta_\Gamma<0$.
We then use this continuation of the scattering matrix in order to define
the meromorphic continuation of $ext$. 
Finally, if $X$ is not the Cayley hyperbolic plane, we drop the assumption  
$\delta_\Gamma<0$.

To be more precise our method for the meromorphic continuations
breaks down at a countable number of points.
Therefore we introduce the set
\begin{equation}\label{acagood}\aca(\sigma):=\{\lambda\in\aca\:|\: \Ree(\lambda)>\delta_\Gamma\:\mbox{or}\:
P_\sigma(\lambda) \hat{J}_{-\lambda}\:\mbox{is regular}\}\ .\end{equation}
Not only the poles of $\hat{J}_{-\lambda}$, but also the poles of $P_\sigma(\lambda)$ are located at integral $\lambda$ (see e.g. the explicit formla (\ref{blanche}) given in the proof of Lemma \ref{casebycase}). Therefore, $\aca(\sigma)$ contains all non-integral
points, and many of the integral points, too.
If $\delta_\Gamma<0$, then we show the meromorphic continuation of $ext$
and of the scattering matrix for $\lambda\in \aca(\sigma)$.
This result is employed in the proof of Proposition \ref{mystic}.
If $\delta_\Gamma\ge 0$, then for simplicity 
we show the meromorphic continuation of
$ext$ and the scattering matrix to the set of all
non-integral $\lambda\in\aca$, only.

\begin{prop}\label{part1}
The scattering matrix
$$\hat{S}_\lambda :C^{\pm\infty}(B,V_B(\sigma_\lambda))\rightarrow C^{\pm\infty}(B,V_B(\sigma_{-\lambda }))$$
and  the extension map 
$$ext: C^{-\infty}(B,V_B(\sigma_\lambda))\rightarrow {}^\Gamma C^{-\infty}(\partial X,V(\sigma_\lambda))$$
have  meromorphic continuations to 
 the set of all non-integral $\lambda\in\aca$. 
In particular, we have    
\begin{equation}\label{forme}
ext  = J_{-\lambda} \circ ext \circ S_{\lambda}, \quad  S_{-\lambda}\circ S_\lambda = \id\ .\end{equation}  
Moreover,
$ext$ and $S_\lambda$ have at most finite-dimensional singularities
at non-integral $\lambda$. 
\end{prop}
\proof
We first assume that $\delta_\Gamma<0$.
We construct the meromorphic continuation of 
$$S_\lambda:C^{\infty}(B,V_B(\sigma_\lambda))\rightarrow C^{\infty}(B,V_B(\sigma_{-\lambda }))\ ,$$ and then we extend this continuation to distributions
by duality using Lemma \ref{lok}. The idea is to set $S_\lambda:=S_{-\lambda}^{-1}$ for $\Ree(\lambda)<-\delta_\Gamma$ and to show that $S_{-\lambda}^{-1}$ forms
a meromorphic family. 

Let $\{U_\alpha\}$ be a finite open covering of $B$ and let $\tilde{U}_\alpha$ be
diffeomorphic  lifts of $U_\alpha$.
Choose a subordinated partition of unity $\phi_\alpha$.
We view $\phi_\alpha$ as a smooth compactly supported function on $\tilde{U}_\alpha$. For $h\in\Gamma$ we set $\phi_\alpha^h(x):=\phi_\alpha(h^{-1}x)$.
Let $1\in L\subset \Gamma$ be a finite subset. Then we define
$\chi\in {}^\Gamma C^\infty(\Omega\times\Omega)$ by 
$$\chi(x,y):=\sum_{g\in\Gamma,h\in L,\alpha} \phi_\alpha(gx)\phi_\alpha^h(gy)\ .$$ 
Let 
$$\hat{J}^{diag}_{\lambda}:C^\infty(B, V_B(\sigma_\lambda))\rightarrow C^\infty(B, V_B(\sigma_{-\lambda}))$$ be the meromorphic family
of operators obtained by multiplying
the distributional kernel of $\hat{J}_{ \lambda}$ by $\chi$.
If $f\in C^\infty(B, V_B(\sigma_\lambda))$, then
$$(\hat{J}^{diag}_{ \lambda})f =\sum_{\alpha,h\in L} \phi_\alpha \hat{J}_\lambda(\phi^h_\alpha f)$$
using the canonical identifications.
Below we shall employ the fact that $\hat{J}^{diag}_{ \lambda}$ depends
on $L$.

Let  $$U:=\{\lambda\in\aca\:|\:\Ree(\lambda)>\delta_\Gamma,\: -\lambda\in\aca(\sigma),\: \lambda\:\mbox{non-integral if}\:\Ree(\lambda)\le 0\}\ .$$ Then $U$ is open and
connected.
For $\lambda\in U$ define 
\begin{equation}\label{klio}R(\lambda):=P_\sigma(\lambda) \hat{J}^{diag}_{-\lambda} \circ \hat{S}_\lambda - \id\ .\end{equation}
 The inverse of the unnormalized scattering matrix for $\lambda\in U$
should be  given by
\begin{equation}\label{finit}\hat{S}_{ \lambda}^{-1}=P_\sigma(\lambda) (\id+ R(\lambda))^{-1}\circ \hat{J}^{diag}_{ -\lambda }\ .\end{equation}
It exists as a meromorphic family if $(\id+ R(\lambda))^{-1}$ does.

We want to apply the  meromorphic Fredholm theory (Proposition \ref{merofred})
in order to invert $\id +R(\lambda)$ for $\lambda\in U$
and to conclude that $(\id+R(\lambda))^{-1}$
is meromorphic. 
 
We check the assumption of Proposition \ref{merofred}.
We choose a Hermitian metric on $V_B(\sigma_0)$ and a volume form on $\Omega$.
The Hilbert space $\cH$ of Proposition \ref{merofred}
is $L^2(B,V_B(\sigma_0))$ defined using these choices.  
The Fr\'echet space $\cF$ is just $C^\infty(B, V_B(\sigma_0))$.
Implicitly, we identify the spaces $C^\infty(B, V_B(\sigma_\lambda))$ with $C^\infty(B, V_B(\sigma_0))$ using a trivialization of the holomorphic
family of bundles $\{V_B(\sigma_\lambda)\}$, $\lambda\in\aca$.

We claim that $R(\lambda)$ is a holomorphic family
of smoothing operators on $U$.
If $\lambda\in U$, then
$\hat{J}^{diag}_{-\lambda}$ as well as
$$P_\sigma(\lambda)\hat{S}_\lambda=res\circ P_\sigma(\lambda)\hat{J}_\lambda\circ ext$$ are regular.
Moreover, $R(\lambda)$ is smoothing by Lemma \ref{off}. 
Hence $R(\lambda)$ is indeed a holomorphic family of smoothing operators 
and this proves the claim. 
 
Next we show that if $L\subset \Gamma$ is sufficiently exhausting, then  
$\id+ R(\lambda)$ is injective for some $\lambda\in U$.
Since $\id+R(\lambda)$ is Fredholm of index zero it is then invertible at this point. 
Here is one of the two instances where we assume $\delta_\Gamma<0$.
We fix some non-integral $\lambda\in \aca$ with $|\Ree(\lambda)|<-\delta_\Gamma$.
Define $\hat{J}^{off}_{-\lambda}:=res\circ\hat{J}_{-\lambda}-\hat{J}^{diag}_{-\lambda}
\circ res$. By Lemma \ref{off} the composition
$R(\lambda) = -  P_\sigma(\lambda) \hat{J}^{off}_{-\lambda} \circ \hat{J}_{\lambda}\circ ext$ 
is a bounded operator on $C^k(B,V_B(\sigma_\lambda))$.
The proof of
Lemma \ref{anal1} shows that for $\Ree(\lambda)>\delta_\Gamma$
the push  down is a continuous map
$\pi_\ast:C^0(\partial X,V(\tilde{\sigma}_{-\lambda}))\rightarrow
C^0(B,V_B(\tilde{\sigma}_{-\lambda}))$.
Thus the adjoint $ext$ restricts to a continuous map between
the Banach spaces of measures of bounded variation
$$ext: M_b(B,V_B(\sigma_\lambda))\rightarrow M_b(\partial X,V(\sigma_\lambda))\ .$$
The operator $\hat{J}_\lambda$ is a singular
integral operator composed with a differential operator
(as explained in Lemma \ref{jjj}).
Thus there is a $k\in\nat_0$ such that
$$\hat{J}_\lambda:M_b(\partial X,V(\sigma_\lambda))\rightarrow C^{k}(\partial X,V(\sigma_{\lambda}))^\prime$$
is continuous.
The scattering matrix $\hat{S}_\lambda$
extends to  a continuous operator
$$\hat{S}_\lambda:M_b(B,V_B(\sigma_\lambda))\rightarrow C^{k}(B,V_B(\sigma_{\lambda}))^\prime$$
and, by dualization and Lemma \ref{lok}, to a continuous map
$$\hat{S}_\lambda:C^k(B,V_B(\sigma_\lambda))\rightarrow  C^{0}(B,V_B(\sigma_{-\lambda})) \ .$$

\begin{lem}\label{converg1} 
If $L$ runs over an increasing sequence of subsets 
exhausting $\Gamma$, then
$(\hat{J}^{diag}_{-\lambda}-\hat{S}_{-\lambda})\to 0$
in the sense of bounded operators from $C^{0}(B,V_B(\sigma_{-\lambda}))$
to $C^{k}(B,V_B(\sigma_{\lambda}))$.
 \end{lem}
\proof
It follows from the estimates proved in Lemma \ref{anal1} that
$(f\mapsto \sum_{h\in L }  \phi^h_\alpha f)$ tends to $(f\mapsto ext(f))$
in the sense of bounded operators  between the
Banach spaces $C^0(B,V_B(\sigma_\lambda))$ and
$M_b(\partial X, V(\sigma_\lambda ))$.
We apply Lemma \ref{dissmo} letting $W$ be a compact neighbourhood
of $\Lambda$ and $U\subset \Omega$ contain $\tilde{U}_\alpha$ for all
$\alpha$. 
It implies that 
$$\sum_{h\in \Gamma\setminus L }\phi_\alpha \hat{J}_{-\lambda}   \phi^h_\alpha\to 0$$  
in the sense of  bounded operators from  $C^0(B,V_B(\sigma_{-\lambda}))$
to $C^{k}(\partial X,V(\sigma_\lambda))$ for all $\alpha$.
The assertion of the lemma now follows.
\hB

By Lemmas \ref{converg1}, \ref{pfun} and Equation (\ref{spex}) the operator $R(\lambda)$ tends to zero in the sense of bounded operators on $C^k(B,V_B(\sigma_\lambda))$
when $L$ runs over an increasing sequence of subsets exausting $\Gamma$.
Thus if $L$ is large enough, then $\id+R(\lambda)$ is injective.

We now have verified the assumptions of Proposition \ref{merofred}.
We conclude that 
the family  $\hat{S}_{\lambda}^{-1}$ is meromorphic for $\lambda\in U$.
If $\hat{S}_{\lambda}^{-1}$ has a singularity and $\lambda$ is non-integral,
then this singularity is finite-dimensional.

Here is the second and main instance, where we need the assumption $\delta_\Gamma<0$. 
Namely, it implies that 
$$\{\lambda\in\aca\:|\:-\lambda\in U\}\cup \{\lambda\in\aca\:|\:\Ree(\lambda)>\delta_\Gamma\}=\aca(\sigma)\ .$$
Furthermore, by Lemma \ref{pfun} we have $S_\lambda=S_{-\lambda}^{-1}$ on $\{\lambda\in\aca\:|\:-\lambda\in U\}\cap \{\lambda\in\aca\:|\:\Ree(\lambda)>\delta_\Gamma\}$. Thus, setting
$S_\lambda:=S_{-\lambda}^{-1}$ for $-\lambda\in U$ we obtain a well-defined continuation of $S_\lambda$ to all of $\aca(\sigma)$.
By duality  this  continuation extends to distributions still having the same finite-dimensional
singularities at non-integral points.

It remains to consider the extension map. 
We  employ the scattering matrix in order to define  for 
$\lambda\in\aca(\sigma)$,
$\Ree(\lambda)<-\delta_\Gamma$
$$
ext_1 := J_{-\lambda} \circ ext \circ S_{\lambda} \ .$$
We claim that $ext=ext_1$. In fact since $res$ is injective on an open subset of $\{|\Ree(\lambda)| < -\delta_\Gamma\}$, the computation 
\begin{eqnarray*}
res\circ ext_1&=& res\circ J_{-\lambda} \circ ext \circ S_\lambda\\
&=&  S_{-\lambda}\circ S_\lambda   \\
&=&\id
\end{eqnarray*}
implies the claim.

We now have constructed a meromorphic continuation of $ext$ to all of 
$\aca(\sigma)$. The relation (\ref{forme})
between the scattering matrix and $ext$ follows by meromorphic continuation. This equation also implies that $ext$  has at most finite-dimensional singularities.
We have finished the proof of Proposition \ref{part1} assuming $\delta_\Gamma<0$.
The identities
$$\hat{S}_\lambda =res\circ \hat{J}_\lambda \circ ext, \quad  S_\lambda\circ S_{-\lambda}=\id $$ extend to all of   $\aca(\sigma)$ by meromorphic continuation.
 
We now show how to drop the assumption $\delta_\Gamma<0$ using the embedding trick
and Assumption \ref{caly}. 
If $X$ is the Cayley hyperbolic plane, then by assumption
$\delta_\Gamma<0$ and the proposition is already proved.
Thus we can assume that $X$ belongs to a series of rank-one symmetric
spaces.
Let $\dots \subset G^n\subset G^{n+1}\subset \dots$ be the corresponding sequence 
of real, semisimple, linear Lie groups inducing embeddings of the corresponding
Iwasawa constituents $K^n\subset K^{n+1}$,
$N^n  \stackrel{\scriptstyle  \subset}{ \scriptstyle \not=}  N^{n+1}$, $ M^n\subset M^{n+1}$
such that $A=A^n = A^{n+1}$.
Then we have totally geodesic embeddings of the symmetric spaces
$X^n\subset X^{n+1}$
inducing embeddings of their boundaries
$\partial X^n\subset \partial X^{n+1}$.
If $\Gamma\subset G^n$ satisfies \ref{asss} then it keeps satisfying \ref{asss}
when viewed  as a subgroup of $G^{n+1}$.
We obtain  embeddings
$\Omega^n\subset  \Omega^{n+1}$ inducing
$B^n\subset  B^{n+1}$
while the limit set $ \Lambda^n$ is identified with $ \Lambda^{n+1}$.
Let $\rho^n(H)=\frac{1}{2}\tr(\ad(H)_{|\naaa^n})$, $H\in\aaaa$.

The exponent of $\Gamma$ now depends on $n$ and is denoted by
$\delta_\Gamma^n$.
We have the relation $\delta_\Gamma^{n+1}=\delta_\Gamma^n-\rho^{n+1}+\rho^n$.
Thus $\delta_\Gamma^{n+m}\to -\infty$ as $m\to 0$ and hence
taking $m$ large enough we can satisfy $\delta_\Gamma^{n+m}<0$.
 The aim of the following discussion is to show how the meromorphic continuation  $ext^{n+1}$  leads to the continuation of  $ext^n$.

Let $\sigma^{n+1}$ be a Weyl-invariant representation of $M^{n+1}$. Then it
restricts to a Weyl-invariant representation of $M^n$.
For any given finite-dimensional representation $\sigma^n$ of $M^n$ we can find
a Weyl-invariant representation $\sigma^{n+1}$ of $M^{n+1}$ such that $\sigma^{n+1}_{|M^n}$
contains $\sigma^n$ as a subrepresentation.

The representation $\sigma^{n+1}_\lambda$ of $P^{n+1}$ restricts to the representation
$(\sigma^{n+1}_{|P^n})_{\lambda+\rho^n-\rho^{n+1}}$ of $P^n$.
This induces an isomorphism of bundles
$$V_{B^{n+1}}(\sigma^{n+1}_\lambda)_{|B^n}=V_{B^n}((\sigma^{n+1}_{|P^n})_{\lambda+\rho^n-\rho^{n+1}})\ .$$
We will omit the subscript ${}_{|P^n}$ and the superscript ${}^{n+1}$ of $\sigma$ in the following discussion.
We obtain a push forward of distributions
$$i_*:C^{-\infty}(B^n,V_{B^n}(\sigma_\lambda))\rightarrow C^{-\infty}(B^{n+1},V_{B^{n+1}}(\sigma_{\lambda+\rho^n-\rho^{n+1}}))\ .$$
For $\phi\in C^{-\infty}(B^n,V_{B^n}(\sigma_\lambda))$ the push forward
$i_\ast(\phi)$ has support in $B^n\subset B^{n+1}$. 
Then also 
$\supp(ext(\phi))\subset  \partial X^n$.
We try to define a pull back
$ext^n(\phi):=i^*\circ ext^{n+1}\circ i_*(\phi)$ as follows.
For $f\in C^\infty(\partial X^{n},V(\tilde{\sigma}_{-\lambda}))$ let $\tilde{f}\in C^\infty(\partial X^{n+1},V(\tilde{\sigma}_{-\lambda-\rho^n+\rho^{n+1}}))$
be an arbitrary extension. Then we set
$$\langle ext^n(\phi),f\rangle := \langle ext^{n+1}\circ i_*(\phi_\lambda),\tilde{f}\rangle\ .$$
\begin{lem}\label{well}
 $ext^n$ is well defined.
\end{lem}
We must show that this definition does not depend on the choice
of the extension of $\tilde{f}$.
It is sufficient to show that if $h\in C^\infty(\partial X^{n+1},V(\tilde{\sigma}_{-\lambda-\rho^n+\rho^{n+1}}))$
vanishes on $\partial X^{n}$, then $\langle ext^{n+1}\circ i_*(\phi_\lambda),h\rangle=0$.
Embed $\phi,h$ into holomorphic families $\phi_\mu, h_\mu$
, $\phi_\mu\in C^{-\infty}(B^n,V_{B^n}(\sigma_\mu))$, $h_\mu\in C^\infty(\partial X^{n+1},V(\tilde{\sigma}_{-\mu-\rho^n+\rho^{n+1}}))$
such that $(h_\mu)_{|\partial X^{n} }=0$.
Then for $\Ree(\mu)$ large enough\linebreak[4] $\pi^{n+1}_*(h_{-\mu})_{|B^n}=0$
and thus $$\langle ext^{n+1}\circ i_*(\phi_\mu),h_{-\mu} \rangle= \langle   i_*(\phi_\mu),\pi^{n+1}_*h_{-\mu}\rangle=0\ .$$
By meromorphic continuation this identity holds for all $\mu$, in particular at $\mu=\lambda$.
\hB
 
If $\lambda\in\aca$ is non-integral, then so is $\lambda+\rho^n-\rho^{n+1}$.
We deduce the properties of $ext^{n}$ from the corresponding 
properties of $ext^{n+1}$.
In particular, $ext^n$ is continuous, meromorphic, and has at most finite-dimensional singularities  at $\lambda$ if $ext^{n+1}$ has these properties at $\lambda+\rho^n-\rho^{n+1}$.
We define the meromorphic continuation of the scattering matrix by
(\ref{scatde}).
Then it is easy to see that the scattering matrix has the properties
as asserted.
This finishes the proof of Proposition \ref{part1}.
\hB

\begin{lem}\label{lead}
If $ext$ is meromorphic at $\lambda\in \aca$ and ${}^\Gamma C^{-\infty}(\Lambda,V(\sigma_\lambda))=0$, then $ext$\
is regular at $\lambda$.
\end{lem}
\proof
If $ext$ is meromorphic, then for any holomorphic family $\mu\to\phi_\mu\in C^{-\infty}(B,V_B(\sigma_\mu))$ the leading singular part of 
$ext(\phi_\mu)$ at $\mu=\lambda$ belongs to ${}^\Gamma C^{-\infty}(\Lambda,V(\sigma_\lambda))$.
In fact $res\circ ext(\phi_\mu)=\phi_\mu$ has no singularity.
If ${}^\Gamma C^{-\infty}(\Lambda,V(\sigma_\lambda))=0$, 
then $ext(\phi_\mu)$ is regular at $\mu=\lambda$ for any  holomorphic family $\mu\to \phi_\mu$.\hB

In Lemma \ref{ghu} below we consider the converse of Lemma \ref{lead}
for non-integral $\lambda\in\aca$ with $\Ree(\lambda)>0$.

\begin{lem}\label{firtu}
If $\Ree(\lambda)>0$ and $ext:C^{-\infty}(B,V_B(\sigma_\mu))\rightarrow {}^\Gamma C^{-\infty}(\partial X,V(\sigma_\mu))$ is meromorphic at $\mu=\lambda$, 
then the order of a singularity of $ext$ at $\lambda$ is at most $1$. 
\end{lem}
\proof
Let $\gamma\in \hat{K}$ be such that there exists an injective $T\in \Hom_M(V_\sigma,V_\gamma)$.
We also can and will require that
the Poisson transform $P_\lambda^T=:P$ is injective
(e.g. by taking $\gamma$ to be the minimal $K$-type of the
principal series representation $\pi^{\sigma,\lambda}$).

Let $f_\mu\in C^\infty(B,V_B(\sigma_\mu))$, $\mu\in \aca$,
be a holomorphic family such that $ext(f_\mu)$ has a pole of order
$n\ge 1$ at $\mu=\lambda$, $\Ree(\lambda)>0$.

We assume that $n\ge 2$ and argue by contradiction.
Let $0\not=\phi \in {}^\Gamma C^{-\infty}(\Lambda,V(\sigma_\lambda))$
be the leading singular part of $ext(f_\mu)$ at $\mu=\lambda$.
Then $(\lambda^2-\mu^2)^{n-1} ext(f_\mu)$ and hence
$(\lambda^2-\mu^2)^{n-1} P_\mu^T ext(f_\mu)$
have  first-order poles, the latter with residue $-(2\lambda)^{n-1} P\phi $.

Since $res_\Omega\circ ext(f_\mu)$ is smooth if $kM\in\Omega M$ we have  $P^T_\mu\circ ext(f_\mu) (ka)=O(a^{\mu-\rho})$.
Moreover $P\phi (ka)=O(a^{-\lambda-\rho})$ and both estimates
hold uniformly for $kM$ in compact subsets of $\Omega$, $\mu$ near $\lambda$, and large $a\in A_+$.
This justifies the following computation using  partial integration:
\begin{eqnarray*}
\infty&=&\lim_{\mu\to\lambda,\: \Ree(\mu)<\Ree(\lambda)} \langle 
(\lambda^2-\mu^2)^{n-1}  P_\mu^T ext(f_\mu), P\phi \rangle_{L^2(Y)}\\
&=&\lim_{\mu\to\lambda,\: \Ree(\mu)<\Ree(\lambda)} \langle
(-\Omega_G+c_\sigma+\lambda^2)^{n-1} P_\mu^T ext(f_\mu), P\phi \rangle_{L^2(Y)}\\
&=& 0
\end{eqnarray*}
This is a contradiction and thus $n=1$.
\hB

Since $\sigma$ is an unitary representation of $M$, we have
for $\lambda\in \imath\aaaa^*$ a positive conjugate linear
pairing $V_{\sigma_\lambda}\otimes V_{\sigma_\lambda} \rightarrow V_{1_{-\rho}}$
and hence a natural $L^2$-scalar product $C^\infty(B,V_B(\sigma_\lambda))\otimes C^\infty(B,V_B(\sigma_\lambda))\rightarrow \C$.
Let $L^2(B,V_B(\sigma_\lambda))$ be associated Hilbert space.
Using Lemma \ref{lok} we see that the adjoint $S^*_\lambda$
with respect to this Hilbert space structure is just $S_{-\lambda}$.

\begin{lem}\label{unitary}
If $\Ree(\lambda)=0$, $\lambda\not=0$, then $S_\lambda$ is regular and unitary.
\end{lem}
\proof
The scattering matrix $S_\lambda$ is meromorphic at non-zero imaginary points $\lambda$. Let now $\lambda$ be imaginary, $S_{\pm\lambda}$ be regular and  
$f\in C^\infty(B,V_B(\sigma_\lambda))$. Then by the functional equation (\ref{forme})
$$\|S_\lambda f\|_{L^2(B,V_B(\sigma_\lambda))}^2=\langle S_{-\lambda}\circ S_\lambda f, f \rangle_{L^2(B,V_B(\sigma_\lambda))} = \|f\|^2_{L^2(B,V_B(\sigma_\lambda))}\ .$$
This equation remains valid at all non-zero imaginary points. Hence, $S_\lambda$
is regular and unitary there.
\hB

 \section{Invariant distributions on the limit set}\label{invvv}

In present section we study the space 
${}^\Gamma C^{-\infty}(\Lambda,V(\sigma_\lambda))$ of
invariant distributions which are
supported on the limit set, mainly for $\Ree(\lambda)\ge 0$. 
We show that nontrivial distributions
of this kind can only exist for a countable set of parameters $\lambda\ge 0$ with possibly finitely many accumulation points.
In particular, we show that
${}^\Gamma C^{-\infty}(\Lambda,V(\sigma_\lambda))=0$ if $ext$ is regular at $\lambda$, $\Ree(\lambda)\ge 0$, and $\chi_{\mu_\sigma+\rho_m-\lambda}$ is non-integral. 

In the course of this paper we will prove the finiteness of the discrete spectrum of $\cZ$ on
$L^2(Y,V_Y(\gamma))$, ore equivalently, that nontrivial invariant distributions
with support on the limit set (with $\Ree(\lambda)\ge 0$) can in fact only exist for a finite set of parameters $\lambda\ge 0$. The proof of the finiteness of the point spectrum
is essentially based on the spectral comparison Proposition \ref{esspec}. But this proposition is not applicable in order to exclude that eigenvalues accumulate at the boundary of the continuous spectrum.
Here we will employ Corollary \ref{nahenull} instead,
and it is important to show that $ext$ is meromorphic at $\lambda=0$ 
for certain $\sigma$.

First we show a variant of Green's formula.
We need nice cut-off functions which exist by the following lemma.

\begin{lem}\label{lll}
There exists a cut-off function $\chi$ such that
\begin{enumerate}
\item $\chi> 0$ on a fundamental domain $F\subset X$,
\item $\supp(\chi)\subset \cup_{g\in L} gF$ for some finite subset
      $L\subset \Gamma$,
\item $\sum_{g\in\Gamma}g^*\chi = 1$,
\item $\sup_{k\in (\clo(F)\cap\Omega)M ,\: a\in A_+} a\:|\nabla^i\chi(ka)|<\infty$, $i\in\nat$,
\item the function $\chi_\infty$ on $K$ defined by $\chi_\infty(k):=\lim_{a\to\infty}\chi(ka)\chi_\infty(k)$ is smooth.
\end{enumerate}
\end{lem}
\proof
Let $W\subset \Omega$ be compact such that $\clo(F)\cap \Omega\subset W$.
Let $\psi\in C^\infty(\partial X)$ be a cut-off function such that
$\psi_{|\clo(F)\cap\Omega}=1$ and $\psi_{|\partial X\setminus W}=0$.
Let $B_R\subset X$, $R\in\R$ be the $R$-ball in $X$ centered at the origin
and choose $R>1$ so large that $F\subset B_R\cup WA_+$.
Let $\sigma\in C^\infty(A_+)$ be a cut-off function such that
$\sigma(r)=1$ for $r>1$ and $\sigma(r)=0$ for $r<1/2$. Finally let $\phi\in C_c^{\infty}(X)$
be a cut-off
function such that $\phi_{|B_R}=1$ and $\phi_{|\partial X\setminus B_{R+1}}=0$.
Then we set $\tilde{\chi}(ka):=\phi(ka) +  \psi(k)\sigma(a)$,
$k\in K$, $a\in A_+$, $ka\in X$.
If we define
$$\chi:=\frac{\tilde{\chi}}{\sum_{g\in \Gamma }g^*\tilde{\chi}}\ ,$$
then $\chi$ obviously satisfies (1),(2), (3), and (5).
It remains to verify (4).
Note that by construction of $\tilde{\chi}$ for all $l\in \nat$ there exists a constant $C<\infty$ such that for all $k\in K$, $a\in A_+$
$$|\nabla^l \tilde{\chi}(ka)| <  C a^{-1}\ .$$
Hence for any finite subset $L\subset \Gamma$  and 
$l\in\nat$ there exists a constant $C<\infty$ such that
for all $k\in  K$,
$g\in L$, $a\in A_+$
$$|\nabla^l g^\ast \tilde{\chi}(ka)| <  C a^{-1}\ .$$
This implies (4).
\hB

Now we consider the variant of Green's formula.
Let $\phi\in {}^\Gamma C^{-\infty}(\Lambda,V(\sigma_\lambda))$. If $\Ree(\lambda)\ge 0$, then
$res\circ \hat{J}_\lambda(\phi)\in C^{\infty}(B,V_B(\sigma_{-\lambda}))$ is well defined
even if $\hat{J}_\mu$ has a pole at $\mu=\lambda$. In the latter case the residue
of $\hat{J}_\mu$ at $\mu=\lambda$ is a differential operator $D_\lambda$ (see Lemma \ref{off})
and $res\circ D_\lambda(\phi)=0$.

\begin{prop}\label{green}
Assume that 
$\Ree(\lambda)\ge 0$ and  $\lambda\not=0$.
If
$\phi\in {}^\Gamma C^{-\infty}(\Lambda,V(\sigma_\lambda))$ and 
$f\in {}^\Gamma C^{-\infty}(\partial X,V(\tilde{\sigma}_\lambda))$
is  such that $f_{|\Omega}$ is  smooth,
 then $$\langle res\circ J_\lambda (\phi),res(f)\rangle =0\ .$$
\end{prop}
\proof
Let $(\gamma,V_\gamma) \in \hat{K}$ be a minimal $K$-type of the principal series representation of $G$
on $C^{\infty}(\partial X, V(\sigma_\lambda))$. Then there is an injective
 $T\in \Hom_M(V_\sigma,V_\gamma)$.
For $\Ree(\lambda)>0$ we define the endomorphism valued function
\begin{equation}\label{cgamma} c_{\gamma}(\lambda):=\int_{\bar{N}} a(\bar{n})^{-(\lambda+\rho)}
\gamma(\kappa(\bar{n}))d\bar{n}\in \End_M(V_\gamma)\ .\end{equation} 
The function $c_{\gamma}(\lambda)$ extends  meromorphically to all of $\aca$.
If $c_\gamma(\lambda)$ is singular, then $\lambda$ is integral. 
We choose $\tilde{T}\in \Hom_M(V_{\tilde{\sigma}},V_{\tilde{\gamma}})$
such that $\langle \gamma(w) T v, c_\gamma(\lambda) \tilde{T}u\rangle =\langle v,u\rangle$
for all $v\in V_\sigma$, $u\in V_{\tilde{\sigma}}$.
This is possible since $\gamma$ is a minimal $K$-type and hence $c_\gamma(\lambda)$ is invertible for $\Ree(\lambda)\ge 0$, $\lambda\not=0$.
 
Let $A=-\Omega_G+c_\sigma+\lambda^2$,
$P=P_\mu$ be the Poisson transform
(associated to $T$ or $\tilde{T}$,
respectively),
and $\chi$ be the cut-off function constructed
in Lemma \ref{lll}.
Note that $A=-\nabla^*\nabla + \cR$ for some selfadjoint endomorphism $\cR$ 
of $V(\gamma)$, where $-\nabla^*\nabla=\Delta$ is the Bochner Laplacian associated to the 
invariant connection $\nabla$ of $V(\gamma)$.
By $B_R$ we denote the metric $R$-ball centered at the origin of $X$.
The following is an application of Green's formula: 
\begin{eqnarray} 
0&=&\langle \chi A P\phi,Pf\rangle_{L^2(B_R)}- \langle \chi  P\phi,A Pf  \rangle_{L^2(B_R)}\label{limiz}\\
&=& \langle A \chi  P\phi,Pf\rangle_{L^2(B_R)}- \langle \chi  P\phi,A Pf  \rangle_{L^2(B_R)} - \langle [ A, \chi]  P\phi,Pf\rangle_{L^2(B_R)}\nonumber                           \\
&=& - \langle \nabla_n \chi  P\phi,Pf\rangle_{L^2(\partial B_R)}+ \langle \chi  P\phi,\nabla_n Pf  \rangle_{L^2(\partial B_R)} - \langle [ A, \chi]  P\phi,Pf\rangle_{L^2(B_R)}\ , \nonumber
\end{eqnarray}
where   $n$ is the exterior unit normal vector field at $\partial B_R$.

For the following discussion we distinguish between the two cases
$\Ree(\lambda)>0$ and $\Ree(\lambda)=0$, $\lambda\not=0$.

We first consider the case $\Ree(\lambda)>0$. 
We employ the following asymptotic behaviour of the Poisson transforms
of $f$ and $\psi$.
For $kM\in \Omega$ we have
\begin{eqnarray}
P(f)(ka)&=&a^{\lambda-\rho} c_\gamma(\lambda) \tilde{T} f(k)+O(a^{\lambda-\rho-\epsilon})\label{k0o}\\
P(\phi)(ka)&=&a^{-\lambda-\rho} \gamma(w)T (\hat{J}_\lambda \phi)(k) + O(a^{-\lambda-\rho-\epsilon})\label{k0o1}
\end{eqnarray}
where $\epsilon>0$.
While (\ref{k0o}) follows from the fact that $f_{|\Omega}$ is smooth,
(\ref{k0o1}) is shown in Lemma \ref{eee}.
The  estimate can be differentiated with respect to $a$ and holds locally uniformly on $\Omega$.
By property  (4) of $\chi$  we see that
$\langle [ A, \chi]  P\phi,Pf\rangle$
is integrable and by  property (3) and the $\Gamma$-invariance
of $f$ and $\phi$ we have
$\langle [ A, \chi]  P\phi,Pf\rangle_{L^2(X)}=0$.
Taking the limit $R\to\infty$ in (\ref{limiz}) we obtain
\begin{eqnarray*}
0&=&  (\lambda+\rho)\int_{\partial X} \chi_\infty(k) \langle\gamma(w)T (\hat{J}_\lambda \phi )(k), c_\gamma(\lambda)\tilde{T} f(k)\rangle  \\
&&+(\lambda-\rho)\int_{\partial X} \chi_\infty(k) \langle\gamma(w)T (\hat{J}_\lambda \phi )(k),c_\gamma(\lambda)\tilde{T}f(k)\rangle  \\
&=&2\lambda  \int_{\partial X} \chi_\infty(k) \langle  (\hat{J}_\lambda \phi )(k),f(k)\rangle \\
&=&2\lambda \langle res\circ\hat{J}_\lambda (\phi),res( f)\rangle  \ .
\end{eqnarray*}
This is the assertion of the proposition for $\Ree(\lambda)>0$.

Now we discuss the case $\Ree(\lambda)=0$ and $\lambda\not=0$.
In this case we have the following asymptotic behaviour
$$P(f)(ka) = a^{\lambda-\rho} c_\gamma(\lambda) \tilde{T} f(k)+ a^{-\lambda-\rho}\gamma(w) \tilde{T} \hat{J}_\lambda f(k) + O(a^{ -\rho-\epsilon})\ .$$
Instead of taking the limit $R\to\infty$ in (\ref{limiz}) we apply
$\lim_{r\to\infty}\frac{1}{r}\int_0^r dR$.
Again we have 
$$\lim_{r\to\infty}\frac{1}{r}\int_0^r \langle [ A, \chi]  P\phi,Pf\rangle_{L^2(B_R)} dR=0\ . $$
Moreover, the asymptotic term $a^{-\lambda-\rho} \tilde{T} \hat{J}_\lambda f(k)$
does not contribute to the limit because of   
$$-2\lambda    \lim_{r\to\infty}\frac{1}{r}\int_0^r R^{-2\lambda} \langle \chi(.R) \gamma(w)T \hat{J}_\lambda \phi  ,\gamma(w) \tilde{T} \hat{J}_\lambda f\rangle_{L^2(\partial X)}  dR=0\ .$$
The contribution of the term  $a^{\lambda-\rho} c_\gamma(\lambda) \tilde{T} f(k)$
leads to
$$0=2\lambda \langle res\circ\hat{J}_\lambda (\phi),res( f)\rangle$$
as in the case $\Ree(\lambda)>0$.
\hB
By the Harish-Chandra isomorphism characters of $\cZ$ are
parametrized by elements of $\haaa_\C^*/W$.
A character $\chi_\lambda$, $\lambda\in \haaa_\C^*$, is called integral,
if $$2\frac{\langle \lambda,\alpha\rangle}{\langle\alpha,\alpha\rangle}\in\Z$$
for all roots $\alpha$ of $(\gaaa,\haaa)$.
 
\begin{lem}\label{th43}
Let $\lambda\in \aca$ be such that 
\begin{itemize}
\item $\Ree(\lambda)\ge 0$ and $\chi_{\mu_\sigma+\rho_m-\lambda}$ is a non-integral character of $\cZ$ or
\item $\Ree(\lambda)<0$ and $\lambda$ is non-integral.
\end{itemize}
Let and $U\subset \partial X$ be open.
If  $\phi\in C^{-\infty}(\partial X, V(\sigma_\lambda))$ satisfies
$\phi_{|U}= 0$,
then $(\hat{J}_\lambda\phi)_{|U}= 0$ implies  $\phi=0$.
\end{lem}
Before turning to the proof note that Lemma \ref{off} implies that  $(\hat{J}_\lambda\phi)_{|U}$
is well-defined even if $\hat{J}_\lambda$ has a pole.\\[0.5cm]\noindent
\proof 
We modify an argument given by van den Ban-Schlichtkrull \cite{vandenbanschlichtkrull89}
for the case $\sigma=1$.    
If $\bar{N}$ is two-step nilpotent, then $\bar{\naaa}=\bar{\naaa}_1\oplus \bar{\naaa}_2$, $[\bar{\naaa}_1,\bar{\naaa}_1]=\bar{\naaa}_2$, where
$\bar{\naaa}_1$ corresponds to the negative of the shorter root $\alpha_1\in\aaaa^*_+$.

\begin{lem}\label{eee}
Let $\lambda\in\aca$, $\gamma$ be a representation of $K$, and 
$\phi\in C^{-\infty}(\partial X,V(\sigma_\lambda))$ satisfy $M \not\in \supp(\phi)$.
Then for $T\in \Hom(V_\sigma,V_\gamma)$, $k\in K\setminus\supp(\phi)M$  the Poisson transform
$P\phi:=P^T_\lambda\phi$ has an asymptotic expansion as $a\to\infty$ of the form
\begin{equation}\label{epan}
(P\phi)(ka)=a^{-(\lambda+\rho)}\gamma(w)T(\hat{J}_\lambda\phi)(k) 
 +\sum_{n\ge 1} a^{-(\lambda+\rho)-n\alpha_1}\psi_n(k)\ .
\end{equation}
Here $\psi_n$ are smooth functions on $K\setminus\supp(\phi)M$. 
The expansion converges uniformly for $a>>0$ and $k$ in compact subsets
of $K\setminus \supp(\phi)M$. 
\end{lem}
\proof
In the following computation we write the pairing of a distribution
with a smooth function as an integral. 
\begin{eqnarray}
(P\phi)(ka)&=&\int_K a(a^{-1}k^{-1}h)^{-(\lambda+\rho)}\gamma(\kappa(a^{-1}k^{-1}h)) T \phi(h) dh\nonumber\\
&=& \int_{\bar{N}} a(  a^{-1}w\kappa(\bar{n}))^{-(\lambda+\rho)}       \gamma(\kappa( a^{-1}w\kappa(\bar{n})))  T\phi(kw\kappa(\bar{n}))
a(\bar{n})^{-2\rho} d\bar{n} \nonumber\\
&=&\int_{\bar{N}} a(a  \bar{n}  a^{-1})^{-(\lambda+\rho)} a^{-(\lambda+\rho)} a(\bar{n})^{\lambda+\rho} \gamma(w) \gamma(\kappa(a\bar{n}a^{-1}))  T\phi(kw\kappa(\bar{n}))
a(\bar{n})^{-2\rho} d\bar{n} \nonumber\\
&=&a^{-(\lambda+\rho)}\gamma(w) \int_{\bar{N}} a(\bar{n})^{\lambda-\rho}  
a(a\bar{n}a^{-1})^{-(\lambda+\rho)} \gamma(\kappa(a\bar{n}a^{-1})) T \phi(kw\kappa(\bar{n})) d\bar{n}\label{tyh}\ .
\end{eqnarray}
For $z\in \R^+$ define $a_z\in A$ through $z=a_z^{-\alpha_1}$. We consider
the map
$\Phi:(0,\infty)\times \bar{N}\ni (z,\bar{n})\mapsto a_z\bar{n}a_z^{-1}\in \bar{N}$  which is can also be written as
$$\Phi(z,exp(X+Y)):= exp(zX+z^2Y), \quad X\in\naaa_1, Y\in\naaa_2\ .$$
Thus $\Phi$ and hence $ (z,\bar{n})\mapsto a(a_z\bar{n}a_z^{-1})^{-(\lambda+\rho)} \gamma(\kappa(a_z\bar{n}a_z^{-1}))$ extend analytically to  $\R\times\bar{N}$.
Taking the Taylor expansion with respect to $z$ at $z=0$ we obtain  
   $$a(a_z\bar{n}a_z^{-1})^{-(\lambda+\rho)} \gamma(\kappa(a_z\bar{n}a_z^{-1}))=\id + \sum_{n\ge 1}A_n(\bar{n}) z^n\ .$$
Here $A_n:\bar{N}\rightarrow \End(V_\gamma)$ are analytic and the series converges
in the spaces of smooth functions on $\bar{N}$ with values in $\End(V_\gamma)$.

Inserting this expansion into (\ref{tyh}) we obtain
$$(P\phi)(ka)=a^{-(\lambda+\rho)} \gamma(w) T (\hat{J}_\lambda \phi)(k) + \sum_{n\ge 1} a^{-(\lambda+\rho)-n\alpha_1} \psi_n(k)\ ,$$
where 
$$\psi_n(k):=\gamma(w)\int_{\bar{N}} A_n(\bar{n})  T a(\bar{n})^{\lambda-\rho }\phi(kw\kappa(\bar{n})) d\bar{n}\ .$$
Note that $k\mapsto \phi(kw\kappa(.))$ is a smooth family of distributions
with compact support in $\bar{N}$.
Thus $\psi_n$ is smooth. This finishes the proof of the lemma. \hB

We now continue the proof of Lemma \ref{th43}.
If $\Ree(\lambda)<0$, then we reduce to the case $\Ree(\lambda)>0$ replacing $\phi$ by $\hat{J}_\lambda(\phi)$. We can do this because $\lambda$ is then non-integral and $\hat{J}_\lambda$  is regular and bijective. 
Note that in this case $\chi_{\mu_\sigma+\rho_m-\lambda}$ is non-integral.
Thus assume that $\Ree(\lambda)\ge 0$.

We choose the representation $\gamma\in \hat{K}$ and $T\in \Hom_M(V_\sigma,V_\gamma)$ such that $P:=P^T_\lambda$ is injective.
The range of $P$ can be identified with the kernel of a certain invariant
differential operator $D:C^{\infty}(X,V(\gamma))\rightarrow C^{\infty}(X,V(\gamma^\prime))$ 
for a suitable representation $\gamma^\prime$ of $K$ (see \cite{bunkeolbrich947}, Sec.3). In particular, $P\phi$ is real-analytic.

We now assume $\phi\not\equiv 0$.
Moreover, without loss of generality we can assume that $ M\in U$.
Since $P\phi$ is real analytic, the expansion (\ref{epan}) does not vanish.
Let $m$ be the smallest integer such that $\psi_m\not\equiv 0$ near $M$
(where $\psi_0:=\gamma(w)T\hat{J}_\lambda \phi)$.
We argue that $m=0$ and thus obtain a contradiction.

To prove that $m=0$ we again argue by contradiction. Assume that $m>0$.
 With respect to the coordinates $k,a$ the operator $D$ has the form
$D=D_0+a^{-\alpha_1}R(a,k)$, where $D$ is a constant coefficient
operator on $A$ and $R$ remains bounded if $a\to\infty$ (see \cite{warner721}, Thm. 9.1.2.4).
Moreover, it is known that $D_0$ coincides with the $\bar{N}$-radial part of $D$.

Choose  $k\in K$  near $1$ and $\sigma^\prime\subset \gamma_{|M}$ such that
that there exists an orthogonal projection $S\in \Hom_M(V_\gamma,V_{\sigma^\prime})$ with 
 $S \gamma(k)\psi_m(k) =:v \not= 0$.
Consider the $\bar{N}$-invariant section $f\in C^{\infty}(X,V(\gamma))$
defined by
$$f(\bar{n}a):=a^{-(\lambda+\rho+m\alpha_1)} S^* v\ .$$
Since $D$ annihilates the asymptotic expansion (\ref{epan}), one can check
that $Df=D_0f=0$ and thus $f=P\phi_1$ for some $\bar{N}$-invariant $\phi_1\in C^{-\infty}(\partial X,V(\sigma_\lambda))$.

Now $f=P^S_{\lambda+m\alpha_1}\delta v$, where $\delta v\in C^{-\infty}(\partial X,V(\sigma^\prime_{\lambda+m\alpha_1}))$
is the delta distribution at $1$ with vector part $v$. 
Since $D$ and $P^S_{\lambda+m\alpha_1}$ are $G$-equivariant and $\delta v$ generates the $G$-module $C^{-\infty}(\partial X,V(\sigma^\prime_{\lambda+m\alpha}))$, we obtain
a non-trivial
intertwining operator $I$ from $C^{-\infty}(\partial X,V(\sigma^\prime_{\lambda+m\alpha}))$ 
to the kernel of $D$, hence to $C^{-\infty}(\partial X,V(\sigma_\lambda))$.
Thus principal series representations $\pi^{\sigma,\lambda}$ and $\pi^{\sigma^\prime,\lambda+m\alpha_1}$ have the same infinitesimal character
$\chi_{\mu_\sigma+\rho_m-\lambda}$. Since by assumption this character is non-integral both principal series are irreducible (see \cite{collingwood85}, 4.3.3).
Hence $I$ is an isomorphism. Because of $m\not=0$ 
we have $\sigma\not=\sigma^\prime$. This implies that  $\pi^{\sigma,\lambda}$
and $\pi^{\sigma^\prime,\lambda+m\alpha_1}$ can not be isomorphic.
This is the contradiction we aimed at.
We conclude that $m=0$.\hB

The above argument can be extended to cover some cases
of $\lambda\in\aaaa^*$, $\Ree(\lambda)\ge 0$, with $\chi_{\mu_\sigma+\rho_m-\lambda}$ integral. 
This would lead to stronger vanishing results 
for ${}^\Gamma C^{-\infty}(\Lambda,V(\sigma_\lambda))$ below.
For example, if $\sigma$ is the trivial $M$-type, then the
assertion of Lemma \ref{th43} holds true for all $\lambda\in \aca$
with $\Ree(\lambda)\ge 0$.
But there exist examples of $\sigma\in\hat{M}$
and $\lambda\in\aaaa^*$ with $\Ree(\lambda)\ge 0$ and 
$\chi_{\mu_\sigma+\rho_m-\lambda}$ integral where
the assertion of Lemma \ref{th43} is false.
The possible failure of the lemma at $\lambda=0$ is connected
with the existence of $L^2(Y,V_Y(\gamma))_{scat}$.

\begin{kor}\label{ujn}
Assume that $\lambda\in \aca$, that $\Ree(\lambda)\ge 0$, and that 
$\chi_{\mu_\sigma+\rho_m-\lambda}$ is non-integral. 
 Then
$$res_\Omega\circ \hat{J}_\lambda:C^{-\infty}(\Lambda,V(\sigma_\lambda))\rightarrow C^{\infty}(\Omega,V(\sigma_{-\lambda}))$$
is injective.
\end{kor}

\begin{lem}\label{ghu}
Assume that $\lambda\in \aca$, that $\Ree(\lambda)\ge 0$, that
$\chi_{\mu_\sigma+\rho_m-\lambda}$ is non-integral,
and that
 $$ext:C^{\infty}(B,V_B(\sigma_\lambda))\rightarrow {}^\Gamma C^{-\infty}(\partial X,V(\sigma_\lambda))$$
is regular. Then ${}^\Gamma C^{-\infty}(\Lambda,V(\sigma_\lambda))=0$ .
\end{lem}
\proof
The assumption of the lemma implies that
$$ext:C^{\infty}(B,V_B(\tilde{\sigma}_\lambda))\rightarrow {}^\Gamma C^{-\infty}(\partial X,V(\tilde{\sigma}_\lambda))$$
is regular. 
Indeed, if this extension is singular, then $\lambda$ is real and
$$ext:C^{\infty}(B,V_B(\sigma_\lambda))\rightarrow {}^\Gamma C^{-\infty}(\partial X,V(\sigma_\lambda))$$
is singular since there is a conjugate linear isomorphism $\sigma\cong\tilde{\sigma}$.

Let $\phi\in{}^\Gamma C^{-\infty}(\Lambda,V(\sigma_\lambda))$.
Then for all $f\in C^{\infty}(B,V_B(\tilde{\sigma}_\lambda))$ we have
by Lemma \ref{green} 
$$0=\langle res\circ \hat{J}_\lambda(\phi),res\circ ext(f)\rangle= \langle res\circ \hat{J}_\lambda(\phi), f \rangle \ .$$
Thus $res\circ \hat{J}_\lambda(\phi)=0$ and by Corollary \ref{ujn} we conclude
$\phi=0$.
\hB

\begin{lem}\label{extregatim}
If $\Ree(\lambda)= 0$ and $\Imm(\lambda)\not=0$, then 
$$ext:C^{-\infty}(B,V_B(\sigma_\lambda))\rightarrow {}^\Gamma C^{-\infty}(\partial X,V(\sigma_\lambda))$$
is regular.
\end{lem}
\proof
Note that $\chi_{\mu_\sigma+\rho_m-\lambda}$ is a non-integral
character of $\cZ$. Since $\lambda\in\aca(\sigma)$, the extension
is meromorphic at $\lambda$. Assume that $ext$ has a pole.
The singular part of $ext$ maps
to distributions which are supported on the limit set $\Lambda$.
Then by Corollary
\ref{ujn} the scattering matrix 
$S_\lambda=res\circ J_\lambda\circ ext$ has a pole at $\lambda$, too.
But this contradicts Lemma \ref{unitary}.
\hB

In the following proposition we formulate
the most complete vanishing result for\linebreak[4] ${}^\Gamma C^{-\infty}(\Lambda,V(\sigma_\lambda))$, $\Ree(\lambda)\ge 0$,
which is stated in the present paper.
In general, it is not optimal for integral infinitesimal characters.
We give upper bounds depending on $\delta_\Gamma$ for the parameters $\lambda$  with ${}^\Gamma C^{-\infty}(\Lambda,V(\sigma_\lambda))\not=0$.
One can then deduce corresponding bounds for the discrete spectrum
of $\cZ$ on $L^2(Y,V_Y(\gamma))$. If one is only interested in
the finiteness of the discrete spectrum, then it is sufficient to
know the rough bound $\lambda \le \rho$, which is implied by the classification
of the unitary representations of $G$ (see the proof of Corollary \ref{nahenull}).

\begin{prop}\label{upperbound}
Let $\Ree(\lambda)\ge 0$. Then any of the following conditions
implies that $${}^\Gamma C^{-\infty}(\Lambda,V(\sigma_\lambda))=0\ .$$
\begin{enumerate}
\item $\Imm(\lambda)\not=0$,
\item $\Ree(\lambda)>\delta_\Gamma$ and $\chi_{\mu_\sigma+\rho_m-\lambda}$ is
a non-integral character of $\cZ$,
\item $$\hspace{-2cm}\Ree(\lambda)>\frac{\rho^2+\rho\delta_\Gamma}{3\rho-\delta_\Gamma}$$
\end{enumerate} 
\end{prop}
\proof
We consider $1$. If $\Ree(\lambda)>0$, then the assertion follows
from Lemma \ref{pokm}. If $\Ree(\lambda)=0$, then we apply Lemmas
\ref{ghu} and \ref{extregatim}.

Sufficiency of condition $2.$  follows from Lemma \ref{ghu} and
\ref{defofext}.

Condition $3.$ is important for $\lambda$ with
$\chi_{\mu_\sigma+\rho_m-\lambda}$ integral. In this
case the relation of the space ${}^\Gamma C^{-\infty}(\Lambda,V(\sigma_\lambda))$ with the singularities
of $ext$ is quite unclear.
Let $\gamma$ denote a minimal $K$-type of the principal
series representation $\pi^{\sigma,\lambda}$ of $G$ on $C^\infty(\partial X,V(\sigma_\lambda))$,
and let $T\in \Hom_M(V_\sigma,V_\gamma)$ be an isometric
embedding. Then the Poisson transform 
$$P^T_\lambda :C^{-\infty}(\partial X,V(\sigma_\lambda))\rightarrow C^\infty(X,V(\gamma))$$ is injective.
Let $\psi\in {}^\Gamma C^{-\infty}(\Lambda,V(\sigma_\lambda))$
and $f:=P^T_\lambda\psi\in C^\infty(Y,V_Y(\gamma))$.
Since $\Ree(\lambda)>0$ we have for any test function $\phi\in C^\infty(\partial X,V(\sigma_{\bar{\lambda}}))$
$$\bar{c}_\sigma(\lambda)(\psi,\phi)=\lim_{a\to\infty} a^{\rho-\bar{\lambda}}\int_K (f(ka), T \phi(k)) dk \ .$$
We will show that condition $3.$ implies that $(\psi,\phi)=0$.
Since $\phi$ was arbitrary, it follows $\psi=0$.

In the following argument we assume that $\Ree(\lambda)<\rho$.
An easy modification works for $\Ree(\lambda)\ge\rho$. 
The asymptotic expansion of $f$ near $\Omega$ obtained
in Lemma \ref{eee} implies that $f\in L^p(Y,V_Y(\gamma))$ for
$p=\frac{2\rho}{\rho+\Ree(\lambda)}+\epsilon$, $\forall \epsilon>0$.  
If we set $\tilde{f}(ka)=a^{-(\rho+\delta_\Gamma+\epsilon)/p}f(ka)$,
then $\tilde{f}\in L^p(X,V(\gamma))$. 
In fact we can estimate
\begin{eqnarray*}
\|\tilde{f}\|_{L^p(X,V(\gamma))}&=&\sum_{g\in\Gamma} \int_{gF} |\tilde{f}(g)|^p dg\\
&\le & \sum_{g\in \Gamma} \min\{a_h\:|\: h\in gF\}^{-(\rho+\delta_\Gamma+\epsilon)}
 \int_{gF} |f(g)|^p dg\\
&\le & C \sum_{g\in \Gamma} a_g^{-(\rho+\delta_\Gamma+\epsilon)}\\
&<& \infty \ .
\end{eqnarray*}\
Let $\chi\in C^\infty_c(0,1)$ be such that $\int_0^1\chi(t) dt=1$.
We extend $\chi$ to $\R$ by zero.
Then we define $\chi_n\in C_c^\infty(A_+)$ by
$\chi_n(a):=a^{-\rho-\bar{\lambda}}\chi(|\log(a)|-n)$.
If we set $\phi_n(ka)= T\chi_n(a)\phi(k)$,
then we can write 
$$\bar{c}_\sigma(\lambda)(\psi,\phi )
=\lim_{n\to \infty} (f,\phi_n) \ .$$

Now let $\tilde{\phi}_n(ka):= a^{(\rho+\delta_\Gamma+\epsilon)/p} \phi_n(ka)$.
Then
\begin{equation}\label{thishere} \bar{c}_\sigma(\lambda)(\psi,\phi)
=\lim_{n\to \infty} (\tilde{f},\tilde{\phi}_n) \ .\end{equation}
We consider $\Psi(ka):=a^{(\rho+\delta_\Gamma+\epsilon)/p}a^{-\rho-\bar{\lambda}}T \phi(k)$
and let $q\in (1,\infty)$ be the dual exponent to $p$.
Note that $\bar{c}_\sigma(\lambda)\not=0$. 
If we assume that 
$\Psi\in L^q(X,V(\gamma))$, then $(\psi,\phi)=0$ follows from
(\ref{thishere}).
The following inequality implies that $\Psi\in L^q(X,V(\gamma))$:
\begin{equation}\label{ineq}q\left(\frac{\rho+\delta_\Gamma+\epsilon}{p}-(\rho+\Ree(\lambda))\right)<-2\rho \ .\end{equation} 
A simple computation shows that (\ref{ineq}) indeed follows from $3$.
\hB 

If $\chi_{\mu_\sigma+\rho_m-\lambda}$ is integral, then 
one can combine the argument of Lemma \ref{th43}, or
the understanding of its failure, respectively, with the
above argument in order to choose a better exponent $p$.
This leads to a stronger vanishing result.
We omit this rather involved discussion,
partly because we do not know how to formulate a general result.
This omission causes some loss of information about the
discrete spectrum of $\cZ$ on $L^2(Y,V_Y(\gamma))$.

In general, we do not know whether $ext$ is meromorphic at $\lambda=0$.
But we can show the following result.

\begin{prop}\label{mystic}
If $\hat J_\lambda$ has a pole at $\lambda=0$, then $ext$ is meromorphic in
a neighbourhood of $0$.
\end{prop}

The proof of the proposition will occupy the remainder of this section.
First we want to fix an important corollary which applies to all $M$-types
$\sigma$.  

\begin{kor}\label{nahenull}
There exists $\epsilon>0$ such that
$ext$ is regular and ${}^\Gamma C^{-\infty}(\Lambda,V(\sigma_\lambda))=0$ for $\lambda\in (0,\epsilon)$.
\end{kor}
\proof
If contrary to the assertion there is a
sequence $\lambda_\alpha>0$ of poles of $ext$ with $\lambda_\alpha\to 0$,
then the residues of $ext$ give elements of ${}^\Gamma C^{-\infty}(\Lambda,V(\sigma_{\lambda_\alpha}))$. Applying to these elements a suitable Poisson
transform $P^T_{\lambda_\alpha}$, $T\in \Hom_M(V_\sigma,V_\gamma)$, we obtain smooth sections in $L^2(Y,V_Y(\gamma))$. Considered as
elements of $L^2(\Gamma\backslash G)\otimes V_\gamma$ they generate unitary
representations of $G$ whose underlying Harish-Chandra modules are dual
to the underlying Harish-Chandra modules of $C^{\infty}(\partial X,V(\sigma_{\lambda_\alpha}))$, at least for $\lambda_\alpha$ sufficiently small. 
Here we have used that principal series representations with non-integral infinitesimal character are irreducible. Since the dual of a unitary representation is
unitary too, we conclude that the principal series representations 
$C^{\infty}(\partial X,V(\sigma_{\lambda_\alpha}))$ carry invariant 
Hermitian scalar products or, in representation theoretic language, belong
to the complementary series. But, by a result of Knapp-Stein \cite{knappstein71}, Par. 14, there is a complementary series for $\sigma$ iff
$P_\sigma(0)=0$ or, equivalently, iff $\hat J_\lambda$ has a pole at $\lambda=0$. But in this case $ext$ is meromorphic at $\lambda=0$ by Proposition
\ref{mystic}. This is in contradiction with the existence of the sequence $\lambda_\alpha$ of poles of $ext$ with $\lambda_\alpha\to 0$.
\hB   

We now turn to the proof of Proposition \ref{mystic}. If $\delta_\Gamma<0$, there
is nothing to show. For $\delta_\Gamma\ge 0$ we are going to use the embedding
trick as in the proof of Proposition \ref{part1}. 

Let $G=G^n\subset G^{n+1}\subset\dots$ be the corresponding series of groups.
Choose $k\in\nat$ such that $\delta_\Gamma^{n+k}=\delta_\Gamma^{n}-\rho^{n+k}+\rho^n<0$. If there is a representation $\sigma^{n+k}$ of $M^{n+k}$ such that $\sigma^n:=\sigma\subset \sigma^{n+k}_{\ |M^n}$ and $-\rho^{n+k}+\rho^n \in \aca (\sigma^{n+k})$,
then $ext^{n+k}$ is meromorphic at $\lambda=-\rho^{n+k}+\rho^n$, which in
turn implies the meromorphy of $ext^n$ at $\lambda=0$.

Recall the definition (\ref{acagood}) of $\aca (\sigma^{n+k})$.
Since $\hat J_\lambda$ has at most first-order poles 
$$P_{\sigma^{n+k}}(-\rho^{n+k}+\rho^n)=0$$ implies that $-\rho^{n+k}+\rho^n \in \aca (\sigma^{n+k})$. It follows from the functional equation (\ref{spex}) that the intertwining operator $\hat J_\lambda$
has a pole at $\lambda=0$ iff 
$P_{\sigma^n}(0)=0$. In this case $\sigma^n$ is irreducible and 
Weyl-invariant.
So the discussion above shows that the
following lemma implies Proposition \ref{mystic}.

\begin{lem}\label{casebycase}
Let $G^n$ one of the following four series of groups :
\begin{eqnarray*}
 Spin(1,2n)&&n\ge 1\ ,\\ 
Spin(1,2n+1)&&n\ge 1\ ,\\
SU(1,n)&&n\ge 2\ ,\\
Sp(1,n)&&n\ge 2\ .
\end{eqnarray*}
If $\sigma^n \in \hat M^n$ (i.e. $\sigma^n$ is irreducible) is
Weyl invariant and satisfies 
$P_{\sigma^n}(0)=0$, then for any $k\in\nat$
there exists $\sigma^{n+k} \in \hat M^{n+k}$ with $\sigma^n\subset \sigma^{n+k}_{\ |M^n}$ and  $P_{\sigma^{n+k}}(-\rho^{n+k}+\rho^n)=0$.
\end{lem} 
\proof
First we need explicit expressions of the Plancherel densities $P_{\sigma^n}$ as
given e.g. in \cite {knapp86}, Prop. 14.26. Consider the Cartan subalgebra
$\haaa_\C^n=\aaaa_\C\oplus\taaa_\C^n$, $\taaa^n$ being a Cartan
subalgebra of $\maaa^n$, and the subset $\Delta_r^n$ of the roots $\Delta^n$
of $\haaa_\C^n$ in $\gaaa^n_\C$ given by
$$ \Delta_r^n:=\{\alpha \in \Delta^n\:|\:\alpha_{|\aaaa} \mbox{ is a root of
$\aaaa$ in }\naaa^n\}\ .$$
For $G^n\not=Spin(1,2n+1)$ there is exactly on real root $\beta_r\in \Delta_r^n$
distinguished by $\beta_{r\:|\taaa^n}=0$. Furthermore, in this case we consider
a special element $m^n$ of the center of $M^n$. For $G=Spin(1,2n)$ the element $m^n$ is the non-trivial element in the kernel of the projection $Spin(1,2n)\rightarrow SO_0(1,2n)$. For $G^n=SU(1,n)$ $(Sp(1,n))$ we use the standard representation
in $Gl(n+1,\C)$ $(Gl(n+1,\bf H))$ in order to fix $m^n$ by
$$ m^n:= \left(
\begin{array}{rrc}
-1&0&0\\
0&-1&0\\
0&0&\id_{n-1}
\end{array}
\right)\ . $$
Since $(m^n)^2=1$ we have $\sigma^n(m^n)=\pm \id$. The embedding
$M^n\hookrightarrow M^{n+k}$ sends $m^n$ to $m^{n+k}$. 

There are nontrivial constants $C(\sigma^n)$ such that
\begin{equation}\label{blanche}
P_{\sigma^n}(\lambda)=C(\sigma^n)f_{\sigma^n}(\lambda)\prod_{\beta\in \Delta_r^n} \langle \mu_{\sigma^n}+\rho_m^n-\lambda,\beta\rangle\ ,
\end{equation}
where $\langle.,.\rangle$ is a Weyl-invariant bilinear scalar product on 
$(\haaa_\C^n)^*$ and
$$ f_{\sigma^n}(\lambda)=\left\{
\begin{array}{ccl}
1&G^n=Spin(1,2n+1)&\\
\tan \left(\pi \frac{\langle\lambda,\beta_r\rangle}{\langle\beta_r,\beta_r\rangle}\right), &G^n\not=Spin(1,2n+1),&
\sigma^n(m^n)=-\ee^{2\pi\imath\frac{\langle\rho^n,\beta_r\rangle}{\langle\beta_r,\beta_r \rangle}}\id\\
\cot \left(\pi \frac{\langle\lambda,\beta_r\rangle}{\langle\beta_r,\beta_r\rangle}\right), &G^n\not=Spin(1,2n+1),&
\sigma^n(m^n)=\ee^{2\pi\imath\frac{\langle\rho^n,\beta_r\rangle}{\langle\beta_r,\beta_r \rangle}}\id
\end{array} 
\right. \ . $$
According to these three possibilities of the form of the Plancherel density we
distinguish between the odd-dimensional, the tan and the cot case. We shall
construct the representation $\sigma^{n+k}$ case by case. The easiest one is 

\noindent
{\bf The tan case}\newline
For any representation $\sigma^n$ fitting into this case we have
$$P_{\sigma^n}(0)=\tan(0)=0\ .$$
Let $\sigma^{n+k}\in \hat M^{n+k}$ be an arbitrary representation satisfying
$\sigma^n\subset \sigma^{n+k}_{\ |M^n}$. Then $\sigma^{n+k}(m^{n+k})=-\id$
iff $\sigma^n(m^n)=-\id$. Thus $\sigma^{n+k}$ belongs to the $\tan$-case iff
$\rho^{n+k}-\rho^n$ is an integer multiple of $\beta_r$ (this is always the case
except when $G^n=SU(1,n)$ and $k$ is odd). We conclude that
$$ P_{\sigma^{n+k}}(-\rho^{n+k}+\rho^n)=\left\{
\begin{array}{c}
C\tan(\pi\frac{\langle\-\rho^{n+k}+\rho^n,\beta_r\rangle}{\langle\beta_r,\beta_r \rangle})\\ 
C^\prime\cot(\pi\frac{\langle\-\rho^{n+k}+\rho^n,\beta_r\rangle}{\langle\beta_r,\beta_r \rangle})
\end{array}\right\}
=0\ .$$

\noindent
{\bf The odd-dimensional case}\newline
As a Cartan subalgebra of $\gaaa^n$ we choose
$$\haaa^n:= \left\{ T_\nu:= \left(
{\scriptsize 
\begin{array}{cccc}
\begin{array}{cc}
0&\nu_0\\ \nu_0&0
\end{array}&&&\\
&\begin{array}{cc}
0&-\nu_1\\ \nu_1&0
\end{array}&&\\
&&\ddots&\\
&&&\begin{array}{cc}
0&-\nu_n\\ \nu_n&0
\end{array}
\end{array}
}
\right)\:\Bigg|\:\nu_i\in\R \right\}\ \ ,$$
where $\aaaa=\{T_\nu\:|\:\nu_i=0,\:i=1,\dots,n\}$ and $\taaa^n=\{T_\nu\:|\:\nu_0=0\}$.
Define $e_i\in (\haaa_\C^n)^*$ by $e_0(T_\nu):=\nu_0$ and 
$e_i(T_\nu):=\imath \nu_i$,
$i=1,\dots,n$. We normalize $\langle.,.\rangle$ such that $\{e_i\}$ becomes an
orthonormal basis of $(\haaa_\C^n)^*$. Sometimes we will write elements of $\imath (\taaa^n)^*$ as $n$-tuples of reals with respect to the basis $e_1,\dots,e_n$. We choose positive roots of $\taaa^n$ in $\maaa^n_\C$ as follows:
$$\Delta_m^{n,+}:=\{e_i\pm e_j \:|\: 1\le i<j\le n\}\ .$$
Then $\rho^n_m=(n-1,n-2,\dots,1,0)$ and the irreducible representations of
$M^n$ correspond to the highest weights
$$ \hat M^n \cong \{\mu_{\sigma^n}=(m_1,\dots,m_n)\:|\:m_1\ge\dots\ge m_{n-1}\ge |m_n|,\: m_i-m_j\in\Z,\: m_n\in\frac{1}{2}\Z\}\ .$$
Furthermore
$$\Delta^n_r=\{e_0\pm e_i\:|\:i=1,\dots,n\}\ .$$
We obtain $\rho^n=ne_0$ and for $\lambda=ze_0$
$$P_{\sigma^n}(\lambda)= C(\sigma^n)\prod_{i=1}^{n}(z^2-(m_i+n-i)^2)\ .$$
We see that $P_{\sigma^n}(0)=0$ iff $m_n=0$. In this case set
$$\mu_{\sigma^{n+k}}:=(m_1,\dots,m_{n-1},m_n=0,\underbrace{0,\dots,0}
_{\mbox{\scriptsize $k$ times}})\ .$$
Then $\sigma^n\subset \sigma^{n+k}_{\ |M^n}$ and the $n$-th factor of $P_{\sigma^{n+k}}$ is given by $z^2-(m_n+n+k-n)^2 = z^2- k^2$. Thus
$P_{\sigma^{n+k}}(-\rho^{n+k}+\rho^n)=P_{\sigma^{n+k}}(-ke_0)=0$.

\noindent
{\bf The cot case}\\
First we observe that if $G=Spin(1,2n)$ and $\sigma^n$ belongs to the $\cot$-case, then $P_{\sigma^n}(0)\not=0$. In fact, in this case $\sigma^n$ is a faithful representation of $M^n=Spin(2n-1)$ and the only root $\beta\in\Delta^n_r$
perpendicular to $\mu_{\sigma^n}+\rho_m^n$ is $\beta_r$. We leave the simple verification to the reader. Since $\cot$ has a pole at $0$ the observation follows
from (\ref{blanche}). An alternative proof is given in \cite{knappstein71}, Prop. 55.
We are left with the discussion of $G^n=SU(1,n)$ and $G^n=Sp(1,n)$.

We start with $G^n=SU(1,n)$.
The group $M^n$ has the form
$$M^n= \left\{\left(\begin{array}{ccc}
z&0&0\\0&z&0\\0&0&B \end{array}\right)\:|\:\:
z\in U(1),B\in U(n-1), z^2\det B=1\right\}\ .$$
We consider the Cartan subalgebra of ${\bf u}(1,n)$
$$\tilde\haaa^n:= 
\left\{ T_\nu:= \left(
{\scriptsize 
\begin{array}{ccccc}
\imath \nu_1&\nu_0&&&\\ 
\nu_0&\imath \nu_1&&&\\
&&\imath\nu_2&&\\
&&&\ddots&\\
&&&&\imath\nu_n
\end{array}
}
\right)\:\Bigg|\:\nu_i\in\R \right\}\ \ ,$$
their subalgebras
$\aaaa:=\{T_\nu\:|\:\nu_i=0,\:i=1,\dots,n\}$ and $\tilde\taaa^n:=\{T_\nu\:|\:\nu_0=0\}$.
Then $\haaa^n:=\{T_\nu\in\tilde\haaa^n\:|\:2\nu_1+\sum_{i=2}^n\nu_i=0\}$ and
$\taaa^n:=\{T_\nu\in\tilde\taaa^n\:|\:2\nu_1+\sum_{i=2}^n\nu_i=0\}$ are Cartan subalgebras 
of $\gaaa^n$ and $\maaa^n$, respectively.
We define elements $\alpha,\beta,e_i\in (\tilde\haaa_\C^n)^*$, $i=2,\dots,n$, by $\alpha(T_\nu):=\nu_0$, $\beta(T_\nu)=\imath\nu_1$ and 
$e_i(T_\nu):=\imath \nu_i$. We extend $\langle.,.\rangle$ to $(\tilde\haaa^n_\C)^*$ such that $\{\alpha,\beta,e_2,\dots,e_n\}$ becomes an orthogonal basis of $(\tilde\haaa_\C^n)^*$ with
$|\alpha|^2=|\beta|^2=1$ and $|e_i|^2=2$.
Sometimes we will write elements of $\imath (\tilde\taaa^n)^*$ as $n$-tuples of reals with respect to the basis $\beta,e_2,\dots,e_n$. We choose the positive roots of $\taaa^n$ in $\maaa^n_\C$ as follows:
$$\Delta_m^{n,+}:=\{e_i-e_j \:|\: 2\le i<j\le n\}\ .$$
Then $\rho^n_m=\frac{1}{2}(0,n-2,n-4,\dots,4-n,2-n)$. 
Furthermore we have
$$\Delta^n_r=\{2\alpha,\alpha\pm(\beta-e_i)\:|\:i=2,\dots,n\}$$
and $\rho^n=n\alpha$.
We represent highest weights of representations of $M^n$ which are functionals
on $\taaa^n$ by elements of $(\tilde\taaa^n_\C)^*$: 
\begin{eqnarray*} 
\hat M^n &\cong& \{\mu_{\sigma^n}=(m_1,\dots,m_n)\:|\:m_2\ge\dots
\ge m_n,\: m_i
\in\Z \}\\
&&\quad\mbox{ modulo translation by elements of the form } 
(2\nu,\nu,\dots,\nu)\  .
\end{eqnarray*}
Then we can compute the scalar products with elements of $\Delta^n_r$ inside $(\tilde\haaa^n_\C)^*$,
and the result will not depend on the chosen representative.
Since $\rho^n=n\alpha$ and $\beta_r=2\alpha$ we see that $\sigma^n$ belongs to the $\cot$-case iff $m_1\equiv n\:(2)$.
We obtain for $\lambda=z\alpha$
$$P_{\sigma^n}(\lambda)= C(\sigma^n)\cot(\frac{\pi}{2}z)2z\prod_{i=2}^{n}(z^2-(m_1-n-2(m_i+1-i))^2)\ .$$
We see that $P_{\sigma^n}(0)=0$ iff for one index $i_0\in\{2,\dots,n\}$ the
following equation holds:
$$\frac{m_1-n}{2}=m_{i_0}+1-i_0\ .$$ 
In this case set
$$\mu_{\sigma^{n+k}}:=(m_1,\dots,m_{i_0},\underbrace{m_{i_0},\dots,m_{i_0}}
_{\mbox{\scriptsize $k$ times}},m_{i_0+1},\dots,m_n)\ .$$
Then $\sigma^n\subset \sigma^{n+k}_{\ |M^n}$.
Depending on the parity of $k$ the representation $\sigma^{n+k}$ belongs to the 
$\tan$-case or $\cot$-case, respectively. In any case, the function $f_{\sigma^{n+k}}$ has a first-order
pole at $\lambda=-\rho^{n+k}+\rho^n=-k\alpha$. But
in addition to the $i_0$-th factor $z^2-(m_1-n-k-2(m_{i_0}+1-i_0))^2$ of the polynomial part of $P_{\sigma^{n+k}}(z\alpha)$ also the $(i_0+k)$-th factor $z^2-(m_1-n-k-2(m_{i_0}+1-i_0-k))^2$ is equal to $z^2- k^2$. Thus
the polynomial part contributes a second order zero at $z=k$ and
$P_{\sigma^{n+k}}(-\rho^{n+k}+\rho^n)=P_{\sigma^{n+k}}(-k\alpha)=0$.

\begin{itemize}\item $G^n=Sp(1,n)$\end{itemize}
The group $M^n$ has the form
$$M^n= \left\{\left(\begin{array}{ccc}
q&0&0\\0&q&0\\0&0&B \end{array}\right)\:|\:\:
q\in Sp(1),B\in Sp(n-1)\right\}\ .$$
As a  Cartan subalgebra of $\gaaa^n$ we choose
$$\haaa^n:= 
\left\{ T_\nu:= \left(
{\scriptsize 
\begin{array}{ccccc}
\imath \nu_1&\nu_0&&&\\ 
\nu_0&\imath \nu_1&&&\\
&&\imath\nu_2&&\\
&&&\ddots&\\
&&&&\imath\nu_n
\end{array}
}
\right)\:\Bigg|\:\nu_i\in\R \right\}\ \ ,$$
where
$\aaaa=\{T_\nu\:|\:\nu_i=0,\:i=1,\dots,n\}$ and $\taaa^n=\{T_\nu\:|\:\nu_0=0\}$.
We define elements $\alpha,\beta,e_i\in (\haaa_\C^n)^*$, $i=2,\dots,n$, by $\alpha(T_\nu):=\nu_0$, $\beta(T_\nu)=\imath\nu_1$ and 
$e_i(T_\nu):=\imath \nu_i$. Then $\{\alpha,\beta,e_2,\dots,e_n\}$ becomes an orthogonal basis of $(\tilde\haaa_\C^n)^*$ and we normalize $\langle.,.\rangle$ such that
$|\alpha|^2=|\beta|^2=1$ and $|e_i|^2=2$.
Sometimes we will write elements of $\imath (\taaa^n)^*$ as $n$-tuples of reals with respect to the basis $\beta,e_2,\dots,e_n$. We choose the positive roots of $\taaa^n$ in $\maaa^n_\C$ as follows
$$\Delta_m^{n,+}:=\{2\beta,e_i\pm e_j, 2e_i \:|\: 2\le i<j\le n\}\ .$$
Then $\rho^n_m=(1,n-1,n-2,\dots,2,1)$ and the irreducible representations of
$M^n$ correspond to the highest weights
$$ \hat M^n \cong \{\mu_{\sigma^n}=(m_1,\dots,m_n)\:|\:m_1\ge 0,m_2\ge\dots\ge m_{n}\ge 0,\: m_i\in\Z\}\ .$$
Furthermore we have
$$\Delta^n_r=\{2\alpha,2(\alpha\pm\beta),\alpha\pm(\beta\pm e_i)\:|\:i=2,\dots,n\}$$
and $\rho^n=(2n+1)\alpha$.
Since $\beta_r=2\alpha$ we see that $\sigma^n$ goes with cot iff $m_1$ is odd.
We obtain for $\lambda=z\alpha$
\begin{eqnarray*}
P_{\sigma^n}(\lambda)&=& C(\sigma^n)\cot(\frac{\pi}{2}z)\:2z\:4(z^2-(m_1+1)^2)\\
&&\quad 
\prod_{i=2}^{n}(z^2-(m_1+1+2(m_i+n+1-i))^2)(z^2-(m_1+1-2(m_i+n+1-i))^2)\ .
\end{eqnarray*}
We see that $P_{\sigma^n}(0)=0$ iff for one index $i_0\in\{2,\dots,n\}$ the
following equation holds:
$$\frac{m_1+1}{2}=m_{i_0}+n+1-i_0\ .$$
We are going to define $\sigma^{n+k}$ by an inductive procedure. It rests on
the following claim.

Let $l\in\nat_0$ and $\mu_{\sigma^{n+l}}=(m^\prime_1,\dots,m^\prime_{n+l})$
be a highest weight of an irreducible representation of $M^{n+l}$ such that
$m_1^\prime=m_1$, $m_{i_0}^\prime=m_{i_0}$ and one of the following conditions
holds:
\begin{enumerate}
\item There exists an index $j_l\in\{i_0+l,\dots,n+l\}$ such that
$$ m_1+1-2(m^\prime_{j_l}+n+l+1-j_l)=2l\ .$$
\item $m_1+1=2l$ .
\item There exists an index $j_l\in\{i_0,\dots,n+l\}$ such that
$$ m_1+1+2(m^\prime_{j_l}+n+l+1-j_l)=2l\ .$$
\end{enumerate}
Then there exists $\sigma_{n+l+1}\in\hat M^{n+l+1}$ with the same properties
($l$ replaced by $l+1$) such that $\sigma^{n+l}\subset\sigma^{n+l+1}_{\ |M^{n+l}}$.

We now prove the claim. If $\sigma^{n+l}$ satisfies condition 1 and $m^\prime_{j_l}>
m^\prime_{j_l+1}$ (by convention $m^\prime_j:=0$ for $j>n+l$) we set
$$\mu_{\sigma^{n+l+1}}:=(m^\prime_1,\dots,m^\prime_{j_l},m^\prime_{j_l}-1, m^\prime_{j_l+1},\dots,m^\prime_{n+l})\:,\ j_{l+1}:=j_l+1\ .$$
Then $\sigma^{n+l+1}$ also satisfies condition 1.

If $\sigma^{n+l}$ satisfies condition 1 and $m^\prime_{j_l}=
m^\prime_{j_l+1}$ we set
$$\mu_{\sigma^{n+l+1}}:=(m^\prime_1,\dots,m^\prime_{j_l},m^\prime_{j_l}, m^\prime_{j_l}, m^\prime_{j_l+2},\dots,m^\prime_{n+l})\ .$$
If in addition $j_l<n+l$, set $j_{l+1}:=j_l+2$. Then again $\sigma^{n+l+1}$ satisfies condition 1. Otherwise we have $m^\prime_{n+l}=0$, hence
$m_1-1=2l$. It follows that $\sigma^{n+l+1}$ satisfies condition 2.

If $\sigma^{n+l}$ satisfies condition 2 or 3, then $\sigma^{n+l+1}$ defined by
$$\mu_{\sigma^{n+l+1}}:=(m^\prime_1,\dots,m^\prime_{n+l},0)$$
satisfies condition 3 with $j_{l+1}=n+l+1$ or $j_{l+1}=j_l$, respectively.

The branching rules for the restriction from $Sp(n+l+1)$ to $Sp(n+l)$
(see e.g. \cite{zelobenko73}, Ch. XVIII) show that in any case $\sigma^{n+l}\subset \sigma^{n+l+1}_{\ |M^{n+l}}$. The claim now follows.

Since $\sigma^n$ satisfies the induction hypothesis of the claim for $l=0$ and
$j_0=i_0$ we can define $\sigma^{n+k}$ inductively. We have to check that
$P_{\sigma^{n+k}}(-\rho^{n+k}+\rho^n)=P_{\sigma^{n+k}}(-2k\alpha)=0$. In fact, the claim ensures that $P_{\sigma^{n+k}}(z\alpha)$ containes in addition to
$z^2-(m_1+1-2(m_{i_0}+n+k+1-i_0))^2$ a second factor which contributes with
$z^2-(2k)^2$. Thus the pole originating from the $\cot$ factor cancels, and we have  $P_{\sigma^{n+k}}(-2k\alpha)=0$. Now the lemma and, hence, 
Proposition \ref{mystic} is proved.
\hB

\section{The essential spectrum}\label{esess}

In the present section we consider the spectral comparison between
$L^2(X,V(\gamma))$ and $L^2(Y,V_Y(\gamma))$.
Let $\cA$ be any commutative algebra of invariant differential operators
on $V(\gamma)$ which is generated by selfadjoint elements and contains $\cZ_\gamma$. As for $\cZ_\gamma$ there are spectral decompositions
of $L^2(X,V(\gamma))$ and $L^2(Y,V_Y(\gamma))$ with respect to $\cA$.
The main result is that the essential spectrum
of $\cA$ on $L^2(X,V(\gamma))$ and $L^2(Y,V_Y(\gamma))$
coincides.

Recall the characterization of the essential spectrum in terms of Weyl sequences. Let $\{A_i\}$ be a finite set of generators of $\cA$.
A character $\lambda$ of $\cA$ belongs to the essential
spectrum of $\cA$ iff there exists a Weyl sequence $\{\phi_\alpha\}\subset C^\infty_c$ (i.e. a sequence without accumulation points in $L^2$) such that
$\max_i \|A_i\phi_\alpha-\lambda(A_i)\phi_\alpha\|_{L^2}\rightarrow 0$ when 
$\alpha\to \infty$.

\begin{prop}\label{esspec}
The essential spectrum of $\cA$ on $L^2(X,V(\gamma))$ and on $L^2(Y,V_Y(\gamma))$
coincides.
\end{prop}
\proof 
The proof of the proposition relies on the transfer of Weyl sequences.
Let $\lambda$ be in the essential spectrum of $\cA$ on $L^2(Y,V_Y(\gamma))$.
Then there is a Weyl sequence $\{\phi_\alpha\}$, satisfying
\begin{itemize}
\item $\|\phi_\alpha\|_{L^2(Y,V_Y(\gamma))}=1$, $\forall \alpha$
\item $\{\phi_\alpha\}$ has no accumulation point in $L^2(Y,V_Y(\gamma))$, and
\item $\max_i \|  A_i\phi_\alpha-\lambda(A_i)\phi_\alpha\|_{L^2(Y,V_Y(\gamma))}\rightarrow 0$ as $\alpha\to \infty$.
\end{itemize}
Using a construction of Eichhorn \cite{eichhorn91}
we can modify this Weyl sequence such that it satisfies in addition
$\|\phi_\alpha\|_{L^2(K,V_Y(\gamma))}\rightarrow 0$ as $\alpha\to \infty$ for any compact $K\subset Y$.

Let $\{V_i\}$ be a finite number of open subsets covering $Y$ at infinity such that each
$V_i$ has a diffeomorphic lift $\tilde{V}_i\subset X$.
Using the method of Lemma \ref{lll} we can choose the $V_i$ such that there exists
a subordinated partition of unity $\{\chi_i\}$ (at infinity of $Y$) 
such that for $A\in \cA$
$$|[A,\chi_i](ka)| \le C(A) a^{-1},\quad \forall kaK\in \tilde{V}_i\ .$$
By taking a subsequence of the Weyl sequence and renumbering the $V_i$ we can assume that
$\|\chi_1 \phi_\alpha\|_{L^2(Y,V_Y(\gamma))}\ge c$ for some $c>0$ independent of $\alpha$.
We set $\psi_\alpha:=\chi_1\phi_\alpha /\|\chi_1 \phi_\alpha\|_{L^2(Y,V_Y(\gamma))}$.
We claim that $\psi_\alpha$ is again a Weyl sequence for $\lambda$.
In fact $\|\psi_\alpha\|_{L^2(Y,V_Y(\gamma))}=1$ by definition, $\| \psi_\alpha\|_{L^2(K,V_Y(\gamma))}\rightarrow 0$ as $\alpha\to \infty$ for any compact $K\subset Y$.
This implies that $\{\psi_\alpha\}$ has no accumulation points.
It remains to verify that for $A\in \cA$
$$\|(A-\lambda(A))\psi_\alpha\|_{L^2(Y)}\rightarrow 0, \quad \alpha\to \infty\ .$$
We have 
$$(A-\lambda(A))\psi_\alpha=\frac{\chi_1}{\|\chi_1 \phi_\alpha\|_{L^2(Y,V_Y(\gamma))}}(A-\lambda(A))\phi_\alpha
+\frac{[A,\chi_1]\phi_\alpha}{\|\chi_1  \phi_\alpha\|_{L^2(Y,V_Y(\gamma))}} \ .$$
Obviously we have
$$\|\frac{\chi_1}{\|\chi_1 \phi_\alpha\|_{L^2(Y,V_Y(\gamma))}}(A-\lambda(A))\phi_\alpha\|_{L^2(Y,V_Y(\gamma))}\rightarrow 0, \quad \alpha\to \infty\ .$$
For any $\epsilon>0$ we can choose $K\subset Y$ compact such that
$\sup_{x\in Y\setminus K} |[A,\chi_1](x)|<\epsilon$.
We then have
\begin{eqnarray*}
\lefteqn{\lim_{\alpha\to\infty}\|\frac{1}{\|\chi_1  \phi_\alpha\|_{L^2(Y,V_Y(\gamma))}}[A,\chi_1]\phi_\alpha\|_{L^2(Y,V_Y(\gamma))}}\hspace{1cm}\\
&\le &\lim_{\alpha\to\infty} \left( \frac{\sup_{x\in K}|[A,\chi_1](x)|}{\|\chi_1  \phi_\alpha\|_{L^2(Y,V_Y(\gamma))}} \|\phi_\alpha\|_{L^2(K,V_Y(\gamma))} +   \frac{\epsilon}{\|\chi_1  \phi_\alpha\|_{L^2(Y,V_Y(\gamma))}}
\|\phi_\alpha\|_{L^2(Y\setminus K,V_Y(\gamma))}\right) \\
&\le& \epsilon/c\ .
\end{eqnarray*}
It follows that
$$\lim_{\alpha\to\infty}\|\frac{[A,\chi_1]\phi_\alpha}{\|\chi_1  \phi_\alpha\|_{L^2(Y,V_Y(\gamma))}} \|_{L^2(Y,V_Y(\gamma))}=0\ .$$
Thus $\{\psi_\alpha\}$ is a Weyl sequence of the algebra $\cA$ 
with respect to the character $\lambda$ which is   supported in $V_1$.
Lifting this sequence to $X$ we obtain a Weyl sequence of $\cA$ to $\lambda$
in $L^2(X,V(\gamma))$.
Thus the essential spectrum of $\cA$ on $L^2(Y,V_Y(\gamma))$ is contained in the essential spectrum of $\cA$ on $L^2(X,V(\gamma))$.

In order to prove the opposite inclusion
we choose a finite cover of infinity of $X$ by sets $W_j$ such that the central projection of
$W_j$ is not surjective onto $\partial X$.
If the character $\lambda$ is in the essential spectrum of $\cA$ on $L^2(X,V(\gamma))$, then repeating the above construction we can find a Weyl sequence $\{\psi_\alpha\}$ of $\cA$ to $\lambda$ with
$\supp(\psi_\alpha)\subset W_1$. There is a $g\in G$ such that $gW_1\subset V_1$.
Then $g^*\psi_\alpha$ can be pushed down to $Y$ and gives a Weyl sequence
for $\cA$ to $\lambda$ on $L^2(Y,V_Y(\gamma))$.
Thus the essential spectrum of $\cA$ on $L^2(X,V(\gamma))$ is contained in the essential spectrum of $\cA$ on $ L^2(Y,V_Y(\gamma))$.
\hB

\section{Relevant generalized eigenfunctions}\label{relsec}

In principle a spectral decomposition of $L^2(Y,V_Y(\gamma))$ with respect to
$\cZ_\gamma$ is a way of expressing elements of $L^2(Y,V_Y(\gamma))$
in terms of generalized eigensection of $\cZ$, i.e., sections
in $C^\infty(Y,V_Y(\gamma))$ on which $\cZ$ acts by a character.
Its turns out that only a small portion of these eigenfunctions
is are needed for the spectral decomposition. We call them relevant.
Relevant eigenfunctions satisfy certain growth conditions.

In order to deal with these growth conditions properly we 
introduce the Schwartz space $S(Y,V_Y(\gamma))$.
This Schwartz space is the immediate generalization of Harish-Chandra's Schwartz space $S(X,V(\gamma))$ to our locally symmetric situation.
The significance of the Schwartz space is the following.
If a generalized eigenfunction is relevant for the spectral decomposition
of $L^2(Y,V_Y(\gamma))$ with respect to $\cZ_\gamma$ then
it defines a continuous functional on the Schwartz space. 
The  Schwartz space is defined as a Fr\'echet space which is contained in $L^2(Y,V_Y(\gamma))$. Then its (hermitian) dual $S(Y,V_Y(\gamma))^\prime$ contains the relevant generalized eigenfunctions.

We now turn to the definition of $S(Y,V_Y(\gamma))$.
The idea is to require the similar conditions as for $S(X,V(\gamma))$
locally at infinity.
Let $W\subset X\cup\Omega$ be any compact subset. For any $A\in\cU(\gaaa)$, $N\in \nat$, we define the seminorm $q_{W,A,N}(f)\in[0,\infty]$ 
of $f\in C^\infty(X,V(\gamma))$ by
\begin{equation}\label{swart}q_{W,A,N}(f):=\int_{(W\cap X)K} \log (\|g\|^N ) |f(Ag)|^2 dg\ .\end{equation}
Here $\|g\|$ denotes the norm of $\Ad(g)$ on $\gaaa$.    

In the following definition we identifiy $C^\infty(Y,V_Y(\gamma))$ with the subspace of $\Gamma$-invariant sections in $C^\infty(X,V(\gamma))$.
\begin{ddd}
The Schwartz space is the space of sections $f\in C^\infty(Y,V_Y(\gamma))$ with $q_{W,A,N}(f)<\infty$ for all $W$, $A$ and $N$ as above. 
The seminorms  $q_{W,A,N}$ define the Fr\'echet
topology of $S(Y,V_Y(\gamma))$.
\end{ddd}
 
\begin{ddd} 
An eigenvector of $\cZ_\gamma$ belonging to
the dual of the Schwartz space $S(Y,V_Y(\gamma))^\prime$ is called
tempered.
\end{ddd}

The main goal of the present section is to provide a list of all
tempered eigenvectors of $\cZ$. Unfortunately there exist exceptional 
characters where such a description is difficult. Fortunately
this exceptional set is at most countable. We will cover the set 
of all characters of $\cZ$ which may provide difficulties
by a countable set $PS$ below.  

The set of $\lambda\in\haaa_\C^*$ of parametrizing integral
characters of $\cZ$ forms a lattice and is hence countable.
We denote that set of integral characters of $\cZ$ by $PS_i$.

We choose a Cartan algebra $\taaa$ of $\maaa$ such that $\aaaa\oplus\taaa=\haaa$
and a positive root system of $\taaa$. Let $\rho_m$ denote half of the sum of the positive roots of $(\maaa,\taaa)$.
For $\sigma\in \hat{M}$ let $\mu_\sigma\in\taaa^*$ be the highest weight.
Then the character $\chi_\lambda$ of $\cZ$ on the principal series representation
$C^\infty(\partial, V(\sigma_\mu))$ is parametrized by
$\lambda:=\mu_\sigma+\rho_m-\mu\in\haaa_\C^*$.
We observe that if $\chi_\lambda\not\in PS_i$, then $\mu\in\acag$. Indeed, a
pole of $P_\sigma$ at $\mu$ implies the integrality $\chi_\lambda$ (compare
(\ref{blanche})), whereas a pole of $\hat J_{-\mu}$ for non-integral $\chi_\lambda$ cancels with a zero of $P_\sigma$ at $\mu$ because of (\ref{spex})
and the irreducibility of $C^{\infty}(\partial X,V(\sigma_{\mu}))$.

If $\chi\not\in PS_i$, then we have a complete understanding
of the corresponding eigenspace 
$$\cE_\chi:= \{f\in  C^\infty_{mg}(X,V(\gamma)), \quad (A-\chi(A))f=0 \quad \forall A\in\cZ \}\ .$$  
Let $(V_\gamma)_{|M}=\oplus_i V_\gamma(\sigma_i)$
denote the decomposition of $(V_\gamma)_{|M}$ into isotypic components.
Set $r_i=[(V_\gamma)_{|M}:V_{\sigma_i}]$ and
choose isomorphisms $T_i\in \Hom_M(\oplus_{j=1}^{r_i} V_{\sigma_i},V_\gamma(\sigma_i))$. Employing the fact that principal series representations with non-integral infinitesimal character are irreducible it is
a simple matter to deduce the following lemma from the results of \cite{olbrichdiss}. 
\begin{lem}\label{oppp} 
If $\chi\not\in PS_i$, then
the eigenspace $\cE_\chi$ is the isomorphic
image of the Poisson transform $$\bigoplus_i P_{\mu_i}^{T_i}:\bigoplus_i \bigoplus_{j=1}^{r_i}C^{-\infty}(\partial X,V(\sigma_{i,\mu_i}))\rightarrow C^\infty_{mg}(X,V(\gamma))\ ,$$
where $\mu_i\in\aca$ is uniquely characterized by $\chi_{\mu_{\sigma_i}+\rho_m-\mu_i}=\chi$, $\Ree(\mu_i)>0$ or $\Ree(\mu_i)=0$ and $\Imm(\mu_i)\ge 0$.  
\end{lem}
  
Let $PS_d$ denote the set of all characters $\chi$ of $\cZ$
such that there exists $\sigma\subset\gamma_{|M}$ and $\mu\in\aca$
with $\Ree(\mu)\ge 0$, $\chi_{\mu_\sigma+\rho_m-\mu}=\chi$
and ${}^\Gamma C^{-\infty}(\Lambda,V(\sigma_\mu))\not= 0$ or $\mu =0$.
By Lemma \ref{ghu} and the fact that $ext$ is meromorphic
on $\acag$ the set $PS_d$ is countable.
Let $PS=PS_i\cup PS_d$. Then $PS$ is countable, too.

We recall the asymptotic expansion of eigenfunctions
(see \cite{wallach88}, Ch. 4, \cite{wallach92}, Ch.11).
Since we want to avoid hyperfunction boundary values
we consider eigensections of moderate growth. A section
$f\in C^\infty(X,V(\gamma))$ is of moderate growth if
for any $A\in \cU(\gaaa)$ there exists $R\in\R$ such that
$$\sup_{g\in G} \|g\|^{-R}|f(Ag)| <\infty\ .$$
Let $C_{mg}^\infty(X,V(\gamma))$ be the space
of all sections of $V(\gamma)$ of moderate growth.

Let $L^+\subset \aaaa^*_+$ be the semigroup generated by the
positive roots of $(\aaaa,\naaa)$ and $0$.  
Let $f\in C^\infty(X,V(\gamma))$ be some $\cZ$-finite section
which is $K$-finite with respect to the action of $K$ by left translations.
Then there is
a finite set of leading exponents $E(f)\subset \aca$
such that
$$f(ka)\stackrel{a\to\infty}{\sim} \sum_{\mu\in E(f)} a^{\mu-\rho} \sum_{\Q\in L^+} a^{-Q} p(f,\mu,Q)(k)(\log(a))\ ,$$
where $p(f,\mu,Q)$ is a polynomial (required to be non-trivial for $Q=0$) on $\aaaa$ with values in $C^\infty(K,V_\gamma)$.
To read this expansion properly consider $f$ as a function on $G$ with
values in $V_\gamma$.
The leading coefficient of the polynomial $p(f,\mu,0)$  has a continuous extension with respect to $f$ which is a $G$-equivariant continuous map
from the closed $G$-submodule of $C^\infty_{mg}(X,V(\gamma))$ generated by $f$ to  $C^{-\infty}(\partial X,V(\gamma_{|M,\mu}))$.
If $f$ is an eigensection of $\cZ$ and $\mu\not=0$, then $p(f,\mu,0)$ is in fact a constant polynomial
\cite{olbrichdiss}, Lemma 4.6.

If $f\in C^\infty(Y,V_Y(\gamma))$ is a tempered generalized eigensection of $\cZ$, then its $\Gamma$-invariance implies that $f\in C_{mg}^\infty(X,V(\gamma))$.
It follows that $p(f,\mu,0)\in {}^\Gamma C^{-\infty}(\partial X,V(\gamma_{|M,\mu}))$ is a $\Gamma$-invariant distribution.
\begin{lem}\label{wo}
If $f$ is a tempered generalized eigensection of $\cZ$, then
$\supp (p(f,\mu,0))\subset\Lambda$ for all $\mu\in E(f)$ with $\Ree(\mu)>0$.
\end{lem}
\proof
We argue by contradiction.
Consider the exponent $\mu\in E(f)$ with the largest real part $\Ree(\mu)>0$
such that $\supp(p(f,\mu,0))\cap\Omega \not= \emptyset$. We assume that such an exponent $\mu$ exists. Note that $p(f,\mu,0)$ is a constant polynomial on $\aaaa$.

We study the support of of $p(f,\mu,0)$ by testing this distribution against
suitable test functions.
Let $F\subset X\cup\Omega$ be a fundamental domain of $\Gamma$
and $\partial F:=F\cap\Omega$.
Since $p(f,\mu,0)$ is $\Gamma$-invariant and since we have
freedom to choose $F$ it suffices to show that $\supp(p(f,\mu,0))\cap\interi(\partial F)=\emptyset$.

Thus let $\phi\in C_c^\infty(\interi(\partial F),V(\gamma_{|M,-\bar{\mu}}))$ be a test function. The application of the distribution boundary value $p(f,\mu,0)$
to $\phi$ can be written as the limit
$$
(\phi,p(f,\mu,0)) = \lim_{a\to \infty} a^{\rho-\mu}
\int_K ( \phi(k), f(ka) )$$
(we write sesquilinear pairings as $(.,.)$).
Let $\chi\in C^\infty_c(0,1)$ be such that $\int_0^1\chi(t) dt=1$.
We extend $\chi$ to $\R$ by zero.
Then we define $\chi_n\in C_c^\infty(A_+)$ by
$\chi_n(a):=a^{-\rho-\bar{\mu}}\chi(|\log(a)|-n)$.
If we set $\phi_n(ka)= \chi_n(a)\phi(k)$,
then we can write 
\begin{equation}\label{mnb}( \phi,p(f,\mu,0))
=\lim_{n\to \infty} (\phi_n,f) \ .\end{equation}
 
If $n$ is sufficiently large, then $\supp(\phi_n)\subset F$.
$\phi_n$ descends to a Schwartz space section $\tilde{\phi}_n\in S(Y,V_Y(\gamma))$. Using the fact that $\Ree(\mu)>0$ one can easily check that $\lim_{n\to \infty}\tilde{\phi}_n=0$ in $S(Y,V_Y(\gamma))$. 
The right-hand side of (\ref{mnb}) can be written as the application
of $f\in S(Y,V_Y(\gamma))^\prime$ to $\tilde{\phi}_n\in S(Y,V_Y(\gamma))$.
It follows that
$$( \phi,p(f,\mu,0))=\lim_{n\to 0}(\tilde{\phi}_n,f)=0\ .$$
Since $\phi$ was arbitrary we conclude that $\supp(p(f,\mu,0))\cap\interi(\partial F)=\emptyset$.
As noted above it follows that $\supp(p(f,\mu,0))\cap\Omega=\emptyset$
and this contradicts our assumption.\hB

The following proposition is a part of our description of the generalized eigenfunctions
of $\cZ$ which are relevant for the spectral decomposition.
We adopt the notation of Lemma \ref{oppp}.
\begin{prop}\label{gener1}
If $\chi\not\in PS$ and $f\in \cE_\chi\cap S(Y,V_Y(\gamma))^\prime$ is a tempered eigenfunction of $\cZ$,
then $f=\bigoplus_i P_{\mu_i}^{T_i}(\phi_i)$, where 
$\Ree(\mu_i)\ge 0$, $\chi_{\mu_{\sigma_i}+\rho_m-\mu_i}=\chi$.
Moreover, $\phi_i\not= 0$ implies that $\Ree(\mu_i)=0$.
\end{prop}
\proof
Assume that $\chi\not\in PS$. Let $0\not=f\in E_\chi$ be a tempered
eigenfunction of $\cZ$. Then by Lemma \ref{oppp} we can represent $f$ as
$\sum_i P_{\mu_i}^{T_i}(\phi_i)$. Here $\phi_i\in  \bigoplus_{j=1}^{r_i}{}^\Gamma C^{-\infty}(\partial X,V(\sigma_{i,\mu_i}))$ is uniquely characterized by $c_\gamma(\mu_i)T_i\phi_i=p(f,\mu_i,0)$. 
If $\phi_i\not=0$, then by Lemma \ref{wo} and the definition of $PS_d\subset PS$  we have $\Ree(\mu_i)=0$.\hB

\section{Wave packets and scalar products}\label{wxa}

In this section we introduce wave packets of Eisenstein series and 
show that they belong to the Schwartz space.
Thus the scalar product of a wave packet with such a tempered generalized
eigenfunction of $\cZ$ makes sense. If the corresponding character
does not belong to the exceptional set $PS$, then we obtain an explicit formula for this scalar product.

The subspace $L^2(Y,V_Y(\gamma))_c\subset L^2(Y,V_Y(\gamma))$ spanned by the wave packets is the absolute contiuous subspace. 
We show that its complement is the discrete subspace $L^2(Y,V_Y(\gamma))_d$
and that there is no singular continuous subspace.

It turns out that the support of the continuous spectrum
of $\cZ$ on $L^2(Y,V_Y(\gamma))$ concides with the support
of the continuous spectrum of $\cZ$ on $L^2(X,V(\gamma))$
which is well known by the Harish-Chandra Plancherel theorem.
The main result of the present section is Theorem \ref{contsp}.

 First we introduce the notion of an Eisenstein series.
Let $\sigma\in \hat{M}$,
$T\in\Hom_M(V_\sigma,V_\gamma)$, and  let $P^T_\mu$, $\mu\in\aca$,
be the associated Poisson transform (see Definition \ref{defofpoi}). 
\begin{ddd}
The Eisenstein series associates to $\phi\in C^{-\infty}(B,V(\sigma_\mu))$   
the eigensection of $\cZ$    
$$E(\mu,\phi,T):= P_\mu^T\circ ext(\phi)\in C^\infty(Y,V_Y(\gamma))$$
corresponding the character $\chi_{\mu_\sigma+\rho_m-\mu}$.
\end{ddd}
The meromorphic continuation of $ext$ immediately implies the following corollary.
\begin{kor} For $\mu\in\acag$ the Eisenstein series
$E(\mu,.,T):C^{-\infty}(B,V(\sigma_\mu))\rightarrow C^\infty(Y,V_Y(\gamma))$ is a meromorphic family  of operators with finite-dimensional singularities.
The Eisenstein series is holomorphic on $\{\Ree(\mu)=0,\mu\not=0\}$.
Moreover, $E(\mu,.,T)$ has at most first-order poles in the set $\{\lambda\in\acag\:|\:\Ree(\mu)>0\}$.
\end{kor}

We now discuss the functional equations satisfied by the Eisenstein series.
This functional equation will be deduced from the corresponding functional equation of the Poisson transform.

Recall the definition (\ref{cgamma}) of $c_\gamma(\lambda)$.
Next we recall the functional equation satisfied by the Poisson transform
which is proved in \cite{olbrichdiss}:
\begin{equation}\label{mi9}c_{\sigma}(\lambda) P^T_\lambda \circ J_{-\lambda} = P^{\gamma(w)c_{ \gamma}(\lambda)T}_{-\lambda }  \ ,\end{equation}
where $w\in N_K(M)$
represents the generator of the Weyl group of $(\gaaa,\aaaa)$.
If we use $J_\lambda\circ ext = ext\circ S_\lambda$, then we obtain
the following corollary.
\begin{kor}\label{funeq}
The Eisenstein series satisfies the functional equation 
$$E(\lambda,c_{\sigma}(\lambda)S_{-\lambda}\phi,T)=E(-\lambda ,\phi,\gamma(w)c_{ \gamma}(\lambda)T)\ .$$
(To be more precise, this is an identity of meromorphic quantities
valid for all $\lambda\in \aca$ where all terms are meromorphic.)
\end{kor}

Now we turn to the definition of the wave packet transform.
Roughly speaking, a wave packet of Eisenstein series is
an average of the Eisenstein series over imaginary parameters
with a smooth, compactly supported weight function.
More precisely, the space of such weight functions
$\cH_0$ is the linear space of smooth families $\aaaa^*_+\ni \mu\mapsto \phi_{\imath\mu}\in C^\infty(B,V_B(\sigma_{\imath\mu}))$ with compact support
with respect to $\mu$. 
Because of the symmetry \ref{funeq} it will be sufficient to consider wave packets on the positive imaginary axis, only.
 
Next we fix a convenient choice of the endomorphisms $T$ entering the definition of the Eisenstein series. 

Let $\tilde{\gamma}$ be the dual representation of $\gamma$ and let
$\tilde{T}\in \Hom_M(V_{\tilde{\sigma}},V_{\tilde{\gamma}})$ be such that $\tilde{T}^*T=\id$.
We set $T(\lambda):=c_\sigma(\lambda)^{-1}T$ and $\tilde{T}(\lambda):=c_{\tilde{\sigma}}(\lambda)^{-1} \tilde{T}$.
There is a conjugate linear isomorphism of $\sigma$ and $\tilde{\sigma}$ and 
hence $\bar{c}_{\sigma}(\bar{\lambda})=c_{\tilde{\sigma}}(\lambda)$.
 
Now we define the wave packet transform on $\cH_0$.
Later we will extend it by continuity to a Hilbert space closure of $\cH_0$.
\begin{ddd}
The wave packet transform  is the  map $E:\cH_0\rightarrow C^\infty(Y,V(\gamma))$
given by 
$$E(\phi) := \int_0^\infty E(\imath\mu,\phi_{\imath\mu},T(\imath\mu))\:d\mu\ .$$
The section $E(\phi)$, $\phi\in\cH_0$, is called a wave packet (of Eisenstein series). 
\end{ddd}

\begin{lem}\label{swa}
If $\phi\in \cH_0$, then $E(\phi)\in S(Y,V_Y(\gamma))$.
\end{lem}
\proof
We reduce the proof of this lemma to the case where $\Gamma$ is trivial.
In this case the assertion is well known \cite{arthur75}.
Let $W\subset X\cup\Omega$ be a compact subset such that $\Gamma W=X\cup\Omega$.
Let $\partial W:=W\cap\Omega$ and let $\chi\in C_c^\infty(\Omega)$
be a cut-off function such that $\chi_{|\partial W}=1$.
Let $\tilde{\phi}_{\imath\mu}=ext\:\phi_{\imath\mu}$ and set
$\psi_{\imath\mu}=\chi \tilde{\phi}_{\imath\mu}$,
$\xi_{\imath\mu}=(1-\chi)\tilde{\phi}_{\imath\mu}$.
Then $\aaaa^*_+\ni\mu\mapsto \psi_{\imath\mu}\in C^\infty(\partial X,V(\sigma_{\imath\mu}))$ is a smooth family with compact support.
It was shown in \cite{arthur75} that 
$$P(\psi):=\int P^T_{\imath\mu}(\psi_{\imath\mu}) d\mu\in S(X,V(\gamma))\ .$$
Since $E(\phi)=P(\psi)+P(\xi)$,
it remains to show that 
$q_{W,A,N}(P(\xi))<\infty$ for all $A\in \cU(\gaaa)$, $N\in\nat$, where $q_{W,A,N}$ is one of the seminorms (\ref{swart}) characterizing the Schwartz space.
By the $G$-equivariance of the Poisson transform we have
$$L_AP(\xi)=\int P^T_{\imath\mu}(\pi^{\sigma,\imath\mu}(A) \xi_{\imath\mu}) d\mu\ .$$
Note that $\xi_{|\partial W} =0$.
The expansion of the Poisson transform
obtained  in Lemma \ref{eee} 
leads to the following decomposition 
into a leading and a remainder term.
For  $k\in \partial WM$
\begin{eqnarray*}
P^T_{\imath\mu}(\pi^{\sigma,\imath\mu}(A) \xi_{\imath\mu})(ka)&=&
  a^{-(\imath\mu+\rho)}\gamma(w)T (\hat{J}_{\imath\mu}(\pi^{\sigma,\imath\mu}(A) \xi_{\imath\mu}))(k)\\
&&+ a^{-(\imath\mu+\rho+\alpha_1)} p_1(a,k,\pi^{\sigma,\imath\mu}(A) \xi_{\imath\mu})\ ,
\end{eqnarray*}
where $p_1(a,k,f)$ is uniformly bounded as $ka\in W$
in terms of the distribution $f$.
Thus for all $N$ we have
$$\{W\ni ka\mapsto |\log(a)|^N \int a^{-(\imath\mu+\rho+\alpha_1)} p_1(a,k,\pi^{\sigma,\imath\mu}(A) \xi_{\imath\mu}) d\mu\}\in L^2(W,V(\gamma))\ .$$ 
This shows that the remainder term leads to something satisfying the 
estimates of the Schwartz space.
We now analyse the leading term.
Since the family $\xi_{\imath\mu}$ has compact support
with respect to $\mu$ and $\hat{J}_{\imath\mu}(\pi^{\sigma,\imath\mu}(A) \xi_{\imath\mu})$ is smooth in $(\mu,k)$ by Lemma
\ref{off} we obtain for all $N\in\nat_0$
$$|\int a^{-(\imath\mu+\rho)} \gamma(w) T  (\hat{J}_{\imath\mu}(\pi^{\sigma,\imath\mu}(A) \xi_{\imath\mu}))(k) d\mu|\le C_N
(1+|\log(a)|)^{-N} a^{-\rho},\quad \forall ka\in W$$
for each $N\in\nat$.
Thus the leading terms satisfies the Schwartz space estimates, too.
This proves the lemma.
\hB

Let $\mu\not=0$, $\Ree(\mu)=0$, and $\psi\in C^{-\infty}(B,V_B(\tilde{\sigma}_{-\imath\mu}))$. Set 
$\tilde{\psi}:=E(-\imath\mu,\psi,\tilde{T}(-\imath\mu))$.
\begin{lem}\label{tttq}
$\tilde{\psi}\in S(Y,V_Y(\gamma))^\prime$.
\end{lem}
\proof
We have $ext (\psi)\in C^{-\infty}(\partial X,V(\tilde{\sigma}_{-\imath\mu}))$.
It known by \cite{arthur75} that $P^{\tilde{T}}_{-\imath\mu}(ext(\psi))\in S(X,V(\gamma))^\prime$. This implies the lemma. \hB

By Lemmas \ref{swa} and \ref{tttq}
the pairing $\langle E(\phi),\tilde{\psi}\rangle$
between the wave packet $E(\phi)$ and the generalized eigensection $\tilde{\psi}$ is well defined.
The following proposition gives an explicit formula for this pairing.
\begin{prop}\label{scalar}
We have $\langle E(\phi),\tilde{\psi}\rangle=\pi\langle \phi_{\imath\mu},\psi\rangle$.
\end{prop}
\proof
Let $\psi_n\in C^{\infty}(B,V_B(\tilde{\sigma}_{-\imath\mu}))$
be a sequence approximating the distribution $\psi$. 
Then by \cite{arthur75} and the continuity of $ext$ we have
$$E(-\imath\mu,\psi_n,\tilde{T}(-\imath\mu))\to E(-\imath\mu,\psi,\tilde{T}(-\imath\mu))$$
in $S(Y,V_Y(\gamma))^\prime$ as $n\to \infty$.
Thus the proposition is a consequence of the following the special case.
\begin{lem}
Let $\psi\in  C^{\infty}(B,V_B(\tilde{\sigma}_{-\imath\mu}))$, then
$\langle E(\phi),\tilde{\psi}\rangle=\pi\langle \phi_{\imath\mu},\psi\rangle$.
\end{lem}
\proof
Let $W\subset\Omega$ be compact. 
The following asymptotic expansions hold  uniformly for $k\in WM$  and $a\in A_+$ large:
\begin{eqnarray}
E(\imath \mu,\phi_{\imath\mu},T(\imath\mu))(ka)& =&  a^{\imath\mu -\rho}\frac{c_\gamma(\imath\mu)}{c_\sigma(\imath\mu)} T ext(\phi_{\imath\mu})(k)\nonumber\\
 &&+a^{-\imath\mu-\rho}\gamma(w)T \frac{c_\sigma(-\imath\mu)}{c_\sigma(\imath\mu)} ext(S_{\imath\mu}\phi_{ \imath\mu})(k) + O(a^{-\rho-\epsilon})\label{w1w}\\
E(-\imath\mu,\psi_{-\imath\mu},\tilde{T}(-\imath\mu))(ka) & = & a^{-\imath\mu -\rho} \frac{c_{\tilde{\gamma}}(-\imath\mu)}{ c_{\tilde{\sigma}}
(-\imath\mu)} \tilde{T} ext(\psi_{-\imath\mu})(k)\nonumber \\
&&+a^{ \imath\mu-\rho}\tilde{\gamma}(w)\tilde{T} \frac{c_{\tilde{\sigma}}(\imath\mu)}{ c_{\tilde{\sigma}}
(-\imath\mu)} ext(S_{-\imath\mu}\psi_{ -\imath\mu})(k) + O(a^{-\rho-\epsilon})\ .\nonumber
\end{eqnarray}
These expansions are immediate consequences of the asymptotic
expansion of the Poisson transform of smooth sections and
can be differentiated with respect to $a$ and differentiated and integrated with respect to $\mu$.
In order to read these formulas appropriately identify sections of $V_Y(\gamma)$
with $\Gamma$-invariant functions on $G$ with values in $V_\gamma$,
sections of $V(\sigma_{\imath\mu})$ with functions on $G$ with values in $V_\sigma$, etc., as usual.

Let $\chi$ be a cut-off function as constructed in Lemma \ref{lll} and $B_R$ the ball of
radius $R$ around the origin of $X$.
If we define $A_{\imath\mu}:=-\Omega_G+\chi_{\mu_\sigma+\rho_m-\imath\mu}(\Omega_G)$, then
$A_{\imath\mu}E(\imath \mu,\phi_{\imath\mu},T(\imath\mu))=0$ and 
$A_{-\imath\mu} E(-\imath\mu,\psi_{-\imath\mu},\tilde{T}(-\imath\mu))=0$.
 
We start with the following identity
\begin{eqnarray*}0&=& \langle\chi A_{ \imath\lambda}E(\imath\lambda,\phi_{\imath\lambda},T(\imath\lambda)), E(-\imath\mu,\psi_{-\imath\mu},\tilde{T}(-\imath\mu))\rangle_{L^2(B_R)}\\&& - \langle\chi E(\imath \lambda,\phi_{\imath\lambda},T(\imath\lambda)),   A_{-\imath\mu} E(-\imath\mu,\psi_{-\imath\mu},\tilde{T}(-\imath\mu)) \rangle_{L^2(B_R)}\ .\end{eqnarray*}
By partial integration as in the proof of Proposition \ref{green}
we obtain
\begin{eqnarray}
\lefteqn{(\mu^2-\lambda^2)\langle \chi E(\imath\lambda,\phi_{\imath\lambda},T(\imath\lambda)), E(-\imath\mu,\psi_{-\imath\mu},\tilde{T}(-\imath\mu))\rangle_{L^2(B_R)}}\hspace{3cm}\\
&=&\langle \chi \nabla_n E(\imath\lambda,\phi_{\imath\lambda},T(\imath\lambda)), E(-\imath\mu,\psi_{-\imath\mu},\tilde{T}(-\imath\mu))\rangle_{L^2(\partial B_R)}\label{w2e}\\
&&-\langle \chi   E(\imath\lambda,\phi_{\imath\lambda},T(\imath\lambda)), \nabla_n E(-\imath\mu,\psi_{-\imath\mu},\tilde{T}(-\imath\mu))\rangle_{L^2(\partial B_R)}\label{w3e}\\
&&-\langle [A_0,\chi] E(\imath\lambda,\phi_{\imath\lambda},T(\imath\lambda)), E(-\imath\mu,\psi_{-\imath\mu},\tilde{T}(-\imath\mu))\rangle_{L^2(B_R)}\ .\label{w4e}
\end{eqnarray}

We insert the asymptotic expansions (\ref{w1w}) which hold on the support of $\chi$.
We obtain with $a_R:=\ee^{R}$  
\begin{eqnarray*}
(\ref{w2e})+(\ref{w3e}) &=& \imath(\lambda+\mu) a_R^{\imath (\lambda -\mu)}\langle\chi \frac{c_\gamma(\imath\lambda)}{c_\sigma(\imath\lambda)}  T ext(\phi_{\imath\lambda}),\frac{c_{\tilde{\gamma}}(-\imath\mu)}{ c_{\tilde{\sigma}}
(-\imath\mu)} \tilde{T} ext(\psi_{-\imath\mu})\rangle \\
&&+\imath (\lambda-\mu) a_R^{\imath(\lambda+\mu)} \langle\chi \frac{c_\gamma(\imath\lambda)}{c_\sigma(\imath\lambda)} T ext(\phi_{\imath\lambda}), \tilde{\gamma}(w) \tilde{T} \frac{c_{\tilde{\sigma}}(\imath\mu) }{c_{\tilde{\sigma}}
(-\imath\mu)}  ext(S_{-\imath\mu}\psi_{-\imath\mu})\rangle\\
&&+\imath(-\lambda+\mu) a_R^{\imath(-\lambda-\mu)} \langle\chi  \gamma(w) T \frac{c_\sigma(-\imath\lambda)}{ c_\sigma(\imath\lambda) } ext(S_{\imath\lambda} \phi_{\imath\lambda}), \frac{c_{\tilde{\gamma}}(-\imath\mu)}{ c_{\tilde{\sigma}}
(-\imath\mu)}\tilde{T} ext(\psi_{-\imath\mu})\rangle\\
&&+\imath(-\lambda-\mu) a_R^{\imath(-\lambda+\mu)} \langle\chi  \frac{c_\sigma(-\imath\lambda) c_{\tilde{\sigma}}(\imath\mu)}{c_\sigma(\imath\lambda )c_{\tilde{\sigma}}(-\imath\mu)} ext(S_{\imath\lambda}\phi_{\imath\lambda}), ext(S_{-\imath\mu}\psi_{-\imath\mu})\rangle\ .\\
&&o(1)\ .
\end{eqnarray*}
The pairings on the right-hand side are defined using the canonical $K$-equivariant
 identification of the bundles $V(\sigma_\lambda)$ with $V(\sigma_0)$.
 We combine the remainder $o(1)$ and the term (\ref{w4e}) to $F(\lambda,\mu,R)$.
Then we can write
\begin{eqnarray}
\lefteqn{ \langle \chi E(\imath\lambda,\phi_{\imath\lambda},T(\imath\lambda)), E(-\imath\mu,\psi_{-\imath\mu},\tilde{T}(-\imath\mu))\rangle_{L^2(B_R)}}\hspace{2cm}\label{e34r}\\
&=& \imath  \frac{a_R^{\imath (\lambda -\mu)}}{ -\lambda+\mu } \langle\chi \frac{c_\gamma(\imath\lambda)}{c_\sigma(\imath\lambda)} T
ext(\phi_{\imath\lambda}),\frac{c_{\tilde{\gamma}}(-\imath\mu)}{ c_{\tilde{\sigma}}
(-\imath\mu)} \tilde{T} ext(\psi_{-\imath\mu})\rangle\label{sing1} \\
&&-\imath  \frac{a_R^{\imath(\lambda+\mu)}}{ \lambda+\mu} \langle\chi ext(\phi_{\imath\lambda}), T^*\frac{c^*_\gamma(\imath\mu)}{c_\sigma(\imath\mu)}\tilde{\gamma}(w) \tilde{T} \frac{c_{\tilde{\sigma}}(\imath\mu)}{ c_{\tilde{\sigma}}
(-\imath\mu) } ext(S_{-\imath\mu}\psi_{-\imath\mu})\rangle\nonumber\\
&&+\imath \frac{a_R^{\imath(-\lambda-\mu)}}{\lambda+\mu} \langle\chi \tilde{T}^*\frac{c^*_{\tilde{\gamma}}(-\imath\mu)}{ c_{\tilde{\sigma}}
(-\imath\mu)}\gamma(w) T \frac{c_\sigma(-\imath\lambda)}{ c_\sigma(\imath\lambda)  }ext(S_{\imath\lambda} \phi_{\imath\lambda}), ext(\psi_{-\imath\mu})\rangle\nonumber\\
&&+\imath \frac{a_R^{\imath(-\lambda+\mu)}}{\lambda-\mu} \langle\chi  \frac{c_\sigma(-\imath\lambda)c_{\tilde{\sigma}}(\imath\mu)}{ c_\sigma(\imath\lambda )c_{\tilde{\sigma}}(-\imath\mu)} ext(S_{\imath\lambda}\phi_{\imath\lambda}), ext(S_{-\imath\mu}\psi_{-\imath\mu})\rangle\label{sing2}\\
&&+\frac{F(\lambda,\mu,R)}{\mu^2-\lambda^2}\ .
\end{eqnarray}

Since $S_{\imath\mu}$ is unitary, 
$$\frac{c_\sigma(-\imath\lambda)c_{\tilde{\sigma}}(\imath\lambda)}{ c_\sigma(\imath\lambda )c_{\tilde{\sigma}}(-\imath\lambda)}  =1\ ,$$
and 
\begin{equation}\label{wonderid}P_{\sigma}(\imath\mu)^{-1}\id_{V_{\tilde{\gamma}}(\tilde{\sigma})}=
c_\sigma(\imath\mu)c_{\tilde{\sigma}}(-\imath\mu)\id_{V_{\tilde{\gamma}}(\tilde{\sigma})}=T^* c_\gamma(\imath\mu)^*c_{\tilde{\gamma}}(-\imath\mu)\tilde{T}\ ,\end{equation}
the singularities of the terms (\ref{sing1}) and (\ref{sing2}) at $\mu=\lambda$ cancel and
 $\frac{F(\lambda,\mu,R)}{\mu^2-\lambda^2}$
is smooth at $\mu=\lambda$. Moreover $F(\lambda,\mu,R)\to 0$ as $R\to \infty$
such that the $C^1$-norm with respect to $\lambda$ remains bounded.
By Lebesgue's theorem about dominated convergence we obtain
$$\lim_{R\to\infty} \int_0^\infty \frac{F(\lambda,\mu,R)}{\mu^2-\lambda^2} d\lambda  = 0\ .$$

By the Lemma of Riemann-Lebesgue 
\begin{eqnarray*}
\lim_{R\to\infty}\int_0^\infty \imath \frac{a_R^{\imath(-\lambda-\mu)}}{\lambda+\mu} \langle\chi \tilde{T}^*\frac{c^*_\gamma(\imath\mu)}{c_\sigma(\imath\mu)}\gamma(w) T \frac{c_\sigma(-\imath\lambda)}{ c_\sigma(\imath\lambda)} ext(S_{\imath\lambda} \phi_{\imath\lambda}), ext(\psi_{-\imath\mu})\rangle d\lambda  &=& 0 \\
\lim_{ R\to\infty}\int_0^\infty  \imath  \frac{a_R^{\imath(\lambda+\mu)}}{ \lambda+\mu} \langle\chi ext(\phi_{\imath\lambda}), T^*\frac{c^*_{\tilde{\gamma}}(-\imath\mu)}{ c_{\tilde{\sigma}}
(-\imath\mu)}\tilde{\gamma}(w) \tilde{T}\frac{c_{\tilde{\sigma}}(\imath\mu)}{ c_{\tilde{\sigma}}
(-\imath\mu)} ext(S_{-\imath\mu}\psi_{-\imath\mu})\rangle d\lambda &=&0\ .
\end{eqnarray*}

We set $s:=\lambda-\mu$. 
We regroup the remaining terms of (\ref{e34r}) to 
\begin{eqnarray*}
\lefteqn{\frac{a_R^{\imath s}-a_R^{-\imath s}}{\imath s}  \langle\chi \frac{c_\gamma(\imath\lambda)}{c_\sigma(\imath\lambda)} T ext(\phi_{\imath\lambda}),\frac{c_{\tilde{\gamma}}(-\imath\mu)}{ c_{\tilde{\sigma}}
(-\imath\mu)} \tilde{T} ext(\psi_{-\imath\mu})\rangle}\hspace{0cm}\\
&&+a_R^{-\imath s } \frac{\langle\chi \frac{c_\gamma(\imath\lambda)}{c_\sigma(\imath\lambda)} T ext(\phi_{\imath\lambda}),\frac{c_{\tilde{\gamma}}(-\imath\mu)}{ c_{\tilde{\sigma}}
(-\imath\mu)} \tilde{T} ext(\psi_{-\imath\mu})\rangle- \langle\chi \frac{c_\sigma(-\imath\lambda)c_{\tilde{\sigma}}(\imath\mu) }{c_\sigma(\imath\lambda )c_{\tilde{\sigma}}(-\imath\mu)} ext(S_{\imath\lambda}\phi_{\imath\lambda}), ext(S_{-\imath\mu}\psi_{-\imath\mu})\rangle}{\imath s}
\end{eqnarray*}
If we integrate the second term with respect to $s$ and perform
the limit $R\to\infty$, then the result vanishes by the Riemann-Lebesgue lemma.
Using the identity of distributions  
$\lim_{r\to\infty}\frac{\sin(rs)}{s }=\pi\delta_0(s)$ and (\ref{wonderid}) the first term gives
$$
 \lim_{R\to\infty} \int_{-\infty}^\infty \frac{a_R^{\imath s}-a_R^{-\imath s}}{2 \imath s}  \langle\chi \frac{c_\gamma(\imath\lambda)}{c_\sigma(\imath\lambda)} T ext(\phi_{\imath\lambda}), \frac{c_{\tilde{\gamma}}(-\imath\mu)}{ c_{\tilde{\sigma}}
(-\imath\mu)} \tilde{T} ext(\psi_{-\imath\mu})\rangle  ds
 =  \pi \langle \phi_{\imath\mu} , \psi \rangle\ .
$$

The limit as $R\to\infty$ of the integral of left-hand side of (\ref{e34r}) with respect to $\lambda$ is equal to
$\langle E(\phi),\tilde{\psi}\rangle$.
This proves the lemma and thus finishes the proof of the proposition. \hB

In order to deal with Hilbert spaces we go over to
employ sesquilinear pairings which will be denoted by $(.,.)$.
The unitary structure of $\sigma$ induces a conjugate
linear isomorphism of $\tilde{\sigma}_{-\imath\mu}$ with
$\sigma_{\imath\mu}$. Analogously the unitary structure of $\gamma$
induces a conjugate linear isomorphism of $\tilde{\gamma}$ with $\gamma$.
We choose $T$ to be unitary, then $\tilde{T}$ corresponds to $T$
under the above identifications.
Moreover, if $\psi\in C^{-\infty}(B,V_B(\tilde{\sigma}_{-\imath\mu}))$
corresponds to $\bar{\psi}\in C^{-\infty}(B,V_B(\sigma_{\imath\mu}))$,
then $\tilde{\psi}$ corresponds to $\tilde{\bar{\psi}}:= E(\imath\mu,\bar{\psi},T(\imath\mu))$ under the isomorphisms above.

Using the sesquilinear pairings on $\gamma$ and $\sigma_{\imath\mu}$ we can rewrite
the result of Proposition \ref{scalar} as
\begin{equation}\label{sscar}( E(\phi),\tilde{\bar{\psi}} )=\pi( \phi_{\imath\mu},\bar{\psi})\ .\end{equation}

We define a scalar product on $\cH_0$ by
$$( \phi,\psi):=\pi \int_0^\infty ( \phi_{\imath\mu},\psi_{\imath\mu}) d\mu$$
and let $\cH$ be the corresponding Hilbert space closure of $\cH_0$.
The following corollary is a consequence of
Proposition \ref{scalar}.
\begin{kor}\label{hermit}
The wave packet transform extends by continuity to an
isometric embedding $E:\cH\hookrightarrow L^2(Y,V_Y(\gamma))$.
\end{kor}

If $A\in\cZ$ and $\phi\in\cH_0$, then we have $AE(\phi)=E(\psi)$ with $\psi_{\imath\mu}=\chi_{\mu_\sigma+\rho_m-\imath\mu}(A)\phi_{\imath\mu}$.
Choose an orthogonal decomposition $(V_\gamma)_{|M}=\oplus_iV_{\sigma_i}$ and
let $T_i\in \Hom_M(V_{\sigma_i},V_\gamma)$ be the corresponding
unitary embeddings. Here some of the $\sigma_i$ may be equivalent.
Let $\cH(i)$ be the Hilbert space corresponding to $\sigma_i$ and
$E_i$ the corresponding wave packet transform.
It is easy to modify the proof of Proposition \ref{scalar} 
in order to show that the ranges
of the $E_i$ are pairwise orthogonal.
We define the unitary embedding $$E_\gamma:=\oplus_i E_i: \cH(\gamma):=\bigoplus_i\cH(i)\hookrightarrow L^2(Y,V_Y(\gamma))\ .$$
Then  $E_\gamma$ represents an absolute-continuous subspace $L^2(Y,V_Y(\gamma))_c\subset L^2(Y,V_Y(\gamma))$ with respect to $\cZ$.

We now prove that the orthogonal complement
of the range of $E_\gamma$ is the discrete subspace
and the corresponding characters belong to $PS$.

The abstract spectral decomposition of $L^2(Y,V_Y(\gamma))$
with respect to the commutative algebra $\cZ$ provides an unitary equivalence
$$\alpha:L^2(Y,V_Y(\gamma))\cong H:= \int_{\haaa_\C^*/W} H_\lambda \kappa(d\lambda)\ ,$$
where the Hilbert space $H_\lambda$ is a $\cZ$-module on which $\cZ$ acts by $\chi_\lambda$.
 A part of the structure of the direct integral is that $H$ 
is a space of sections $ \haaa_\C^*/W\ni\lambda\mapsto \psi_\lambda\in H_\lambda$   such that
$\clo\{\psi_\lambda | \psi\in H\} = H_\lambda$ ($\kappa$-almost everywhere), and such that the scalar
products $\haaa_\C^*/W\ni \lambda\mapsto(\psi_\lambda,\psi_\lambda^\prime)\in\C$, $\psi,\psi^\prime\in H$,
are measurable functions. 
Then scalar product on $H$ is given by 
$$(\psi,\psi^\prime)=
\int_{\haaa_\C^*/W}(\psi_\lambda,\psi_\lambda^\prime)\kappa(d\lambda)\ .$$

Fix a base $\{X_i\}$ of $\gaaa$ and let $I_N$, $N\in\nat_0$, denote
the set of all multiindices $i=(i_1,\dots,i_{\dim(\gaaa)})$, $|i|\le N$.
Let $\chi$ be the cut-off function constructed in Lemma \ref{lll}.
For $f\in S(Y,V_Y(\gamma))$ we define
$$\|f\|^2_N:=\sum_{i\in I_N}\int_G \chi(gk) |\log(a(g))|^N |f(X_i g)|^2 dg \ .$$
Note that if $f\in S(Y,V_Y(\gamma))$, then $\|f\|_N<\infty$.
By $S^N(Y,V_Y(\gamma))$ we denote the closure of the Schwartz space
$S(Y,V_Y(\gamma))$ with respect to $\|.\|_N$.
Then $S^N(Y,V_Y(\gamma))$ is a Hilbert space contained in $L^2(Y,V_Y(\gamma))$.
\begin{lem}\label{komppp}
If $N$ is sufficiently large, then
the inclusion
$$S^N(Y,V_Y(\gamma))\hookrightarrow L^2(Y,V_Y(\gamma))$$
is Hilbert-Schmidt.
\end{lem}
\proof
This follows from the results of \cite{bernstein88}.
In order to provide some details we employ the
notions "space of polynomial growth" and "comparable scale functions" 
introduced in \cite{bernstein88}.
Fix some base point $y\in Y$. 
Let the scale function 
$r:\Gamma\backslash G\rightarrow \R^+$ be given by
$r(\Gamma g)=\dist_Y(y,\Gamma g K)$, where $\dist_Y$
denotes the Riemannian distance in $Y$.
Then $\Gamma\backslash G$ is a space of polynomial growth with respect to $r$,
i.e., if $A\subset G$ is a compact neighbourhood of the identity,
then there exist constants $d\ge 0$, $C>0$ such that for any $R>0$
there exists a set of $\le C(1+R)^d$ points $x_i$ of $\Gamma\backslash G$
such that  $\{r(x)<R\}\subset \cup_i x_i A$. This follows essentially
from the fact that $G$ itself is of polynomial growth and 
that $\dist_X(x,.)$ and $\dist_Y(y,.)$ are comparable
when restricted to a compact fundamental domain 
$F\subset X\cup\Omega$ containing the lift $x$ of $y$.
Note that $g\mapsto r(\Gamma g)$ and $g\mapsto a(g)$ are comparable on this fundamental domain, too.

If we choose $N>\max\{d,\dim(\gaaa)\}$, then the assertion of the lemma
is proved by the Proposition of \cite{bernstein88}, Sec. 3.4. \hB

In the following we choose $N$ suffiently large.
It follows by a theorem of Gelfand/Kostyuchenko (see \cite{bernstein88}) that
the composition
$$S^N(Y,V_Y(\gamma))\hookrightarrow L^2(Y,V_Y(\gamma))\rightarrow \int_{\haaa_\C^*/W} H_\lambda \kappa(d\lambda)$$
is pointwise defined, i.e., there exists a collection of continuous maps
$$\alpha_\lambda:S^N(Y,V_Y(\gamma))\rightarrow H_\lambda,\quad \lambda\in\haaa^*_\C/W$$ such that
for $\phi\in S^N(Y,V_Y(\gamma))$ we have $\alpha(\phi)_\lambda=\alpha_\lambda(\phi)$.

Let $S^N(Y,V_Y(\gamma))^*$ denote the Hermitean dual of $S^N(Y,V_Y(\gamma))$.
Then we have inclusions
$$S(Y,V_Y(\gamma))\hookrightarrow S^N(Y,V_Y(\gamma))\hookrightarrow L^2(Y,V_Y(\gamma))\hookrightarrow S^N(Y,V_Y(\gamma))^*\hookrightarrow S(Y,V_Y(\gamma))^*\ .$$
By changing $\alpha_\lambda$ on a set of $\lambda$'s of measure
zero (mod $\kappa$) we
can assume that for all $\lambda\in\haaa_\C^*/W$ the map
$$\alpha_\lambda:S(Y,V_Y(\gamma))\rightarrow H_\lambda$$ 
is a morphism of $\cZ$-modules.
Let 
$$\beta_\lambda:H_\lambda\rightarrow S^N(Y,V_Y(\gamma))^*$$ denote
the adjoint of $\alpha_\lambda$. Since
$$\beta_\lambda:H_\lambda\rightarrow S(Y,V_Y(\gamma))^*$$
is a morphism of $\cZ$-modules we see that $\beta_\lambda(H_\lambda)$
consists of tempered eigensections of $\cZ$ corresponding to the character 
$\chi_\lambda$.

\begin{prop}\label{ortho}
Let $\psi\in L^2(Y,V_Y(\gamma))$ be represented by
$\alpha(\psi)$ such that $\alpha(\psi)_\lambda=0$ for all $\lambda\in PS$.
Assume further that $(E_\gamma(\phi),\psi)=0$ for all wave packets
$E_\gamma(\phi)$, $\phi=\oplus_i\phi_i\in \oplus_i\cH_0(i) $.
Then $\psi=0$.
\end{prop}
\proof
Let $\phi=\oplus\phi_i\in \oplus \cH_0(i)$. Then we can write
\begin{eqnarray*}
0&=&( E_\gamma(\phi),\psi)\\
&=&\int_{\haaa_\C^*/W}(\alpha_\lambda(E_\gamma(\phi)),\alpha(\psi)_\lambda) \kappa(d\lambda)\\
&=&\int_{\haaa_\C^*/W} (E_\gamma(\phi),\beta_\lambda\alpha(\psi)_\lambda)\kappa(d\lambda)\ .
\end{eqnarray*}
We claim that for $\lambda\not\in PS$ we can write 
$$\beta_\lambda\alpha(\psi)_\lambda=\sum_i E_\gamma(\imath\mu_i,\psi_{i,\imath\mu_i},T_i(\imath\mu_i))\ ,$$
where $\psi_{i,\imath\mu_i}\in C^{-\infty}(B,V_B(\sigma_{i,\imath\mu_i}))$, 
$\Imm(\mu_i)=0$ if $\psi_{i,\imath\mu_i}\not=0$, and $\lambda=\mu_{\sigma_i}+\rho_m-\imath\mu_i$.
In fact, by Proposition \ref{gener1} we have 
$\beta_\lambda\alpha(\psi)_\lambda=\sum_i P_{\mu_i}^{T_i}(\tilde{\psi}_{i,\imath\mu_i})$, where $\tilde{\psi}_{i,\imath\mu_i}\in {}^\Gamma C^{-\infty}(\partial X,V(\sigma_{i,\imath\mu_i}))$.
If $\Imm(\mu_i)\not=0$ and $\tilde{\psi}_{i,\imath\mu_i}\not=0$, then $\supp(\tilde{\psi}_i)\in \Lambda$ by Lemma \ref{wo}. But then $\lambda \in PS_d$ and this case was excluded.
By the functional equation of the Eisenstein series
(Corollary \ref{funeq}) and since $\aaaa^*\ni\mu_i\not=0$ 
(because of $\lambda\not\in PS$) we can assume that $\Ree(\mu_i)>0$
for all relevant $i$.
By Proposition \ref{upperbound} we have 
$\tilde{\psi}_{i,\imath\mu_i}=ext (\Psi_{i,\imath\mu_i})$
with $\Psi_{i,\imath\mu_i}=res (\tilde{\psi}_{i,\imath\mu_i})$.
Putting $\psi_{i,\imath\mu_i}=c_\sigma(\imath\mu)\Psi_{i,\imath\mu_i}$  
we obtain the claim.

We consider $\mu_i$ as a function of $\lambda$. Then using (\ref{sscar})
we obtain
\begin{eqnarray*}
0&=&\sum_i\int_{\haaa_\C^*/W\setminus PS} ( E_\gamma(\phi),E(\imath\mu_i,\psi_{i,\imath\mu_i},T_i(\imath\mu_i))) \kappa(d\lambda)\\
&=&\pi\sum_i \int_{\haaa_\C^*/W\setminus PS} (\phi_{i,\imath\mu_i},\psi_{i,\imath\mu_i}) \kappa(d\lambda)\ .
\end{eqnarray*}
Let $f_i\in C_c^\infty(\haaa_\C^*/W\setminus PS)$. Then $\{ (0,\infty) \ni \mu\rightarrow   f_i(\lambda(\mu))\phi_{i,\imath\mu}\}\in\cH_0(i)$ and thus
\begin{equation}\label{opop}0=\sum_i \int_0^\infty f_i(\mu_i) (\phi_{i,\imath\mu_i},\psi_{i,\imath\mu_i}) \kappa(d\lambda)\ .\end{equation}
We conclude that $( \phi_{i,\imath\mu_i},\psi_{i,\imath\mu_i})
=0$ for almost all $\lambda$ (mod $\kappa$).
We now trivialize the family of bundles 
by identifying $V_B(\sigma_{i,\imath\mu})$ with $V_B(\sigma_{i,0})$
in some holomorphic manner. We choose a countable dense set $\{\phi_j\}\subset C^\infty(B,V_B(\sigma_{i,0}))$.
Viewing $\mu\mapsto \psi_{i,\imath\mu_i}$ as a family of sections
in $C^{-\infty}(B,V_B(\sigma_{i,0}))$ we form
$B_j:=\{\lambda|(\phi_j,\psi_{i,\imath\mu_i})\not=0\}\subset \haaa^*_\C/W$.
Then $\kappa(B_j)=0$. Moreover let $U:=\cup_j B_j$.
Then $\kappa(U)=0$ and we have 
$(\phi_j,\psi_{i,\imath\mu_i})=0$ for all $\lambda\in \haaa_\C^*/W_+\setminus U$
and all $j$.  Thus $\psi_{i,\imath\mu_i}=0$ for $\lambda\in \haaa_\C^*/W_+\setminus U$.
Hence $\alpha(\psi)_\lambda=0$ for almost all $\lambda$ (mod $\kappa$).
Hence $\psi=0$. \hB
The following theorem is the immediate consequence of Proposition \ref{ortho}.

\begin{theorem}\label{contsp}
The wave packet transform $E_\gamma$ is
an unitary equivalence of $\cH$ with the absolute-continuous
subspace of $L^2(Y,V_Y(\gamma))$.
The orthogonal complement of the absolute-continuous
subspace is the discrete subspace $L^2(Y,V_Y(\gamma))$ 
and the corresponding eigencharacters belong to $PS$.
\end{theorem}

\section{The discrete spectrum}

In Section \ref{wxa} we obtained a complete description of the 
continuous subspace $L^2(Y,V_Y(\gamma))$ in terms of 
the wave packet transform. In the present section
we study the othogonal complement $L^2(Y,V_Y(\gamma))_d$ of the continuous subspace. The discrete subspace
decomposes further into a cuspidal, residual, and 
scattering component (see Definition \ref{t6r}). 
We show that the residual and the scattering
component are finite-dimensional and that the cuspidal
component is either trivial or infinite-dimensional.
The notions of the cuspidal and the residual components are similar
to the corresponding notions known in the finite volume
case. The appearence of the scattering component is a new phenomenon which does
not occur in the finite volume case. We
give examples where the scattering component is non-trivial.

In order to study the eigenspaces of $\cZ$ on $L^2(Y,V_Y(\gamma))$
we decompose them further with respect to the full algebra
of invariant differential operators $\cD_\gamma$ on $V(\gamma)$.
We first recall some facts concerning $\cD_\gamma$ 
(see \cite{olbrichdiss}).
The algebra $\cD_\gamma$ is in general a non-commutative 
finite extension of $\cZ_\gamma$.
The right action of $\cU(\gaaa)^K$ on $C^\infty(X,V(\gamma))$
induces a surjective homomorphism $\cU(\gaaa)^K\rightarrow \cD_\gamma$.
Hence any representation of $D_\gamma$ can be lifted to
a representation of $\cU(\gaaa)^K$. For $\sigma\in\hat{M}$ and $\lambda\in\aca$ there is a representation of $\chi_{\sigma,\lambda}$
of $\cU(\gaaa)^K$ into $\End_M(V_\gamma(\sigma))$ which descends
to $\cD_\gamma$ such that its restriction to $\cZ_\gamma$
induces $\chi_{\mu_\sigma+\rho_m-\lambda}$. 
Here $V_\gamma(\sigma)$ denote the $\sigma$-isotypic component
of $V_\gamma$.
The representation
$\chi_{\sigma,\lambda}$ is characterized by
\begin{equation}\label{nearby}D\circ P^T_\lambda=P^{\chi_{\sigma,\lambda}(D)\circ T}_\lambda,\quad T\in \Hom_M(V_\sigma,V_\gamma),\:\: D\in D_\gamma\ ,\end{equation}
and where $P^T_\lambda$ denotes the Poisson transform.
To be more precise, the representations $\chi_{\sigma,\lambda}$ depend
holomorphically on $\lambda$ and (\ref{nearby}) is an identity
between holomorphic families of maps. 

By $\cE_{\sigma,\lambda}(X,V(\gamma))$ we denote the space
of all $f\in C^\infty(X,V(\gamma))$ which under $\cD_\gamma$ generate a quotient of the representation
$\chi_{\sigma,\lambda}$.
Let $\cE_{\sigma,\lambda}(Y,V_Y(\gamma))$ be the subspace of $\cE_{\sigma,\lambda}(X,V(\gamma))$ of $\Gamma$-invariant sections.
Then $\cZ$ acts on $\cE_{\sigma,\lambda}(X,V(\gamma))$ and $\cE_{\sigma,\lambda}(Y,V_Y(\gamma))$ by the character $\chi_{\mu_\sigma+\rho_m-\lambda}$. 

Define $PS_{res}(\sigma)\subset \aca$  by
$PS_{res}(\sigma):=\{\mu \:|\:\Ree(\mu)>0,\: {}^\Gamma C^{-\infty}(\Lambda,V(\sigma_\mu))\not= 0\}$.
By Proposition \ref{upperbound}, $3.$ we have $PS_{res}(\sigma)\subset (0,\rho)$.

\begin{prop}\label{ddsp}
\begin{enumerate}
\item $PS_{res}(\sigma)$ is a finite set.
\item $\dim {}^\Gamma C^{-\infty}(\Lambda,V(\sigma_\mu)) <\infty$ for all $\mu\in PS_{res}(\sigma)$.
\item The singularities of 
$ext:C^{-\infty}(B,V_B(\sigma_\mu))\rightarrow {}^\Gamma C^{-\infty}(\partial X,V(\sigma_\mu))$ are isolated in $\{\mu\in\aca\:|\:\Ree(\mu)\ge 0\}$.
\end{enumerate}
\end{prop}
\proof
Let $L^2(X,V(\gamma))_d$ denote the discrete subspace of $L^2(X,V(\gamma))$.
\begin{lem}\label{mmm2}
There exist finitely many irreducible (hence finite-dimensional)
mutually inequivalent representations $(\chi_i,W_i)$ of $D_\gamma$
such that 
\begin{equation}\label{mmm1}
L^2(X,V(\gamma))_d\cong\bigoplus_{i}V_{\pi_i}\otimes W_i
\end{equation}
as a $G\times D_\gamma$-module. Here $V_{\pi_i}$ are representations
of the discrete series $\hat{G}_d$ of $G$. Given $\sigma\in\hat{M}$ and  $\lambda\in\aca$, $\Ree(\lambda)>0$, let $\gamma$ be a minimal $K$-type of the principal series
representation $C^\infty(\partial X,V(\sigma_\lambda))$ of $G$.
Then  $\chi_i\not\cong \chi_{\sigma,\lambda}$, $\forall i$.
\end{lem}
\proof
The Harish-Chandra Plancherel Theorem for $L^2(G)$ implies
that
$$L^2(X,V(\gamma))_d=\bigoplus_{\pi\in \hat{G}_d} V_\pi\otimes \Hom_K(V_\pi,V_\gamma)\ .$$
Since the representations $V_\pi$ are irreducible and mutually non-equivalent
the same is true for the $\cU(\gaaa)^K$-modules $\Hom_K(V_\pi,V_\gamma)$
(see \cite{wallach88}, 3.5.4.). Since there is only a finite number
of $\pi\in \hat{G}_d$ with $\Hom_K(V_\pi,V_\gamma)\not=0$, equation
(\ref{mmm1}) follows.

Let $\gamma$ now be a minimal $K$-type of
$C^\infty(\partial X,V(\sigma_\lambda))$. 
Argueing by contradiction we assume that  $\chi_{\sigma,\lambda}\cong\chi_i$
for some $i$. Then there is an embedding 
$$V_{\pi_i}\otimes W_i\hookrightarrow \cE_{\sigma,\lambda}(X,V(\gamma))\cap L^2(X,V(\gamma))=:L^2_{\sigma,\lambda}\ .$$
There exists  $f\in L^2_{\sigma,\lambda}$ such that $f(1)\not=0$.
Let $f_\gamma$ be its projection onto the (left) $K$-type $\gamma$.
Then $f_\gamma\not=0$. But $f_\gamma$ is the Poisson transform
of some $\phi\in C^\infty(\partial X,V(\sigma_\lambda))(\gamma)$
(see \cite{olbrichdiss}, Thm. 3.6, and \cite{minemura92}).
Thus $f_\gamma=P^T_\lambda(\phi)$ and asymptotically
$$f_\gamma(ka)\stackrel{a\to\infty}{\sim} c_\sigma(\lambda) T\phi(k) a^{\lambda-\rho}\ .$$ 
Since $\Ree(\lambda)>0$ we have $c_\sigma(\lambda)\not=0$. We conclude that $f_\gamma\not\in L^2$. This is a contradiction to $f_\gamma\in L^2_{\sigma,\lambda}$. \hB

For the following two lemmas let $\gamma$ be a minimal $K$-type of
$C^\infty(\partial X,V(\sigma_\lambda))$.
Then its
multiplicity is one.
By Frobenius reciprocity $[\gamma:\sigma]=1$ and  
$\chi_{\sigma,\lambda}$ is a one-dimensional representation.
\begin{lem}\label{mko1}
Let $\Ree(\lambda)>0$. If there exist $\sigma^\prime\in\hat{M}$, $\lambda^\prime\in\aca$, $\Ree(\lambda^\prime)\ge 0$, such that
$\cE_{\sigma,\lambda}(X,V(\gamma))\cap \cE_{\sigma^\prime,\lambda^\prime}(X,V(\gamma))\not= 0$, then $\Ree(\lambda^\prime)\ge\Ree(\lambda)$.
\end{lem}
\proof
As in the proof of Lemma \ref{mmm2} there exists $$0\not=f_\gamma\in \left(\cE_{\sigma^\prime,\lambda^\prime}(X,V(\gamma))\cap \cE_{\sigma^\prime,\lambda^\prime}(X,V(\gamma))\right)(\gamma)\ .$$
Again $f_\gamma=P^T_\lambda(\phi)=P^{T^\prime}_{\lambda^\prime}(\phi^\prime)$
for $\phi\in C^\infty(\partial X,V(\sigma_\lambda))(\gamma)$,
$\phi^\prime\in C^\infty(\partial X,V({\sigma}^\prime_{\lambda^\prime}))(\gamma)$.
Thus on the one hand we have 
$$f_\gamma(ka)\stackrel{a\to\infty}{\sim} c_\sigma(\lambda) T\phi(k) a^{\lambda-\rho}\ ,$$
and on the other hand for any $\epsilon>0$
$$f_\gamma(ka)\stackrel{a\to\infty}{\sim} o(a^{\lambda^\prime-\rho+\epsilon})\ .$$
This implies $\Ree(\lambda^\prime)\ge\Ree(\lambda)$.\hB
\begin{lem}\label{o9o9o}
If $\lambda\in\aaaa^*$, $\lambda>0$, then there exists a commutative algebra extension $\cA\subset\cD_\gamma$ of $\cZ_\gamma$ which is generated by selfadjoint elements such that the character $(\chi_{\sigma,\lambda})_{|\cA}$
does not belong to the spectrum of $\cA$ on $L^2(X,V(\gamma))$.
\end{lem}
\proof
Consider the finite set
$B_{\sigma,\lambda}:=\{(\sigma^\prime,\lambda^\prime)\:|\: \sigma^\prime
\subset \gamma_{|M}, \lambda^\prime\in\aca,  (\chi_{\sigma^\prime,\lambda^\prime})_{|\cZ}=(\chi_{\sigma,\lambda})_{|\cZ}\}$.
If $(\sigma^\prime,\lambda^\prime)\in B_{\sigma,\lambda}$, then
$\lambda^\prime\in\aaaa^*$ since $(\chi_{\sigma,\lambda})_{|\cZ}$
is a real character. 

We define
$B_{\sigma,\lambda}^+:=\{(\sigma^\prime,\lambda^\prime)\in B_{\sigma,\lambda}\:|\: \Ree(\lambda^\prime)=0\}$.
Note that the $\chi_i$ (the $\chi_i$ have been introduced in Lemma \ref{mmm2}) and the $\chi_{\sigma^\prime,\lambda^\prime}$, $(\sigma^\prime,\lambda^\prime)\in B_{\sigma,\lambda}^+ $ are $*$-representations of $\cD_\gamma$. While the $\chi_i$ were already 
irreducible the $\chi_{\sigma^\prime,\lambda^\prime}$, $(\sigma^\prime,\lambda^\prime)\in B_{\sigma,\lambda}^+$,
can be completely decomposed into irreducible components.
Let $\chi^+$ denote the representation of $\cD_\gamma$,
which is obtained by taking the direct sum of the $\chi_i$
and mutually inequivalent representatives of the irreducible components
of $\chi_{\sigma^\prime,\lambda^\prime}$, $(\sigma^\prime,\lambda^\prime)\in B_{\sigma,\lambda}^+$.

Set $I_{\sigma,\lambda}:=\ker\chi_{\sigma,\lambda}\subset \cD_\gamma$ and
$I^+:=\ker\chi^+\subset \cD_\gamma$.
We claim that $I^+\not\subset I_{\sigma,\lambda}$.
Argueing by contradiction we assume that $I^+\subset I_{\sigma,\lambda}$.
Let $R^+$ denote the range of 
the representation $\chi^+$.
Since $\chi^+$ is the direct sum of mutually inequivalent representations
of $\cD_\gamma$ the commutant of $R^+$ is generated by the projections
onto these components. Thus  $R^+$ is a finite
direct sum $\oplus_j \Mat(l_j,\C)$ of matrix algebras.
$R^+$ admits a character $\kappa:R^+\rightarrow \C$ such that the composition
$\cD_\gamma\rightarrow R^+\stackrel{\kappa}{\rightarrow} \C$
coincides with $\chi_{\sigma,\lambda}$.
If $l_j>1$, then the restriction of $\kappa$ to the summand $\Mat(l_j)$
vanishes. Thus the representation $\chi_{\sigma,\lambda}$
of $\cD_\gamma$ must be one of the one-dimensional components defining $\chi^+$.
By Lemma \ref{mmm2} we have $\chi_{\sigma,\lambda}\not=\chi_i$
for all $i$. Hence there exists $(\sigma^\prime,\lambda^\prime)\in B_{\sigma,\lambda}^+$ such that $\chi_{\sigma^\prime,\lambda^\prime}$
contains $\chi_{\sigma,\lambda}$ as an irreducible component. 
But then $\cE_{\sigma,\lambda}(X,V(\gamma))\subset \cE_{\sigma^\prime,\lambda^\prime}(X,V(\gamma))$.
By Lemma \ref{mko1} we conclude $\Ree(\lambda^\prime)\ge \lambda>0$.
This is in conflict with the definition of
$B_{\sigma,\lambda}^+$.

Since the ideal $I^+$
is a $*$-ideal there exists a selfadjoint $A\in \cD_\gamma$ with
$A\in I^+\setminus I_{\sigma,\lambda}$.
Let $\cA$ be the algebra generated by $A$ and $\cZ_\gamma$.

Let $f\in S(X,V(\gamma))$ be an eigenfunction of $\cA$
corresponding to $(\chi_{\sigma,\lambda})_{|\cA}$. The Harish-Chandra
Plancherel theorem for the Schwartz space $S(X,V(\gamma))$
(see \cite{arthur75})
implies that $f=f_c+f_d$, where $f_d\in L^2(X,V(\gamma))_d$
and $$f_c\in\sum_{(\sigma^\prime,\lambda^\prime)\in B_{\sigma,\lambda}^+}
\cE_{\sigma^\prime,\lambda^\prime}(X,V(\gamma))\ .$$
It follows that $$f=\frac{1}{\chi_{\sigma,\lambda}(A)} A f = \frac{1}{\chi_{\sigma,\lambda}(A)} (Af_c + A f_d)=0\ .$$
We conclude that there are no tempered eigenfunctions of $\cA$
for the character $\chi_{\sigma,\lambda}$. Thus $\chi_{\sigma,\lambda}$
is not on the spectrum of $\cA$ on $L^2(X,V(\gamma))$. 
This finishes the proof of the lemma.\hB

Now we finish the proof of Proposition \ref{ddsp}.
Fix $\sigma\in \hat{M}$. Let $\gamma$ be a minimal $K$-type
of the principal series representation associated to $\sigma$.
We choose an $M$-equivariant embedding $T:V_\sigma\rightarrow V_\gamma$.
If $\mu,\lambda\in\aaaa^+$, $\lambda>\mu>0$, and if
$\phi\in {}^\Gamma C^{-\infty}(\Lambda,V(\sigma_\lambda))$, $\psi\in
{}^\Gamma C^{-\infty}(\Lambda,V(\sigma_\mu))$, then 
$P^T_\lambda(\phi)\perp P^T_\mu(\psi)$ in $L^2(Y,V_Y(\gamma))$.
If $\Ree(\mu)>0$, then $c_\sigma(\mu)\not=0$  and the Poisson transform
$P^T_\mu$ is injective.

By Corollary \ref{nahenull} there exists $\epsilon>0$ such that
$PS_{res}(\sigma)\cap(0,\epsilon)=\emptyset$.
In the following we argue by contradiction.
Assume that 
$\oplus_{\lambda\in [\epsilon,\rho]} {}^\Gamma C^{-\infty}(\Lambda,V(\sigma_\lambda))$ is infinite-dimensional.

Since $[\epsilon,\rho]$ is compact this would imply
that there exists a sequence $\mu_i$ and $\psi_i\in {}^\Gamma C^{-\infty}(\Lambda,V(\sigma_\mu))$, such that
$\mu_i\to \mu\in [\epsilon,\rho]$ and the $\psi$
are pairwise orthonormal. Let $\chi:=\chi_{\sigma,\mu}$ and
let $\cA$ be the algebra constructed for $\chi$ in Lemma \ref{o9o9o}
such that $\chi$ does not belong to the essential
spectrum of $\cA$ on $L^2(X,V(\gamma))$. But our assumption
implies that $\chi$ belongs to the essential spectrum of $\cA$ on
$L^2(Y,V_Y(\gamma))$. In fact for any $A\in\cA$ we have
$$\lim_{i\to\infty}\|(A-\chi(A))\psi_i\|\le \lim \lim_{i\to\infty}\|(A-\chi_{\sigma,\mu_i}(A))\psi_i\|+\lim_{i\to\infty}\|(\chi_{\sigma,\lambda}(A)-\chi(A))\psi_i\|=0\ .$$
By Proposition \ref{esspec} we conclude that $\chi$ belongs to the essential spectrum of $\cA$ on $L^2(Y,V_Y(\gamma))$. But this contradicts our construction of $\cA$.

Thus $\oplus_{\lambda\in [\epsilon,\rho]} {}^\Gamma C^{-\infty}(\Lambda,V(\sigma_\lambda))$ is finite-dimensional.
This shows $1.$ and $2.$ of Proposition \ref{ddsp}.
Assertion $3.$ follows from $1.$ and Lemma \ref{lead}.\hB

We now investigate the fine structure of the discrete subspace
$L^2(Y,V_Y(\gamma))_d$ for any $\gamma\in\hat{K}$.
By definition $L^2(Y,V_Y(\gamma))_d$ is the closure
of the subspace $L^2(Y,V_Y(\gamma))_\cZ$ of $\cZ$-finite
vectors. As explained in Section \ref{relsec},
if $f\in L^2(Y,V_Y(\gamma))_\cZ$, then it has an asymptotic expansion at infinity. Since by Theorem \ref{contsp} the eigencharacters of $\cZ$ on $L^2(Y,V_Y(\gamma))_d$ belong to $PS$ they are real. 
This implies (see \cite{knapp86}, Ch.8) that the set of leading exponents 
$E(f)$ is contained in $\aaaa^*$.
To be precise with the zero exponent we distinguish two types of leading exponents $\mu=0$ which we
write as $0_0,0_1$.
We say that $f$ has the leading exponent $\mu=0_1$ ($\mu=0_0$)
if $p(f,0,0)$ is a non-constant (constant) polynomial on $\aaaa$. 

We use the leading exponents in order to define a filtration $F_*$ of $L^2(Y,V(\gamma))_{\cZ}$.
For any exponent $\mu$  we set
$$F_\mu L^2(Y,V(\gamma))_{\cZ}:=\{f\in L^2(Y,V(\gamma))_{\cZ}\:|\:\mu \ge\lambda \quad \forall \lambda \in E(f)\}\ .$$ 
Then $F_\mu L^2(Y,V(\gamma))_{\cZ}$ is the subspace of $L^2(Y,V(\gamma))_{\cZ}$
on which the boundary value map $p(.,\mu,0)$ is well-defined 
(if $\mu=0_1$, then we consider the leading coefficient $p_1(.,\mu,0)=:p(.,0_1,0)$ of $p(.,\mu,0)$). We define for any leading exponent
$$L^2(Y,V_Y(\gamma))_{\cZ}^\mu:=F_\mu  L^2(Y,V(\gamma))_{\cZ}\cap \ker(p(.,\mu,0))^\perp\ .$$
\begin{ddd}\label{t6r}
We define 
\begin{eqnarray*}
L^2(Y,V_Y(\gamma))_{res}&:=&\bigoplus_{\mu>0} L^2(Y,V_Y(\gamma))_{\cZ}^\mu\\
L^2(Y,V_Y(\gamma))_{cusp}&:=&\bigoplus_{\mu<0} L^2(Y,V_Y(\gamma))_{\cZ}^\mu\\
L^2(Y,V_Y(\gamma))_{scat}&:=&\bigoplus_{\mu=0_0,0_1} L^2(Y,V_Y(\gamma))_{\cZ}^\mu\ .
\end{eqnarray*}
\end{ddd}
We apriori have $$L^2(Y,V_Y(\gamma))_d=\overline{L^2(Y,V_Y(\gamma))_{res}\oplus L^2(Y,V_Y(\gamma))_{cusp}\oplus L^2(Y,V_Y(\gamma))_{scat}}\ ,$$
but by $1.$ of Theorem \ref{poinye} the subspaces defined above are already closed. 
We now describe these spaces in detail.
\begin{theorem}\label{poinye}
\begin{enumerate}
\item 
The spectrum of $\cZ$ on $L^2(Y,V_Y(\gamma))_d$ is finite.
In particular $L^2(Y,V_Y(\gamma))_{*}$, $*\in\{res,cusp,scat\}$,
are closed subspaces.
\item 
The space $L^2(Y,V_Y(\gamma))_{res}$ is finite-dimensional.
There is an embedding $$L^2(Y,V_Y(\gamma))_{res}\hookrightarrow \bigoplus_{\sigma\subset \gamma_{|M}} \bigoplus_{\mu \in PS(\sigma)} {}^\Gamma C^{-\infty}(\Lambda,V(\sigma_{\mu}))\otimes \Hom_M(V_\sigma,V_\gamma)\ .$$
\item 
The space $L^2(Y,V_Y(\gamma))_{cusp}$ is trivial or infinite-dimensional.
More precisely, for any discrete series representation $\pi\in\hat{G}_d$ there exists
an infinite-dimensional subspace $V_{\pi,\Gamma}\subset{}^\Gamma V_{\pi,-\infty}$
($V_{\pi,-\infty}$ denotes the distribution vector globalization) such that 
for any $\gamma\in\hat{K}$
\begin{equation}\label{y7y7}L^2(Y,V_Y(\gamma))_{cusp}\cong\bigoplus_{\pi\in\hat{G}_d} V_{\pi,\Gamma}\otimes \Hom_K(V_\pi,V_\gamma)\ .\end{equation}
\item 
There exists an embedding 
$$L^2(Y,V_Y(\gamma))_{scat}\hookrightarrow
\bigoplus_{\sigma\subset \gamma_{|M}} \bigoplus_{0_0,0_1} {}^\Gamma C^{-\infty}(\Lambda,V(\sigma_0))\otimes \Hom_M(V_\sigma,V_\gamma)\ .$$
The space $L^2(Y,V_Y(\gamma))_{scat}$ is finite-dimensional.  
\end{enumerate}
\end{theorem}
\proof
The assertion $1.$ follows from $2.$ and (\ref{y7y7}).
Assertion $2.$ follows from Lemma \ref{wo} and Proposition \ref{ddsp}.
In fact, $\oplus_{\mu>0} p(.,\mu,0)$ defines the embedding
$$L^2(Y,V_Y(\gamma))_{res}\hookrightarrow \bigoplus_{\sigma\subset \gamma_{|M}} \bigoplus_{\mu \in PS(\sigma)} {}^\Gamma C^{-\infty}(\Lambda,V(\sigma_\mu))\otimes \Hom_M(V_\sigma,V_\gamma)\ .$$

We now prove $3$.
For $\pi\in \hat{G}_d$ set $V_{\pi,\Gamma}:=\Hom_G(V_\pi^*,L^2(\Gamma\backslash G))$.
For any Banach representation of $G$ on $V$ let $V_\infty$ denote the space of smooth vectors.
We have
\begin{eqnarray*}
\Hom_G(V_\pi^*,L^2(\Gamma\backslash G))&=&\Hom_G((V_\pi)^*_\infty,L^2(\Gamma\backslash G)_\infty)\\
&\subset&\Hom_G((V_\pi)^*_\infty,C^\infty(\Gamma\backslash G))\\
&=&\Hom_\Gamma((V_\pi)^*_\infty,\C)\\
&=&{}^\Gamma V_{\pi,-\infty}\ .
\end{eqnarray*}
Let $I_\pi$ be the natural injection
$$I_\pi:V_{\pi,\Gamma}\otimes \Hom_K(V_\pi,V_\gamma)\hookrightarrow (L^2(\Gamma\backslash G)\otimes V_\gamma)^K=L^2(Y,V_Y(\gamma))\ .$$
The range of $I_\pi$ consists of eigenfunctions of $\cZ$ which are matrix coefficients
of the discrete series representation $V_\pi$. 
Matrix coefficients of discrete series representations 
are characterizwed by the fact that all their leading exponents are negative.  
We conclude that $I_\pi$ injects into $L^2(Y,V(\gamma))_{cusp}$.
Conversely, if $f\in L^2(Y,V(\gamma))_{cusp}$, then all its leading
exponents are negative and hence $f$ is a finite sum of matrix coefficients of discrete series
representations. Given $\gamma$ the set $\{\pi_i\}$ of discrete series representations of $G$ satisfying $I_\pi\not=0$
is finite.
Thus 
\begin{equation}\label{scd}L^2(Y,V(\gamma))_{cusp}=\bigoplus_i \im(I_{\pi_i})\end{equation} and the spectrum of $\cZ$ on
$L^2(Y,V(\gamma))_{cusp}$ is finite.
This finishes the proof of $1$.

To prove $3.$ it remains to show that
$V_{\pi,\Gamma}$ is infinite-dimensional. 
Let $\gamma$ be a minimal $K$-type of the discrete series representation
$V_\pi$ ( or more general a $K$-type occuring with multiplicity one in $V_\pi$).
Let $\chi$ denote the corresponding character of $\cD_\gamma$
induced by the representation of $\cU(\gaaa)^K$ on $\Hom_K(V_\pi,V_\gamma)$.
Furthermore, let $\chi_i$ denote the irreducible representations of $\cD_\gamma$
introduced in Lemma \ref{mmm2}.
Without loss of generality we can assume that $\chi=\chi_1$.
Let $\{\chi_1,\dots,\chi_r\}$ denote the subset of these representations satisfying $(\chi_i)_{|\cZ}=\chi_{|\cZ}$.
We claim that there exists an abelian extension $\cA\subset\cD_\gamma$
of $\cZ_\gamma$ which separates $\chi$ from all the characters occcuring in $(\chi_i)_{|\cA}$, $\forall i>1$.  Let $\chi^+$ denote the sum of all $\chi_i$ with $i>1$. Since the $\chi_i$ are all mutually inequivalent, the range of $\chi^+$ is a finite sum of matrix algebras $\oplus_{i>1}\Mat(l_i,\C)$.
Let $I^+\subset \cD_\gamma$, $I\subset \cD_\gamma$ denote the kernels of $\chi^+$, $\chi$. We claim that $I^+\not\subset I$. 
Assuming the contrary $R^+$ would admit an character $\kappa:R^+\rightarrow \C$
such that $\chi$ is given by $\chi:\cD_\gamma\rightarrow R^+\stackrel{\kappa}{\rightarrow} \C$. But this is impossible by the definition of $\chi^+$. Choose a selfadjoint $A\in I^+\setminus I$ and let $\cA=\cZ_\gamma[A]$.
Then $\chi_i(A)=0$ for all $i>1$ but $\chi(A)\not=0$.

Now $\chi_{|\cA}$ belongs to the essential spectrum of $\cA$ on $L^2(X,V(\gamma))$.
By Proposition \ref{esspec} the character $\chi_{|\cA}$ belongs to the essential spectrum of $\cA$ on $L^2(Y,V_Y(\gamma))$. 
The characters $\chi_{|\cZ}$ and $\chi_{|\cA}$ are separated from
the continuous spectrum on $L^2(X,V(\gamma))$.
By Theorem \ref{contsp} the character $\chi_{|\cA}$ 
is also separated from the continuous spectrum of $\cA$ on $L^2(Y,V_Y(\gamma))_c$ and from the spectrum on $L^2(Y,V_Y(\gamma))_{scat}$.
Since the discrete spectrum of $\cZ$ and of $\cA$ on $L^2(Y,V_Y(\gamma))$
is finite the eigenspace of $\cA$ in $L^2(Y,V_Y(\gamma))$ according to
$\chi_{|\cA}$ must be infinite-dimensional.
Since $L^2(Y,V_Y(\gamma))_{res}$ can only contribute an finite-dimensional
subspace to this eigenspace the eigenspace of $\cA$ to $\chi_{|\cA}$
in $L^2(Y,V_Y(\gamma))_{cusp}$ is infinite-dimensional. By (\ref{scd})
this eigenspace is just given by $I_\pi (V_{\pi,\Gamma}\otimes \Hom_K(V_\pi,V_\gamma))$. It follows that $\dim\:V_{\pi,\Gamma}=\infty$.
This finishes the proof of $3.$.
 
We now prove $4$.
There is an embedding 
$$p(.,0_0,0)\oplus p(.,0_1,0):L^2(Y,V_Y(\gamma))_{scat}\hookrightarrow
{}^\Gamma C^{-\infty}(\partial X,V(\gamma_{|M,0}))\oplus {}^\Gamma C^{-\infty}(\partial X,V(\gamma_{|M,0}))\ .$$
We prove the assertion about the support.

We show that if $f\in L^2(Y,V_Y(\gamma))_{\cZ}^{0_1}$,
then $\supp(p(f,0_1,0))\subset \Lambda$.
The argument is similar to the one used in the proof of Lemma
\ref{wo}.
Let $U\subset \bar{U}\subset \Omega$
be open.  If $\phi\in C_c^\infty(U,V(\gamma_{|M,0}))$, we have
$$(\phi,p(f,0_1,0))=\lim_{a\to\infty} a^{\rho-\mu}|\log(a)|^{-1} \int_K (\phi(k),f(ka)) dk\ .$$
Constructing the sequence $\phi_n$ with $\supp(\phi_n)\subset UMA^+K$ as in the proof of Lemma \ref{wo} but using $\chi_n(a):=|\log(a)|^{-1}a^{-\rho-\bar{\mu}}\chi(|\log(a)|-n)$
we can write
$$(\phi,p(f,0_1,0))=\lim_{n\to\infty}  (\phi_n,f)\ .$$
By construction $\phi_n\to 0$ weakly in $L^2(UMA^+K,V(\gamma))$
and $f_{|UMA^+K}\in L^2(UMA^+K,V(\gamma))$. 
We obtain $(\phi,p(f,0_1,0))=0$.
Since $U$ and $\phi$ were arbitrary,
this proves that $\supp(p(f,0_1,0))\subset \Lambda$.
The same argument works in the case $\mu=0_0$.

The following lemma (for $\lambda=0$) implies that $L^2(Y,V_Y(\gamma))_{scat}$ is finite-dimensional. For $\lambda>0$ it provides an alternative proof of Proposition \ref{ddsp}, $2$.  
\begin{lem}
If $\lambda\ge 0$, then 
${}^\Gamma C^{-\infty}(\Lambda,V)(\sigma_\lambda))$ is finite-dimensional.
\end{lem}
\proof
We first prove an apriori estimate of the order of elements
of ${}^\Gamma C^{-\infty}(\Lambda,V(\sigma_\lambda))$, i.e.
we show that there exists a $k\in\nat$ such that $f\in {}^\Gamma C^{-\infty}(\Lambda,V(\sigma_\lambda))$ implies
$f\in C^k(\partial X,V(\tilde{\sigma}_{-\lambda}))^\prime$.
The inclusion ${}^\Gamma C^{-\infty}(\Lambda,V(\sigma_\lambda))\hookrightarrow 
C^k(\partial X,V(\tilde{\sigma}_{-\lambda}))^\prime$ induces on ${}^\Gamma C^{-\infty}(\Lambda,V(\sigma_\lambda))$ the structure of a Banach space.
Since ${}^\Gamma C^{-\infty}(\Lambda,V(\sigma_\lambda))$ is a closed
subspace of the Montel space $C^{-\infty}(\partial X,V(\sigma_\lambda))$
it must be finite-dimensional. It remains to prove the apriori estimate
of the order.

If $\lambda\not=0$  or if $\lambda=0$ and the principal series
representation $\pi^{\sigma,\lambda}$ is irreducible,  
then let $\gamma$ be a minimal $K$-type of the principal series representation
$\pi^{\sigma,\lambda}$ and 
fix an unitary $T\in \Hom_M(V_\sigma,V_\gamma)$.
If $\lambda=0$ 
and $\sigma=\sigma^\prime\oplus\sigma^{\prime w}$ is reducible, then let
$\gamma$ be the sum of two copies of a minimal $K$-type of the principal
series representation $\pi^{0,\sigma^\prime}$.
In this case we let $T:=T^\prime_0\oplus T^\prime_1$, where
$T^\prime_0\in \Hom_M(V_{\sigma^\prime},V_\gamma)$ and $T^\prime_1\in \Hom_M(V_{\sigma^{w \prime}},V_\gamma)$ are unitary.
In the remaining case where $\lambda=0$, $\sigma$ is irreducible and $\pi^{\sigma,\lambda}$ splits as a sum of two irreducible representations
we let $\gamma$ be the sum of minimal $K$-types of these representations.
In this case we let $T\in \Hom_M(V_\sigma,V_\gamma)$ denote
the "diagonal" embedding of $\sigma$ into $\gamma$.
In any case let
$P:=P^T_\lambda$ denote the associated injective Poisson transform.

Consider $f\in {}^\Gamma C^{-\infty}(\Lambda,V(\sigma_\lambda))$.
The asymptotic expansion (\ref{epan}) shows that $Pf$
is bounded along $\Omega$. Using the $\Gamma$-invariance we conclude
that $Pf$ is a uniformly bounded section of $V(\gamma)$.
Let $\chi$ be the infinitesimal character of $\cZ$ on the 
principal series representation $\pi^{\sigma,\lambda}$.
Let $C^\infty_{mg}(X,V(\gamma))_\chi$ denote the corresponding eigenspace.
As a topological vector space $C^\infty_{mg}(X,V(\gamma))_\chi$ is  a direct
limit of Banach spaces 
$$C^\infty_{R}(X,V(\gamma))_\chi:=\{f\in C^\infty_{mg}(X,V(\gamma))_\chi\:|\:
\sup_{g\in G}\|g\|^{-R} |f(g)| <\infty\}\ .$$
In particular, $Pf\in C^\infty_{0}(X,V(\gamma))_\chi$.
The range of the Poisson transform $P$ is a closed $G$-submodule $\cM$
of $C^\infty_{mg}(X,V(\gamma))_\chi$. 

We claim that there is a boundary value map $\beta$
defined on $\cM$ which is continuous and inverts $P$.
Before proving the claim we finish the proof of the apriori estimate assuming the claim. On the one hand the topological vector space $\cM$ is the direct limit of the Banach spaces
$\cM_R$, where $\cM_R:=\cM\cap C^\infty_{R}(X,V(\gamma))_\chi$.
On the other hand $C^{-\infty}(\partial X, V(\sigma_\lambda))$
is the direct limit of Banach spaces $C^k(\partial X,V(\tilde{\sigma}_{-\lambda}))^\prime$.
Since $\beta$ is continuous for any $R\ge 0$ there exists a $k\in\nat$ such that
$\beta(\cM_R)\subset C^k(\partial X,V(\tilde{\sigma}_{-\lambda}))^\prime$.
Since $Pf\in \cM_0$ this yields the apriori estimate we looked for.

We now show the existence of the boundary value $\beta$. It is intimately related with the leading asymptotic coefficient $p(\phi,\lambda,0)$, $\phi\in \cM$.
Let first $\lambda>0$.
Then we can define $\beta(\phi)$, $\phi\in \cM$, by
$$( \beta(\phi),\psi) := c_\sigma(\lambda)^{-1}\lim_{a\to\infty} a^{-\lambda+\rho} \int_K (\phi(ka),T\psi(k)) dk, \quad \psi\in C^\infty(\partial X,V(\sigma_{-\lambda}))\ .$$

Let now $\lambda=0$.
If $c_\sigma(\mu)$ has a pole at $\mu=0$, then
$\pi^{\sigma,\lambda}$ is irreducible an we define
$\beta(\phi)$, $\phi\in \cM$, by
$$(\beta(\phi),\psi) := \frac{1}{2\res_{\mu=0}c_\sigma(\mu)} \lim_{a\to\infty}\log(a)^{-1} a^{\rho} \int_K (\phi(ka),T\psi(k) ) dk, \quad \psi\in C^\infty(\partial X,V(\sigma_{0})) .$$
If $c_\sigma(\mu)$ is regular at $\mu=0$ and 
$\sigma=\sigma^\prime\oplus\sigma^{\prime w}$, then we define $\beta$ by
$$(\beta(\phi),\psi) := \frac{1}{c_\sigma(0)} \lim_{a\to\infty}
a^{\rho} \int_K (\phi(ka),T\psi(k) ) dk, \quad \psi\in C^\infty(\partial X,V(\sigma_{0})) .$$ 
In the remaining case  $c_\sigma(\mu)$ is regular at $\mu=0$ 
and $\pi^{\sigma,0}$ is reducible. Let $\gamma=\gamma_1\oplus\gamma_2$ and
$t_i\in\Hom_M(V_\sigma,V_{\gamma_i})$, $i=1,2$, be such that $T=t_1\oplus t_2$.
Let $c_{\gamma_i,\sigma}(\mu)$ denote the value of $c_{\gamma_i}(\mu)$
on the range of $t_i$. Note that $c_{\gamma_i,\sigma}(0)\not=0$.
We define $\beta$ by 
$$(\beta(\phi),\psi) :=   \lim_{a\to\infty}
a^{\rho} \int_K (\phi(ka),(t_1 c_{\gamma_1,\sigma}(0)^{-1}\oplus t_2 c_{\gamma_2,\sigma}(0)^{-1}) \psi(k) ) dk, \quad \psi\in C^\infty(\partial X,V(\sigma_{0})) .$$
One can check in each case that $\beta$ inverts the Poisson transform $P$
(e.g. use the method of the proof of \cite{olbrichdiss}, Lemma 4.31).
The fact that $\beta$ is continuous follows from the globalization
theory for Harish-Chandra modules. 
\hB
This finishes the proof of $4.$ and hence of the theorem.
\hB

It is clear by Theorem \ref{contsp} 
that $L^2(Y,V_Y(\gamma))_c$ is always non-trivial.
By Theorem \ref{poinye} there are examples 
with $L^2(Y,V_Y(\gamma))_{cusp}$ non-trivial.
If $\delta_\Gamma>0$, then the Patterson-Sullivan
measure leads to a non-trivial element in $L^2(Y)_{res}$.
Thus examples with non-trivial $L^2(Y,V_Y(\gamma))_{res}$ exist.

We now give an example with $L^2(Y,V_Y(\gamma))_{scat}\not=0$.
Let $\Gamma\subset SO(1,2)$ be a cocompact Fuchsian group.
Consider $\Gamma\subset SO(1,3)$ in the standard way.
Then $Y$ is a $3$-dimensional hyperbolic manifold of the type considered
in the present paper. $Y$ has two ends, i.e., $B$ has two
components. It was shown by Mazzeo-Phillips \cite{mazzeophillips90},
Corollary 3.20, that the dimension of the space of harmonic, square-integrable one-forms is $\ge \sharp\{\mbox{ends of $Y$}\}-1=1$. The character of $\cZ$ corresponding to harmonic one-forms is the boundary of the continuous
spectrum of $\cZ$ on one-forms. Thus square-integrable harmonic one-forms
are elements of $L^2(Y,V_Y(\gamma))_{scat}$.

As indicated after the proof of Proposition \ref{upperbound} 
one can obtain vanishing results for the discrete spectrum.
One can bound the residual spectrum in terms of $\delta_\Gamma$.
For certain $K$-types (e.g. for the trivial one) one
can show that $L^2(Y,V_Y(\gamma))_{scat}$ vanishes (see the remark following the proof of Lemma \ref{th43}).

\bibliographystyle{plain}

\end{document}